\documentclass[letterpaper]{amsart}

\newtheorem{theorem}{Theorem}[section]
\newtheorem{lemma}[theorem]{Lemma}

\newtheorem{proposition}[theorem]{Proposition}
\newtheorem{corollary}[theorem]{Corollary}

\theoremstyle{definition}
\newtheorem{definition}[theorem]{Definition}
\newtheorem{example}[theorem]{Example}

\theoremstyle{remark}
\newtheorem{remark}[theorem]{Remark}

\usepackage[pdftex]{graphicx}
\usepackage{color}
\usepackage{amsmath}
\usepackage{amsfonts}
\usepackage{amssymb}
\usepackage{epsfig}
\usepackage{rotating}
\usepackage{graphics}
\newcommand{\be}{\begin{equation}}

\newcommand{\ee}{\end{equation}}

\newcommand{\al}{\alpha}

\newcommand{\gam}{\gamma}

\newcommand{\om}{\omega}

\newcommand{\nen}{{\mathcal N}}

\newcommand{\si}{\sigma}

\newcommand{\la}{\lambda}

\newcommand{\eps}{\varepsilon}

\newcommand{\dz}{\wedge}

\newcommand{\ba}{\begin{array}}

\newcommand{\ea}{\end{array}}

\newcommand{\beq}{\begin{eqnarray}}

\newcommand{\eeq}{\end{eqnarray}}

\newtheorem{lm}{lemma}

\newtheorem{thee}{theorem}

\newtheorem{proo}{proposition}

\newtheorem{co}{corollary}

\newtheorem{rem}{remark}

\newtheorem{deff}{definition}

\newcommand{\bd}{\begin{deff}}

\newcommand{\ed}{\end{deff}}

\newcommand{\bl}{\begin{lm}}

\newcommand{\el}{\end{lm}}

\newcommand{\bp}{\begin{proo}}

\newcommand{\ep}{\end{proo}}

\newcommand{\bt}{\begin{thee}}

\newcommand{\et}{\end{thee}}

\newcommand{\bc}{\begin{co}}

\newcommand{\ec}{\end{co}}

\newcommand{\brm}{\begin{rem}}

\newcommand{\erm}{\end{rem}}

\newcommand{\der}{{\rm d}}

\usepackage{t1enc}
\def\frak{\mathfrak}

\newcommand{\newc}{\newcommand}

\let\ccdot\cdot
\def\cdot{\hbox to 2.5pt{\hss$\ccdot$\hss}}

\newc{\aR}{\mbox{\boldmath{$ R$}}}
\newc{\aS}{\mbox{\boldmath{$ S$}}}
\newc{\aT}{\mbox{\boldmath{$ T$}}}
\newc{\aW}{\mbox{\boldmath{$ W$}}}

\newc{\aK}{\mbox{\boldmath{$ K$}}}
\newc{\aL}{\mbox{\boldmath{$ L$}}}
\newcommand{\bbC}{\mathbb{C}}

\usepackage{amssymb}
\usepackage{amscd}
\newcommand{\del}{D}
\newcommand{\bel}{\delta}
\newcommand{\pel}{\partial}
\newcommand{\vel}{\triangle}

\newcommand{\Rho}{{\mbox{\sf P}}}

\newcommand{\hook}{\raisebox{-0.35ex}{\makebox[0.6em][r]
{\scriptsize $-$}}\hspace{-0.15em}\raisebox{0.25ex}{\makebox[0.4em][l]{\tiny
 $|$}}}

\newcommand{\bma}{\begin{pmatrix}}
\newcommand{\ema}{\end{pmatrix}}

\newcommand{\te}{\theta}
\newcommand{\bet}{\beta}

\newcommand{\nw}{{\stackrel{\scriptscriptstyle{W}}{\nabla}}\phantom{}}
\newcommand{\nwh}{{\stackrel{\scriptscriptstyle{HW}}{\nabla}}\phantom{}}
\newcommand{\gw}{{\stackrel{\scriptscriptstyle{W}}{\Gamma}}\phantom{}}
\newcommand{\ow}{{\stackrel{\scriptscriptstyle{W}}{\Omega}}\phantom{}}
\newcommand{\rw}{{\stackrel{\scriptscriptstyle{W}}{R}}\phantom{}}
\newcommand{\row}{{\stackrel{\scriptscriptstyle{W}}{\Rho}}\phantom{}}
\newcommand{\aw}{{\stackrel{\scriptscriptstyle{W}}{A}}\phantom{}}

\newc{\obstrn}[2]{B^{#1}_{#2}}

\newcommand{\rpl}                         % +) or <+
{\mbox{$
\begin{picture}(12.7,8)(-.5,-1)
\put(0,0.2){$+$}
\put(4.2,2.8){\oval(8,8)[r]}
\end{picture}$}}

\newcommand{\lpl}                         % (+ or +>
{\mbox{$
\begin{picture}(12.7,8)(-.5,-1)
\put(2,0.2){$+$}
\put(6.2,2.8){\oval(8,8)[l]}
\end{picture}$}}

\usepackage{ifthen}

\newc{\tensor}[1]{#1}
\newc{\Mvariable}[1]{\mbox{#1}}
\newc{\down}[1]{{}_{#1}}
\newc{\up}[1]{{}^{#1}}

\newc{\JulyStrut}{\rule{0mm}{6mm}}
\newc{\midtenPan}{\mbox{\sf S}}
\newc{\midten}{\mbox{\sf T}}
\newc{\midtenEi}{\mbox{\sf U}}
\newc{\ATen}{\mbox{\sf E}}
\newc{\BTen}{\mbox{\sf F}}
\newc{\CTen}{\mbox{\sf G}}

\def\sideremark#1{\ifvmode\leavevmode\fi\vadjust{\vbox to0pt{\vss% the remark
 \hbox to 0pt{\hskip\hsize\hskip1em%                          will appear only
 \vbox{\hsize3cm\tiny\raggedright\pretolerance10000%          on the side
 \noindent #1\hfill}\hss}\vbox to8pt{\vfil}\vss}}}%
\newcommand{\edz}[1]{\sideremark{#1}}

\newcommand{\bgw}{{\textstyle \bigwedge}}

\newcommand{\Span}{\mathrm{Span}}
\newcounter{romenumi}
\newcommand{\labelromenumi}{(\roman{romenumi})}

\begin{document}
\title{Sharp version of the Goldberg-Sachs theorem}
\vskip 1.truecm
\author{A. Rod Gover} \address{Department of Mathematics, University of
  Auckland, Private Bag 92019, Auckland, New Zealand}
\email{rgover@auckland.ac.nz}
\author{C. Denson Hill} \address{Department of Mathematics, Stony
  Brook University, Stony Brook, N.Y. 11794, USA}
\email{dhill@math.sunysb.edu}  
\author{Pawe\l~ Nurowski} \address{Instytut Fizyki Teoretycznej,
Uniwersytet Warszawski, ul. Hoza 69, Warszawa, Poland, and Instytut Matematyczny PAN, ul. Sniadeckich 8, Warszawa, Poland}
\email{nurowski@fuw.edu.pl} 
%\thanks{This research was supported by the Polish grant 2 P03B 12724}

\date{\today}

\begin{abstract}
We reexamine from first principles the classical Goldberg-Sachs
theorem from General Relativity. We cast it into the form valid for
complex metrics, as well as real metrics of any signature. We obtain
the sharpest conditions on the derivatives of the curvature 
that are sufficient for the implication 
(integrability of a field of alpha
planes)$\Rightarrow$(algebraic degeneracy of the Weyl tensor). With
every integrable field of alpha planes we associate a natural
connection, in terms of which these conditions have a very simple form.

\vskip5pt\centerline{\small\textbf{MSC classification}: 83C05,
  83C60}\vskip15pt
\vskip5pt\centerline{\small\textbf{Key words}: Goldberg-Sachs theorem,
algebraically special fields, Newman-Penrose formalism}\vskip15pt
\end{abstract}
\maketitle
%*************
\tableofcontents
\newcommand{\bbS}{\mathbb{S}}
\newcommand{\bbR}{\mathbb{R}}
\newcommand{\sog}{\mathbf{SO}}
\newcommand{\slg}{\mathbf{SL}}
\newcommand{\og}{\mathbf{O}}
\newcommand{\soa}{\frak{so}}
\newcommand{\sla}{\frak{sl}}
\newcommand{\sua}{\frak{su}}
\newcommand{\dr}{\mathrm{d}}
\newcommand{\sug}{\mathbf{SU}}
\newcommand{\gat}{\tilde{\gamma}}
\newcommand{\Gat}{\tilde{\Gamma}}
\newcommand{\thet}{\tilde{\theta}}
\newcommand{\Thet}{\tilde{T}}
\newcommand{\rt}{\tilde{r}}
\newcommand{\st}{\sqrt{3}}
\newcommand{\kat}{\tilde{\kappa}}
\newcommand{\kz}{{K^{{~}^{\hskip-3.1mm\circ}}}}
\newcommand{\bv}{{\bf v}}
\newcommand{\di}{{\rm div}}
\newcommand{\curl}{{\rm curl}}
\newcommand{\cs}{(M,{\rm T}^{1,0})}
\newcommand{\tn}{{\mathcal N}}
%*************
\section{Introduction}\label{itrr}
The original Goldberg-Sachs theorem of General Relativity \cite{gs} is
a statement about Ricci flat 4-dimensional Lorentzian
manifolds. Nowadays it is often stated in the following, slightly
stronger, form:
\begin{theorem}\label{gstt}
Let $({\mathcal M},g)$ be a 4-dimensional Lorentzian manifold which satisfies
the Einstein equations $Ric(g)=\Lambda g$.
Then it locally admits a congruence of null and
shearfree geodesics if and only if its Weyl tensor is algebraically special. 
\end{theorem}
If $({\mathcal M},g)$ is conformally flat, then such a spacetime admits infinitely
many congruences of null and shearfree geodesics.

This theorem proved to be very useful in General Relativity,
especially during the `golden era' of General Relativity in
the 1960s, when the important Einstein
spacetimes, such as Kerr-Newman, were constructed.

Remarkably, years after the Lorentzian version was first stated, 
it was pointed out that the theorem has a Riemannian analog
\cite{prza}. This gives a very powerful local result
in 4-dimensional Riemannian geometry, which can be stated as follows 
\cite{optical,phd}:
\begin{theorem}\label{gstte}
Let $({\mathcal M},g)$ be a 4-dimensional Riemannian manifold which satisfies
the Einstein equations $Ric(g)=\Lambda g$.
Then it is locally a hermitian manifold if and only if its Weyl tensor is algebraically special.  
\end{theorem}
Note that the notion of a congruence of null and
shearfree geodesics, in the Lorentzian case, is replaced by the
notion of a complex surface with an orthogonal 
complex structure, in the Riemannian case. Also in this case, if
$({\mathcal M},g)$ is conformally flat, it admits infinitely many
local hermitian structures.

Theorem \ref{gstte} was in particular used by LeBrun \cite{lb} to
obtain all compact complex surfaces, which admit an Einstein metric 
that is hermitean but not K\"ahler, (see also \cite{clb,lb1}).

The only other signature which, in addition to the Lorentzian and Euclidean signatures, a four
dimensional metric may have, is the `split signature': $(+,+,-,-)$. It
is again remarkable, that the Goldberg-Sachs theorem has also its split
signature version. Here, however, the situation is more complicated and
the theorem should be split into two statements:
\begin{theorem}\label{gstts}
Let $({\mathcal M},g)$ be a 4-dimensional manifold equipped with a
split signature metric which satisfies
the Einstein equations $Ric(g)=\Lambda g$.
If in addition $({\mathcal M},g)$ is either locally a pseudohermitian manifold,
or it is locally foliated by real 2-dimensional totally
null submanifolds, then $({\mathcal M},g)$ has an algebraically special Weyl
tensor.
\end{theorem}
\begin{theorem}\label{gstts1}
Let $({\mathcal M},g)$ be a 4-dimensional manifold equipped with a
split signature metric which satisfies
the Einstein equations $Ric(g)=\Lambda g$ and which is conformally non flat.
If in addition $({\mathcal M},g)$ has an algebraically special Weyl
tensor with a multiple principal totally null field of 2-planes
having locally constant real index, then it is either locally a pseudohermitian manifold,
or it is locally foliated by real 2-dimensional totally
null submanifolds.
\end{theorem}

In these two theorems the term `pseudohermitian manifold' means: `a complex
manifold with a complex structure which is an orthogonal transformation
for the split signature metric $g$'. The more complicated terms such as `multiple principal totally null field of 2-planes
having locally constant real index' will be explained in Section \ref{se2}.

All four theorems have in common the part concerned with
the Einstein assumption and algebraic speciality of the Weyl
tensor. But they look quite different on the other side of the equivalence. The similarity in the
first part suggests that also the second part should have a unified
description. This is indeed the case. As will be shown in 
the sequel, 
these theorems are consequences, or better said, 
appropriate interpretations, of the following complex theorem \cite{plhac,plro}:
\begin{theorem}\label{cgsc}
Let $({\mathcal M},g)$ be a 4-dimensional manifold equipped with a
complex valued metric $g$ which is Einstein. Then the following two conditions are equivalent:
\begin{itemize}
\item[(i)]
$({\mathcal M},g)$ 
admits a complex two-dimensional totally null distribution
$\nen\subset {\rm T}^\bbC\mathcal M$, which is
integrable in the sense that $[\nen,\nen]\subset\nen$.
\item[(ii)] The Weyl tensor of $({\mathcal M},g)$ is algebraically special.
\end{itemize} 
\end{theorem}  
\section{Convenient sharper versions}
Our motivation for reexamining these theorems is as follows:

First, as remarked e.g. by Trautman \cite{at}, all the theorems have
an \emph{aesthetic}
defect. This is due to the fact that both equivalence conditions, such
as (i) and (ii) in Theorem \ref{cgsc}, are \emph{conformal} properties of $({\mathcal M},g)$; 
the Einstein assumption does not share this symmetry. 
Of course, a way out is to replace the Einstein assumption  
by an assumption about $({\mathcal M},g)$ being
\emph{conformal} to Einstein, see e.g. \cite{govnur}. Thus, in the 
complex version of the theorem the assumption
should be: $({\mathcal M},g)$ is conformal to Einstein. 

This leads to the
question about the weakest conformal assumption involving (the derivatives
of) the Ricci part of
the curvature that is sufficient to ensure the
thesis of the Goldberg-Sachs theorem. Several authors have proposed their
assumptions here (see \cite{KT,pr,plprz1,plprz2,robsch}). 
For example the authors of
\cite{KT,pr,robsch} use
an assumption, which involves   
contractions of (the derivatives of) the Ricci tensor with the vectors
spanning the totally null distribution $\nen$. 

Trautman in \cite{at} 
has a different point of view. He proposes that there should be a conformally
invariant assumption which does not refer to the thesis of the
theorem. Trautman conjectures that a
proper replacement for the assumption is: $({\mathcal M},g)$ is 
\emph{Bach flat}. This, in four dimension, is certainly conformal,
does not refer to $\nen$, and is neccessary for $g$ to be conformal to
Einstein. 

In this paper, among other things, we show that the approach of
\cite{KT,pr,robsch} is the proper one. In particular in Section
\ref{cout} we show that, in the case of a Riemannian signature 
metric, Trautman's conjecture is not true. 

Our new analysis of the Goldberg-Sachs theorem starts with Theorem
\ref{gst}. Its proof shows that it is rather hard to find a single
curvature condition, different than the conformally Einstein one, which
would guarantee equivalence in the thesis of Goldberg and
Sachs. This proof also clearly shows that it is the implication $(algebraical~
speciality)\Rightarrow(integrability~ of~ totally~ null~ 2-planes)$ that
causes the difficulties. Then in Section \ref{lpo} we give various
generalizations of the Goldberg-Sachs theorem to the conformal
setting, starting with the conformal replacement of the assumption of
Theorem \ref{gst} which implies $(algebraical~
speciality)\Rightarrow(integra$- $bility~ of~$ $totally ~null~ 2-planes)$. This
culminates in a slight improvement of the theorem of Penrose and
Rindler \cite{pr}, which we give in our Theorem \ref{pol}, and in
Theorems 
\ref{pol1} and \ref{pol2}, which treat more special cases. These three
theorems we consider as the sharpest conformal improvement of the classical
Goldberg-Sachs theorem, in a sense that they include both implications 
$(algebraical~
speciality)\Rightarrow(integrability~ of~ totally ~null ~2-planes)$ and $(algebraical~
speciality)\Leftarrow(integrabi$- $lity~ of~ totally~ null~ 2-planes)$. In Section \ref{rc} the
real versions of theorems from Section \ref{lpo} are considered, the
most striking of them being:
\begin{theorem}\label{polo}
Let $\mathcal M$ be a 4-dimensional oriented manifold with a (real) metric $g$
of Riemannian signature, whose selfdual part of the Weyl tensor is nonvanishing. Let $J$ be a metric compatible almost complex
structure on $\mathcal M$ such that its holomorphic distribution 
${\mathcal N}={\rm T}^{(1,0)}{\mathcal M}$ is selfdual. 
Then any two of the following imply the third:
\begin{itemize}
\item[(0)] The Cotton tensor of $g$ is degenerate on ${\mathcal N}$,
  $A_{|\mathcal N}\equiv 0$.    
\item[(i)] $J$ has vanishing Nijenhuis tensor on $\mathcal M$, meaning
  that $({\mathcal M},g,J)$ is a hermitean manifold. 
\item[(ii)] The selfdual part of the Weyl tensor is algebraically special on
  $\mathcal M$ with $\mathcal N$ as a
  field of multiple principal selfdual totally null 2-planes.
\end{itemize}
\end{theorem}   
This theorem in its (more complicated) Lorentzian version is present in \cite{KT,pr,robsch}. The Riemannian version is implicit there, once one understands the relation between fields of totally null 2-planes and almost hermitian structures, as for example, explained in \cite{phd,optical}, (see also \cite{apo} where 
these developments are related to global issues on compact Riemannian
manifolds.)

When one is only interested in the implication $(algebraical~
speciality)\Rightarrow(integ$- $rability~ of~ totally ~null ~2-planes)$, 
our proposal for the sharpest version of the Goldberg-Sachs theorem,
 is given in Theorem \ref{gsccc}. This gets its final and very elegant
 (but equivalent) version in
Theorem \ref{ghj}. This last theorem utilizes a new object which we introduce in
this paper, namely a connection, which is naturally associated with
each \emph{integrable} field of totally null 2-planes $\nen$. We call this
connection the \emph{characteristic connection} of a field of totally null
  2-planes. 

If $\nen$ satisfies the integrability conditions
$[\nen,\nen]\subset\nen$, we prove in Theorem \ref{cch1} the existence
of a class of connections $\nw$, which are characterized by the following two
conditions:
$$
\begin{aligned}
& \nw_X\nen\subset \nen\\ 
&\nw_Xg=-B(X) g
\end{aligned}\quad\quad
{\rm for ~all}~ X\in {\rm T}\nen.
$$
These connections are \emph{not} canonical - they define the 1-form $B$ only
partially. However, they naturally restrict to a \emph{unique}
(partial) connection
$\check{\nabla}$ on $\nen$. This by definition is the characteristic
connection of $\nen$. In general this connection is complex. It is defined
everywhere on $\mathcal M$, but it only enables one to differentiate
vectors from $\nen$ along vectors from $\nen$. Thus the connection
$\check{\nabla}$ is effectively 2-dimensional, and as such, its
curvature $\check{R}^A_{~BCD}$ has only
one independent component. It follows that 
$$\check{R}^A_{~BCD}=4\Psi_1\delta^A_{~B}\epsilon_{CD},$$
where $\Psi_1$ is the Weyl tensor component whose nonvanishing is the
obstruction to the algebraic speciality of the metric. The symbol
$\delta^A_{~B}$ is the Kronecker delta (i.e. the identity) on $\nen$
and the $\epsilon_{CD}$ is the 2-dimensional antisymmetric tensor. The Ricci tensor 
$\check{R}_{AB}=\check{R}^C_{~ACB}$ for $\check{\nabla}$ is then 
$\check{R}_{AB}=4\Psi_1\epsilon_{AB}$ and is \emph{antisymmetric}.

Now
the replacement for the Einstein condition in the Goldberg-Sachs
theorem, in its $(integrability ~of ~\nen)\Rightarrow (algebraical
~speciality)$ part, is 
$$\check{\nabla}_{[A}\check{\nabla}_{B]}\check{R}_{CD}\equiv 0,$$
as is explained in Theorem \ref{ghj}.

An interesting situation occurs in the Riemannian (and also in the
split signature) case. There, the reality conditions imposed on the
1-form $B$ defining the class of connections $\nw$, choose a prefered
connection from the class. This connection yields more information
than the partial connection. Using this connection we get Theorem
\ref{scsi}, which is a slightly more elegant (pseudo)hermitian version
of the signature independent Theorem \ref{ghj}.

% Instead we concentrate on
%obtaining the sharpest statement in 
%\emph{only one direction} in the equivalence ${\rm (i)}\Leftrightarrow
%{\rm (ii)}$, namely ${\rm (i)}\Rightarrow
%{\rm (ii)}$. We do it in three steps. First we weaken the
%Einstein assumption of Theorem \ref{cgsc} to an 
%assumption about the vanishining of a fewer number of components of the Ricci
%tensor. This is summarized in Theorem \ref{gscc}. Then we use Theorem
%\ref{gscc} and replace the conditions 
%on the Ricci tensor with a requirement that the Cotton
%tensor of the metric has a certain conformal property. This leads to 
%Theorem \ref{gsccd}. Then the assumptions
%of Theorem \ref{gsccd} are still weakened to the ultimate stage, which
%can not be made any better, and which gives the weakest curvature
%assumptions that imply  ${\rm (i)}\Rightarrow
%{\rm (ii)}$. This is given in Theorem \ref{west} and its obvious real
%versions.

\section{Totally null 2-planes in four dimensions}\label{se2}
To discuss the geometrical meaning of the complex version of the
Goldberg-Sachs theorem we recall the known \cite{koptra} properties of 
totally null 2-planes as we range over the possible signatures of 4-dimensional
metrics.   

Let $V$ be a 4-dimensional \emph{real} vector space equipped with a
metric $g$, of some signature. Given $V$ and $g$ we consider their
complexifications. Thus we have $V^\bbC$ and the metric $g$ which is
extended to act on complexified vectors of the form $v_1+i v_2$,
$v_1,v_2\in V$, via:
$g(v_1+iv_2,v'_1+iv'_2)=g(v_1,v'_1)-g(v_2,v'_2)+i(g(v_1,v'_2)+g(v_2,v'_1))$.

Let ${\mathcal N}$ be a 2-complex-dimensional vector subspace in $V^\bbC$, 
${\mathcal N}\subset V^\bbC$, with the property that $g$ identically
vanishes on $\nen$, $g_{|\nen}\equiv 0$. In other words: $\nen$ is a 2-complex-dimensional vector
subspace of $V^\bbC$ such that for all $n_1$ and $n_2$ from $V^\bbC$
we have $g(n_1,n_2)=0$. This is the definition of $\nen$ being \emph{totally null}.

Such $\nen$s exist irrespectively of the signature of $g$. In fact, let 
$(e_1,e_2,e_3,e_4)$ be an orthonormal basis for $g$ in $V$. Then, if the
metric has signature $(+,+,+,+)$, an example of $\nen$ is given by 
$$\nen_E=\Span_\bbC(e_1+i e_2, e_3+ie_4).$$
If the metric has Lorentzian signature $(+,+,+,-)$ then
we chose the basis so that
$g(e_1,e_1)=g(e_2,e_2)=g(e_3,e_3)=1=-g(e_4,e_4)$, and as an example of
$\nen$ we take 
$$\nen_L=\Span_\bbC(e_1+i e_2,e_3+e_4).$$
In the case of split signature $(+,+,-,-)$ we have
$g(e_1,e_1)=g(e_2,e_2)=1$, $g(e_3,e_3)=g(e_4,e_4)=-1$, and we
distinguish two different classes of 2-dimensional totally null
$\nen$s. As an example of the first class we take 
$$\nen_{S_c}=\Span_\bbC(e_1+ie_2,e_3+ie_4),$$
and as an example of the second class we take 
$$\nen_{S_r}=\Span_\bbC(e_1+e_3,e_2+e_4).$$

If $V$ is a \emph{complex} 4-dimensional vector space with
a \emph{complex} metric $g$, the notion of a totally null
2-dimensional vector subspace $\mathcal N$ still makes sense: these
are simply 
2-dimensional complex vector subspaces $\nen\subset V$ for which $g_{|\mathcal
  N}\equiv 0$. 

Irrespective of the fact if the 2-dimensional totally null vector
space $\nen$ is defined in terms of a complex vector space $V$ with a
complex metric, or in terms of $(V^\bbC,g)$ in which $V$ is real
and $g$ is the complexified real metric $g$, choosing an 
orientation in $V$, one can check that $\nen$ is always either
\emph{selfdual} or \emph{antiselfdual} (see e.g. \cite{twistor}). By
this we mean that we always have 
\begin{itemize}
\item either: $*(n_1\dz n_2)=n_1\dz n_2$ for all $n_1,n_2\in\nen$,
\item or: $*(n_1\dz n_2)=-n_1\dz n_2$ for all $n_1,n_2\in\nen$,
\end{itemize} 
where $*$ denotes the Hodge star operator. Thus the property of being
selfdual or antiselfdual (partially) characterizes totally null
2-planes.

In case of \emph{real} $V$, irrespective of the metric signature,
totally null spaces in $V^\bbC$
may be further characterized by their \emph{real index} \cite{koptra}. This is
defined as follows:

Given a vector subspace $\nen\subset V^\bbC$ one considers its
\emph{complex conjugate} 
$$\bar{\nen}=\{ w\in V^\bbC ~|~\bar{w}\in\nen\}.$$
Then the intersection $\nen\cap\bar{\nen}$ is the complexification of a
real vector space, say ${\mathcal K}$, and the real index of $\nen$ is
by definition the real dimension of ${\mathcal K}$, or 
the complex dimension of $\nen\cap\bar{\nen}$, which is the
same. 

In our examples above, $\nen_E$ and $\nen_{S_c}$ have real index
\emph{zero}, $\nen_L$ has real index \emph{one} and $\nen_{S_r}$ has
real index \emph{two}. These are examples of a general fact,
discussed in any dimension in \cite{koptra}, which when specialized to
a \emph{four} dimensional $V$, reads: 
\begin{itemize}
\item[-] If $g$ has Euclidean signature, $(+,+,+,+)$, then 
every 2-dimensional totally null space $\nen$ in the
complexification $V^\bbC$ has real index \emph{zero};  
\item[-] If $g$ has Lorentzian signature, $(+,+,+,-)$, then 
every 2-dimensional totally null space $\nen$ in the
complexification $V^\bbC$ has real index \emph{one};
\item[-] If $g$ has split signature, $(+,+,-,-)$, then 
a 2-dimensional totally null space $\nen$ in the
complexification $V^\bbC$ has either real index \emph{zero} or
\emph{two};
\item[-] In either signature the spaces of all $\nen$s with indices
  \emph{zero} or \emph{one} are generic - they
  form real 2-dimensional manifolds; In the split signature the spaces of
  all $\nen$s
  with index \emph{two} are special - they form a real manifold of
  dimension \emph{one}.          
\end{itemize}  

If we have a 2-dimensional totally null $\nen$ with real index
\emph{zero} then $V^\bbC=\nen\oplus\bar{\nen}$. This enables us to
equip the real vector space $V$ with a complex structure $J$, by
declaring that the holomorphic vector space $V^{(1,0)}$ of this
complex structure is $\nen$. In other words, $J$ is defined as a linear
operator in $V$ such that, after complexification,
$J(\nen)=i\nen$. Due to the fact that $\nen$ is totally null, the so
defined $J$ is hermitian, $g(Jv_1,Jv_2)=g(v_1,v_2)$ for all $v_1,v_2 \in V$. Thus \emph{a
  totally null $\nen$ of real index zero in dimension four defines a
  hermitian structure $J$ in the corresponding 4-dimensional real
  vector space $(V,g)$}. Also the converse is true. For if we have
$(V,g,J)$ in real dimension four, we define $\nen$ by $\nen=V^{(1,0)}$,
i.e. we declare that $\nen$ is just the holomorphic vector space for
$J$. Due to the fact that $J$ is hermitian, and because of the assumed
Euclidean or split signature of the metric, $\nen$ is totally null and has
real index zero. This proves the following
\begin{proposition}\label{eus}
There is a one to one correspondence between (pseudo)hermitian structures $J$
in a four dimensional real vector space $(V,g)$, equipped with a metric of
either Euclidean or split signature, and 2-dimensional totally
null planes $\nen\subset V^\bbC$ with real index zero.
\end{proposition}

In the Lorentzian case, where all $\nen$s have index one, every $\nen$
defines a 1-real-dimensional vector space ${\mathcal K}$. This is
spanned by a real vector, say $k$, which is \emph{null}, as it is a vector
from $\nen$. 
The space  ${\mathcal
  K}^\perp$ orthogonal to ${\mathcal K}$ includes ${\mathcal K}$, ${\mathcal K}\subset{\mathcal
  K}^\perp$. Its complexification $({\mathcal
  K}^\perp)^\bbC=\nen+\bar{\nen}$. The 
\emph{quotient space} ${\mathcal H}={\mathcal K}^\perp/{\mathcal K}$
has real dimension two, and acquires a complex structure in a
similar way as $V$ did in the Euclidean/split case. Indeed, we define
$J$ in $\mathcal H$ by declaring that its 
holomorphic space ${\mathcal H}^{(1,0)}$ coincides with the
2-dimensional complex vector space
$(\nen+\bar{\nen})/(\nen\cap\bar{\nen})$. 
This shows that \emph{a 2-dimensional totally null $\nen$, in the
complexification of a Lorentzian 4-dimensional $(V,g)$, defines a real null direction $k$
in $V$ together with a complex structure $J$ in the quotient space
${\mathcal K}^\perp/{\mathcal K}$, ${\mathcal K}=\bbR k$}. One can easily see that also the converse is true, and
we have the following 
\begin{proposition}
There is a one to one correspondence between 2-dimensional totally
null planes $\nen$, in the complexification of a four
dimensional oriented and time oriented Lorentzian vector space $(V,g)$, and null directions ${\mathcal K}=\bbR k$ in
$V$ together with their associated complex structures $J$ in ${\mathcal K}^\perp/{\mathcal K}$.
\end{proposition}

The last case, in which the signature of $g$ is split, $(+,+,-,-)$, and
in which the $\nen$s have real index 2, provides us with a
\emph{real} 2-dimensional totally null plane in $V$. Thus we have
\begin{proposition}
There is a one to one correspondence between 2-dimensional totally
null planes $\nen$ with real index two, in the complexification of a four
dimensional split signature vector space $(V,g)$, and real totally null 2-planes in $V$.
\end{proposition} 

We now pass to the analogous considerations on 4-manifolds. 
Thus we consider a 4-dimensional manifold $\mathcal M$, with a metric
$g$, equipped in addition with a smooth distribution $\nen$ of complex
totally null 2-planes $\nen_x$, $x\in\mathcal M$, of a \emph{fixed} 
index. Applying the above propositions we see that, depending on the index of $\nen$, such
an $\mathcal M$ is equipped either with an \emph{almost hermitian structure}
$({\mathcal M},g,J)$ (in case of index 0), or with an \emph{almost optical
structure} $({\mathcal M},g,{\mathcal K},J_{{\mathcal
    K}^\perp/{\mathcal K}})$ (in case of index 1), or
with a \emph{real distribution of totally null 2-planes} (in case of
index 2). The interesting question about the integrability conditions
for these three different real structures has a uniform answer in
terms of the integrability of the complex distribution
$\nen$. Actually, by inspection of the three cases determined by the real
indices of $\nen$, one proves the following 
\cite{optical}
\begin{proposition}
Let $M$ be a 4-dimensional real manifold and $g$ be a real metric on
it. Let $\nen$ be a complex 2-dimensional distribution on $\mathcal M$
such that $g_{|\nen}\equiv 0$. Then the integrability condition,
$$[\nen,\nen]\subset \nen,$$    
for the distribution $\nen$ is equivalent to
\begin{itemize}
\item[-] the Newlander-Nirenberg integrability condition for the corresponding $J$,
  if $\nen$ has index zero;
\item[-] the geodesic and shear-free condition for the corresponding
  real null direction field $k$, if $\nen$ has index one. In this case
  the 3-dimensional space of integral curves of $k$ has (locally) the 
  structure of 3-dimensional CR manifold.
\item[-] the classical Fr\"obenius integrability for the real distribution 
corresponding to $\nen$, if $\nen$ has index two. In this case we have
a foliation of $\mathcal M$ by 2-dimensional real manifolds
corresponding to the leaves of $\nen$.
\end{itemize}
\end{proposition}   

Returning to the complex Goldberg-Sachs theorem \ref{cgsc}, we see that one part of its
thesis, which is concerned with the integrabilty condition
$[\nen,\nen]\subset\nen$, has a very nice geometric interpretation in
each of the real signatures. In particular, in the real index zero case,
the theorem gives if and only if conditions for the local existence of a
hermitian structure on a 4-manifold \cite{optical,phd}.  

\section{Signature independent Newman-Penrose formalism}\label{npf}
The purpose of this section is to establish a version of the 
Newman-Penrose formalism \cite{np} - a very convenient tool
to study the properties of 4-dimensional manifolds equipped with
a metric - in such a way that it will be usable in the following 
different settings. These are:
\begin{itemize}
\item[(a)] $\mathcal M$ is a \emph{complex} 4-dimensional manifold, and $g$
  is a \emph{holomorphic} metric on $\mathcal M$,
\item[(b)] $\mathcal M$ is a \emph{real} 4-dimensional manifold, and $g$ is
  a \emph{complex valued} metric on $\mathcal M$,
\item[(c)] $\mathcal M$ is a \emph{real} 4-dimensional manifold, and $g$
  is:
\begin{itemize}
\item[(ci)] real of Lorentzian signature,
\item[(cii)] real of Euclidean signature,
\item[(ciii)] real of split signature,
\item[(civ)] a complexification of a real metric having one of 
the above signatures. 
\end{itemize}
\end{itemize}
The classical Newman-Penrose formalism was devised for the case where
$\mathcal M$ is real, and $g$ is Lorentzian. Although the generalization of the
formalism, applicable to all the above settings, is implicit in the
formulation given in the Penrose and Rindler monograph \cite{pr}, one
needs to have some experience to use it in the cases (cii) and (ciii). For
this reason we decided to derive the formalism from first
principles, emphasizing from the very begining how to apply it to the
above different situations. To achieve
our goal of very easy applicability of this formalism to these different
situations, we have introduced 
a convenient notation, in various instances quite different from the
Newman-Penrose original. Since the Newman-Penrose formalism
proved to be a great tool in the study of Lorenztian 4-manifolds, we
believe that our formulation, explained here from the basics, will 
help the community of mathematicians working with 4-manifolds having 
metrics of Euclidean or split signature to appreciate this tool.

From now on $({\mathcal M},g)$ is a 4-dimensional
\emph{real} or \emph{complex} manifold
equipped with a \emph{complex} valued metric. This means that the
metric $g$ is a nondegenerate symmetric bilinear 
form, $g:{\rm T}^\bbC{\mathcal M}\times {\rm T}^\bbC{\mathcal M}\to
\bbC$, with values in the complex numbers \cite{twistor}. 

Given $g$ we use a (local) \emph{null} coframe
$(\theta^1,\theta^2,\theta^3,\theta^4)=(M,P,N,K)$ on $\mathcal M$ in 
which $g$ is 
\be
g=g_{ab}\theta^a\theta^b=2 (M P+N K).\label{npg}\ee
Here, and in the following, formulae like $\theta^a\theta^b$ denote 
the symmetrized
tensor product of the complex valued 1-forms $\theta^a$ and 
$\theta^b$: 
$\theta^a\theta^b=\tfrac12(\theta^a\otimes \theta^b+\theta^b\otimes \theta^a).$

\begin{remark}\label{r1}
Note that our setting, although in general complex, includes all the 
\emph{real} cases. These cases correspond to metrics $g$ such that
$g(X,Y)$ is \emph{real} for all \emph{real} vector fields $X,Y\in {\rm
  T}\mathcal M$. In other words, in such cases the metric $g$
restricted to the tangent space ${\rm T}\mathcal M$ of $\mathcal M$ is
real. If $\mathcal M$ is equipped with a metric $g$ satisfying
  this condition, then we always locally have a null coframe
  $(\theta^1,\theta^2,\theta^3,\theta^4)=(M,P,N,K)$ in which 
\begin{itemize}
\item[($E$)] $P=\bar{M}$ and $K=\bar{N}$ if the
  metric $g_{|T\mathcal M}$ has \emph{Euclidean} signature,
\item[($S_c$)] $P=\bar{M}$ and $K=-\bar{N}$, if the
  metric $g_{|T\mathcal M}$ has \emph{split} signature,
\item[($L$)] $P=\bar{M}$, $N=\bar{N}$ and $K=\bar{K}$, if the
  metric $g_{|T\mathcal M}$ has \emph{Lorentzian} signature.
\end{itemize}
%We have chosen the complex setting to make our considerations
%signature independent when $g_{|T\mathcal M}$ happens to be real. 
\end{remark} 
\begin{remark}\label{r2}
The main statement above about the cases ($E$), ($S_c$) and
($L$) can be rephrased as follows: In the complexification of the
cotangent space of $T^{*\bbC}\mathcal M$, one can introduce three
different \emph{real structures} by appropriate \emph{conjugation
  operators}: `bar'. On the basis of the 1-forms
$(\theta^1,\theta^2,\theta^3,\theta^4)=(M,P,N,K)$ these are defined 
according to:
\begin{itemize}
\item[($E$)] $\bar{M}=P$, $\bar{P}=M$, $\bar{N}=K$ and
  $\bar{K}=N$. With this choice of the conjugation, $g_{|T\mathcal M}$
  is real and has \emph{Euclidean} signature.  
\item[($S_c$)] $\bar{M}=P$, $\bar{P}=M$, $\bar{N}=-K$ and
  $\bar{K}=-N$. With this choice of the conjugation, $g_{|T\mathcal M}$
  is real and has \emph{split} signature. 
\item[($L$)] $\bar{M}=P$, $\bar{P}=M$, $\bar{N}=N$ and
  $\bar{K}=K$. With this choice of the conjugation, $g_{|T\mathcal M}$
  is real and has \emph{Lorentzian} signature.   
\end{itemize}    
Note also that the labels $a=1,2,3,4$ of the null coframe components
$\theta^a$, behave in the following way under these conjugations:
\begin{itemize}
\item[($E$)] $\bar{1}\to 2$, $\bar{2}\to 1$, $\bar{3}\to 4$, $\bar{4}\to
  3$ in the Euclidean case,
\item[($S_c$)] $\bar{1}\to 2$, $\bar{2}\to 1$, $\bar{3}\to -4$, $\bar{4}\to
  -3$ in the split case,
\item[($L$)] $\bar{1}\to 2$, $\bar{2}\to 1$, $\bar{3}\to 3$, $\bar{4}\to
  4$ in the Lorentzian case.
\end{itemize}
These transformations of indices under the respective complex
conjugations will be important when we perform complex
conjugations on multiindexed quantities, such as for example,
$R_{abcd}$. In particular, the above transformation of indices imply,
for example, that in the ($S_c$) case $\bar{R}_{1323}=R_{2414}$, 
$\bar{R}_{1321}=-R_{2412}$, and so on. 
\end{remark}

\begin{remark}\label{r3}
We denoted the split signature case by the letter $S$ with a subscript
$c$ to distinguish this case from the case $S_r$ in which the field of
2-planes anihilating the coframe 1-forms $P$ and $K$ in ($S_c$)
is \emph{totally real}. It is well known \cite{koptra}, that if the
metric $g_{|{\rm T}\mathcal M}$ has split signature, one can choose
a \emph{totally real} null coframe on $\mathcal M$, such that 
\begin{itemize}
\item[($S_r$)] $\bar{M}=M$, $\bar{P}=P$, $\bar{N}=N$,
  $\bar{K}=K$.
\end{itemize}
This situation, although less generic \cite{koptra} than ($S_c$) is worthy of consideration,
since in the integrable case of the Goldberg-Sachs theorem it leads to
the foliation of $\mathcal M$ by \emph{real} 2-dimensional leaves,  
corresponding to the distribution of totally null 2-planes.    
\end{remark}

Given a null coframe $(\theta^a)$ we calculate the differentials of
its components \be
\der\theta^a=-\tfrac12
c^a_{~bc}\theta^b\dz\theta^c.\label{ncn}\ee Following Newman and Penrose \cite{np}, and the
tradition in General Relativity literature \cite{KSMH}, we will assign
Greek letter names to the coefficient functions $c^a_{~bc}$. As 
is well known these coefficients naturally split onto two groups with 12 complex
coefficients in each group. They correspond to two spin connections
associated with the metric $g$. The 12 coefficients from the first
group will be denoted by
$\al,\bet,\gam,\lambda,\mu,\nu,\rho,\si,\tau,\eps,\kappa,\pi$. The 12 
coefficients from the second
group will be denoted by putting primes on the same Greek
letters. The
`primed' and `unprimed' quantities, as describing two different
spinorial connections, will be treated as independent objects in the
complex setting. Their relations to the complex conjugation in the
real settings will be described in Reamark \ref{r4}. 
This said, we write the four equations (\ref{ncn}) as:
\begin{eqnarray}
\der\te^1&=& (\al-\bet') \te^1\dz\te^2 + ( \gam-\gam' - \mu)
\te^1\dz\te^3 + 
(\eps-\eps'  - \rho') \te^1\dz\te^4 -\nonumber\\&& \la \te^2\dz\te^3 - 
  \si' \te^2\dz\te^4 + (\pi - \tau') \te^3\dz\te^4\nonumber\\
\der\te^2&=& (\bet-\al') \te^1\dz\te^2 - \la' \te^1\dz\te^3 - 
  \si \te^1\dz\te^4 + \nonumber\\&&(\gam' - \gam - \mu') \te^2\dz\te^
     3 + (\eps' - \eps - \rho) \te^2\dz\te^4 + (\pi' - \tau) \te^
     3\dz\te^4\label{ca0}\\
\der\te^3 &=&(\rho'- \rho) \te^1\dz\te^2 + (\al' + \bet - \tau) \te^
     1\dz\te^3 - 
  \kappa \te^1\dz\te^4 +\nonumber\\&& (\al + \bet' - \tau') \te^2\dz\te^3 -
  \kappa' \te^2\dz\te^4 - (\eps' + \eps) \te^3\dz\te^4\nonumber\\
\der\te^4&=&(\mu - \mu') \te^1\dz\te^2 - 
  \nu' \te^1\dz\te^3 - (\al' +\bet + \pi') \te^1\dz\te^4 - \nonumber\\&&
  \nu \te^2\dz\te^3 - (\al +\bet' + \pi) \te^2\dz\te^4 -(\gam' + 
     \gam) \te^3\dz\te^4.\nonumber
\end{eqnarray}
This notation for the coefficient functions $c^a_{~bc}$, although ugly
at first sight, has many
advantages. One of them is the already mentioned property of
separating the two spin connections associated with the metric $g$ by
associating them with the respective `primed' and `unprimed' objects. More explicitly, 
defining the Levi-Civita connection 1-forms $\Gamma^a_{~b}$ by  
\begin{eqnarray}
\der\theta^a+\Gamma^a_{~b}\dz\theta^b=0\label{ca1}\\
\Gamma_{ab}=-\Gamma_{ba},\quad\quad \Gamma_{ab}=g_{ac}\Gamma^c_{~b},\nonumber
\end{eqnarray}
we get the following expressions for $\Gamma_{ab}$:
%/home/pawel/notebooks/stonybrook/godlbergsachs/nplorentzgspaperroddenny.nb
\begin{eqnarray}
\tfrac12(\Gamma_{12}+\Gamma_{34})&=&\al' \te^1 + \bet' \te^2 + \gam'\te^3
  +\eps' \te^4\nonumber\\
\Gamma_{13} &=& \la' \te^1 + \mu' \te^2 + \nu' \te^3 + \pi' \te^4\label{spc}\\
\Gamma_{24} &=& \rho' \te^1 + \si' \te^2 + \tau' \te^3 + \kappa' \te^4.\nonumber
\end{eqnarray}
\begin{eqnarray}
\tfrac12(-\Gamma_{12}+\Gamma_{34})&=&\bet \te^1 + \al \te^2 + \gam\te^3
  +\eps 
\te^4\nonumber\\
\Gamma_{23} &=& \mu \te^1 + \la \te^2 + \nu \te^3 + \pi \te^4\label{spcp}\\
\Gamma_{14} &=& \si \te^1 + \rho \te^2 + \tau \te^3 + \kappa \te^4.\nonumber
\end{eqnarray} 
The two spin connections correspond to
$\chi'=(\Gamma_{24},\tfrac12(\Gamma_{12}+\Gamma_{34}),\Gamma_{13})$ and 
$\chi=(\Gamma_{14},\tfrac12(-\Gamma_{12}+\Gamma_{34}),\Gamma_{23})$,
respectively.

\begin{remark}\label{r4}
The above notation is an adaptation of the \emph{Lorentzian version} of the
Newman-Penrose formalism. This can be easily seen, taking into
account the reality conditions discussed in Remarks \ref{r1},
\ref{r2}. In particular, in the Lorentzian case (L), the complex
conjugation defined in Remark \ref{r2}, applied to the quantities
$\al,\bet,\gam,\dots$, yields:
\begin{itemize}
\item[($L$)] $\quad\quad\quad\quad\bma\bar{\al}&\bar{\bet}&\bar{\gam}&\bar{\eps}\\
\bar{\la}&\bar{\mu}&\bar{\nu}&\bar{\pi}\\
\bar{\rho}&\bar{\si}&\bar{\tau}&\bar{\kappa}
\ema=\bma\al'&\bet'&\gam'&\eps'\\
\la'&\mu'&\nu'&\pi'\\
\rho'&\si'&\tau'&\kappa'
\ema$.
\end{itemize}
Thus in the Lorentzian case the complex conjugation changes `unprimed'
Greek letters into `primed' ones and \emph{vice
  versa}. Therefore in this signature the `primed' Greek letter quantities
are totally determined by the `unprimed' ones. The situation is
drastically different in the two other real signatures. There the
`primed' Greek letter quantities are \emph{independent} of the
`unprimed' ones. On the other hand in these two cases, there are some relations between
the quantities \emph{within} each of the `primed' and `unprimed'
family. In the Euclidean case they are given by 
\begin{itemize}
\item[($E$)] $\quad\quad\quad\quad\bma\bar{\al}&\bar{\bet}&\bar{\gam}&\bar{\eps}\\
\bar{\la}&\bar{\mu}&\bar{\nu}&\bar{\pi}\\
\bar{\rho}&\bar{\si}&\bar{\tau}&\bar{\kappa}
\ema=\bma-\bet&-\al&-\eps&-\gam\\
\si&\rho&\kappa&\tau\\
\mu&\la&\pi&\nu
\ema$,
\end{itemize}
with the same relations after the replacement of all `unprimed'
quantities by
their `primed' counterparts on both sides. 

In the split signature cases, we have  
\begin{itemize}
\item[($S_c$)] $\quad\quad\quad\quad\bma\bar{\al}&\bar{\bet}&\bar{\gam}&\bar{\eps}\\
\bar{\la}&\bar{\mu}&\bar{\nu}&\bar{\pi}\\
\bar{\rho}&\bar{\si}&\bar{\tau}&\bar{\kappa}
\ema=\bma-\bet&-\al&\eps&\gam\\
-\si&-\rho&\kappa&\tau\\
-\mu&-\la&\pi&\nu
\ema$,
\end{itemize} and
\begin{itemize}
\item[($S_r$)] $\quad\quad\quad\quad\bma\bar{\al}&\bar{\bet}&\bar{\gam}&\bar{\eps}\\
\bar{\la}&\bar{\mu}&\bar{\nu}&\bar{\pi}\\
\bar{\rho}&\bar{\si}&\bar{\tau}&\bar{\kappa}
\ema=\bma\al&\bet&\gam&\eps\\
\la&\mu&\nu&\pi\\
\rho&\si&\tau&\kappa
\ema$,
\end{itemize}
again with the identical relations for the `primed' quantities. 
\end{remark} 

Now we pass to the `prime'--`unprime' decomposition of the
curvature. The Riemann tensor 
coefficients $R^a_{~bcd}$ are defined
by Cartan's second structure equations:
\begin{eqnarray}
&&\der\Gamma^a_{~b}+\Gamma^a_{~c}\dz\Gamma^c_{~b}=\tfrac12 R^a_{~bcd}\theta^c\dz\theta^d.\label{ca2}
\end{eqnarray}
Due to our conventions,  modulo symmetry, the only
nonzero components of the metric are $g_{12}=g_{34}=1$. The inverse 
of the metric, $g^{ab}$, again modulo symmetry, 
has $g^{12}=g^{34}=1$ as the only nonvanishing components. 
The Ricci tensor is defined as
  $R_{ab}=R^c_{~acb}$. Its scalar is: $R=R_{ab}g^{ab}$, and its
tracefree part is: $\check{R}_{ab}=R_{ab}-\tfrac14 R g_{ab}$. Using the 
metric $g_{ab}$ we also define $R_{abcd}=g_{ae}R^e_{~bcd}$. This is 
further used to define the covariant components of the Weyl tensor $C^a_{~bcd}$ via:
$$
C_{abcd}=R_{abcd}-\tfrac{1}{12} R
  (g_{ac}g_{db}-g_{ad}g_{cb})+\tfrac12(g_{ad}\check{R}_{cb}-g_{ac}\check{R}_{db}+g_{bc}\check{R}_{da}-g_{bd}\check{R}_{ca}).$$
In the context of the present paper, in which the conformal properties
matter, it is convenient to use the Schouten tensor $\Rho$, with help
of which we can write the above displayed equality as
\be
C_{abcd}=R_{abcd}+g_{ad}\Rho_{cb}-g_{ac}\Rho_{db}+g_{bc}\Rho_{da}-g_{bd}\Rho_{ca}
.\label{wyel}
\ee
The Schouten tensor $\Rho$ is a `trace-corrected' Ricci tensor,
with the explicit relation given by
$$\Rho_{ab}=\tfrac12 R_{ab}-\tfrac{1}{12}R g_{ab}.$$

In the Newman-Penrose formalism, the 10 components of the Weyl tensor
are encoded in 10 complex quantities $\Psi_0,\Psi_1,\Psi_2,\Psi_3,\Psi_4$ and
$\Psi_0',\Psi_1',\Psi_2',\Psi_3',\Psi_4'$. Five of them have `primes',
to emphasize that they are associated with the `primed' spin
connection. Another way of understanding this notation 
is to say that the `unprimed'
$\Psi$s are five components of the self-dual part of the Weyl tensor, and the
`primed' $\Psi$s are the components of the anti-self-dual part of the Weyl.

The Ricci and Schouten tensors are mixed `prime'-`unprime' objects,
and as such are not very nicely denoted in the `prime' vs `unprime'
setting. For this reason, when referring to $R_{ab}$, $\check{R}_{ab}$ and
$\Rho_{ab}$, we will not use the Newman-Penrose notation, and will
express these objects using the standard four-dimensional indices
$a=1,2,3,4$, as e.g. in $12(\Rho_{12}+\Rho_{34})=2(R_{12}+R_{34})=R$.

Having said all of this we express Cartan's second  
structure equations (\ref{ca2}), and in particular the curvature
coefficients $R^a_{~bcd}$, in terms of $\Psi$s, $\Psi'$s, $\Rho$ and the
null coframe $(\theta^a)$ as follows:  
\begin{eqnarray}
&&\tfrac12\der(\Gamma_{1 2} +  \Gamma_{34})+\Gamma_{24}\dz\Gamma_{13}=
\nonumber\\ 
&& 
- \Psi_3' \theta^1\dz\theta^3 + 
 \Psi_1' \theta^2\dz\theta^4+\tfrac12 (2\Psi_2' - \Rho_{12}-\Rho_{34})
 (\theta^1\dz\theta^2+\theta^3\dz\theta^4) +\nonumber\\
&&\Rho_{23} \theta^2\dz\theta^3 -\Rho_{14} \theta^1\dz\theta^4 - 
 \tfrac12(\Rho_{12}-\Rho_{34})(\theta^1\dz\theta^2-
 \theta^3\dz\theta^4)\nonumber\\
&&\nonumber\\
&&\der\Gamma_{13}+(\Gamma_{12}+\Gamma_{34})\dz\Gamma_{13}=\nonumber\\ 
&&\Psi_4' \theta^1\dz\theta^3  + (\Psi_2' + \Rho_{12} + \Rho_{34})
  \theta^2
\dz\theta^4-\Psi_3'(\theta^1\dz\theta^2+\theta^3\dz\theta^4) + \label{ca21}\\ 
 && \Rho_{33} \theta^2\dz\theta^3 + \Rho_{11} \theta^1\dz\theta^4-
\Rho_{13}( \theta^1\dz\theta^2 -\theta^3\dz\theta^4)\nonumber\\
&&\nonumber\\
&&\der\Gamma_{24} +\Gamma_{24}\dz(\Gamma_{12}+\Gamma_{34})=\nonumber\\ 
&& (\Psi_2' + \Rho_{12} + 
    \Rho_{34}) \theta^1\dz\theta^3+ 
 \Psi_0' \theta^2\dz\theta^4+\Psi_1'( \theta^1\dz\theta^2+ \theta^3\dz\theta^4) 
+
\nonumber\\
&&\Rho_{22} \theta^2\dz\theta^3  + \Rho_{44} \theta^1\dz\theta^4 + \Rho_{24}( \theta^1\dz\theta^2- \theta^3\dz\theta^4), \nonumber
\end{eqnarray}
with analogous equations for the `unprimed' objects:
\begin{eqnarray}
&&\tfrac12\der(-\Gamma_{1 2} +  \Gamma_{34})+\Gamma_{14}\dz\Gamma_{23}=
\nonumber\\ 
&& 
- \Psi_3 \theta^2\dz\theta^3 + 
 \Psi_1 \theta^1\dz\theta^4-\tfrac12 (2\Psi_2 - \Rho_{12}-\Rho_{34})
 (\theta^1\dz\theta^2-\theta^3\dz\theta^4) +\nonumber\\ 
&&\Rho_{13} \theta^1\dz\theta^3 -\Rho_{24} \theta^2\dz\theta^4 + 
 \tfrac12(\Rho_{12}-\Rho_{34})(\theta^1\dz\theta^2+
 \theta^3\dz\theta^4)\nonumber\\
&&\nonumber\\
&&\der\Gamma_{23}+(-\Gamma_{12}+\Gamma_{34})\dz\Gamma_{23}=\nonumber\\ 
&&\Psi_4 \theta^2\dz\theta^3  + (\Psi_2 + \Rho_{12} + \Rho_{34})
  \theta^1
\dz\theta^4+\Psi_3(\theta^1\dz\theta^2-\theta^3\dz\theta^4) + \label{ca22}\\ 
 && \Rho_{33} \theta^1\dz\theta^3 + \Rho_{22} \theta^2\dz\theta^4+
\Rho_{23}( \theta^1\dz\theta^2 +\theta^3\dz\theta^4)\nonumber\\
&&\nonumber\\
&&\der\Gamma_{14} +\Gamma_{14}\dz(-\Gamma_{12}+\Gamma_{34})=\nonumber\\ 
&& (\Psi_2 + \Rho_{12} + 
    \Rho_{34}) \theta^2\dz\theta^3+ 
 \Psi_0 \theta^1\dz\theta^4-\Psi_1(\theta^1\dz\theta^2-\theta^3\dz\theta^4) 
+
\nonumber\\
&&\Rho_{11} \theta^1\dz\theta^3  + \Rho_{44} \theta^2\dz\theta^4 - 
\Rho_{14}(\theta^1\dz\theta^2+ \theta^3\dz\theta^4). \nonumber
\end{eqnarray}
Note that in the first part (\ref{ca21}) of the structure equations,  the full traceless part of the
Schouten tensor $\Rho$, represented by its nine components 
$\Rho_{11}$, $\Rho_{13}$, $\Rho_{14}$, $\Rho_{22}$, $\Rho_{23}$,
$\Rho_{24}$, $\Rho_{33}$, $\Rho_{44}$ and $\Rho_{12}-\Rho_{34}$,
stays with the basis of the selfdual 2-forms:
\be\Sigma=(\theta^2\dz\theta^3,\theta^1\dz\theta^4,\theta^1\dz\theta^2-\theta^3\dz\theta^4).\label{sela}\ee In
the second part (\ref{ca22}) of the structure equations, the full traceless 
part of the Schouten tensor $\Rho$ appears again, but now at the
basis of the antiselfdual 2-forms:
\be
\Sigma'=(\theta^1\dz\theta^3,\theta^2\dz\theta^4,\theta^1\dz\theta^2+\theta^3\dz\theta^4).\label{asel}\ee On
the other hand the selfdual and the antiselfdual parts of the Weyl
tensor, corresponding to the respective $\Psi$s and $\Psi'$s, are
separated: in equations (\ref{ca21}) we only have  $\Psi'$s, 
whereas in (\ref{ca22}) we only have $\Psi$s. The trace of the Schouten tensor
$2(\Rho_{12}+\Rho_{34})$, proportional to the Ricci scalar $R$,
appears in both sets of equations, always together with the
respective Weyl tensor components $\Psi_2$ and $\Psi_2'$. It is also
worthwhile to mention that if one uses the following basis 
$$E_-=\bma 0&0\\1&0\ema,\quad E_0=\bma 1&0\\0&-1\ema,\quad E_+=\bma 0&-1\\0&0\ema,$$
of the Lie algebra $\sla(2)$, and if one defines 
$$\Gamma=\Gamma_{14}E_-+\tfrac12(-\Gamma_{12}+\Gamma_{34})E_0
+\Gamma_{23}E_+,$$
$$\Gamma'=\Gamma_{24}E_-+\tfrac12(\Gamma_{12}+\Gamma_{34})E_0
+\Gamma_{13}E_+,$$
then the left hand sides of equations (\ref{ca21})-(\ref{ca22}) appear
in the formulae
\begin{eqnarray*}
&&\der\Gamma+\Gamma\dz\Gamma=\bma \tfrac12\der(-\Gamma_{1
  2} +
\Gamma_{34})+\Gamma_{14}\dz\Gamma_{23}&-\der\Gamma_{23}-(-\Gamma_{12}+\Gamma_{34})\dz\Gamma_{23}\\
\der\Gamma_{14}
+\Gamma_{14}\dz(-\Gamma_{12}+\Gamma_{34})&-\tfrac12\der(-\Gamma_{1
  2} +
\Gamma_{34})-\Gamma_{14}\dz\Gamma_{23}
\ema,
\end{eqnarray*}
\begin{eqnarray*}
&&\der\Gamma'+\Gamma'\dz\Gamma'=\bma \tfrac12\der(\Gamma_{1
  2} +
\Gamma_{34})+\Gamma_{24}\dz\Gamma_{13}&-\der\Gamma_{13}-(\Gamma_{12}+\Gamma_{34})\dz\Gamma_{13}\\
\der\Gamma_{24}
+\Gamma_{24}\dz(\Gamma_{12}+\Gamma_{34})&-\tfrac12\der(\Gamma_{1
  2} +
\Gamma_{34})-\Gamma_{24}\dz\Gamma_{13}
\ema.
\end{eqnarray*}
This explains the term `spin connections' assigned to the previously 
defined quantities $\chi$ and $\chi'$. It also justifies the 
`prime'-`unprime' notation, which is rooted in the speciality of
4-dimensions, stating that for $n\geq 3$ the Lie algebra 
${\frak so}(n,\bbC)$ is not simple only when $n=4$, and in that case it has 
the symmetric split: 
${\frak so}(4,\bbC)={\frak sl}(2,\bbC)\oplus{\frak sl}(2,\bbC)$. This
enables us to split the ${\frak so}(4,\bbC)$-valued Levi-Civita
connection into the well defined ${\frak sl}(2,\bbC)$-valued `primed'
and `unprimed' parts, which are totally independent. In real
signatures we have an analogous split for ${\frak so}(4-p,p)={\frak
  g}\oplus{\frak g}'$, $p=0,1,2$, where now $\frak g$ and ${\frak g}'$
are two copies of the appropriate real form of ${\frak
  sl}(2,\bbC)$. This again enables us to split the Levi-Civita connection into the
`primed' and `unprimed' connections, with the appropriate reality
conditions, as in ($E$), ($S_c$), ($S_r$) or ($L$).

Comparing equations (\ref{spc})-(\ref{spcp}) with
(\ref{ca21})-(\ref{ca22}), one finds relations between the curvature
quantities $\Rho$, $\Psi$ and $\Psi'$ and the first derivatives 
of the connection coefficients
$\alpha,\beta,\dots,\alpha',\beta',\dots$. These relations are called
the Newman-Penrose equations \cite{np}. We present them in the
Appendix. In these equations, and in the rest of the paper, we denote 
the vector fields dual on $\mathcal M$ to the null coframe $(M,P,N,K)$ by the
respective symbols $(\bel,\pel,\vel,\del)$. Thus we have e.g. $\bel\hook
M=1$, and zero on all the other coframe components, $\del\hook N=0$,
etc. Also, when applying these vector fields to functions on $\mathcal
M$ we omit parentheses. Thus, instead of writing $\del(\al)$ to denote
the derivative of a connection coefficient $\al$ in the direction of
the basis vector field $\del$, we simply write $\del\al$.

In addition to the Newman-Penrose equations we will also need the
commutators of the basis vector fields. These are given by the
formulae dual 
to equations (\ref{ca0}), and read:
\begin{eqnarray}
&&[\bel,\pel]=(\bet'-\al)\bel+(\al'-\bet)\pel+(\rho-\rho')\vel+(\mu'-\mu)\del\nonumber\\
&&[\bel,\vel]=(\mu+\gam'-\gam)\bel+\la'\pel+(\tau-\al'-\bet)\vel+\nu'\del\nonumber\\
&&[\pel,\vel]=\la \bel+(\mu'+\gam-\gam')\pel+(\tau'-\al-\bet')\vel+\nu\del\nonumber\\
&&[\bel,\del]=(\rho'+\eps'-\eps)\bel+\si \pel+\kappa\vel+(\al'+\bet+\pi')\del\label{comw}\\
&&[\pel,\del]=\si'\bel+(\rho+\eps-\eps')\pel+\kappa'\vel+(\al+\bet'+\pi)\del\nonumber\\
&&[\vel,\del]=(\tau'-\pi)\bel+(\tau-\pi')\pel+(\eps'+\eps)\vel+(\gam'+\gam)\del\nonumber
\end{eqnarray}

The Newman-Penrose equations are supplemented by the second Bianchi
identities, which are crucial for the proof of the Goldberg-Sachs
theorem. These are relations between the first derivatives of the
curvature quantities $\Psi$, $\Psi'$ and $\Rho$ and the connection
coefficients. These Bianchi identities are also presented in the Appendix.

\section{Generalizations of the Goldberg-Sachs theorem for complex metrics}\label{prof}
The thesis of the Goldberg-Sachs theorem can be restated in the
language of the Newman-Penrose formalism as follows:

To interpret the integrability condition $[\nen,\nen]\subset\nen$ on the totally
null distribution $\nen$, we align the Newman-Penrose coframe
$(\theta^1,\theta^2,\theta^3,\theta^4)=(M,P,N,K)$ in such a way that
the two \emph{null} and \emph{mutually orthogonal} frame vectors
$e_1=m=\bel$ and $e_4=k=\del$ span $\nen$,
$\nen=\Span_\bbC(\bel,\del)$. Such a coframe on $({\mathcal M},g)$
will be called a coframe \emph{adapted to} $\tn$.

Then the integrability of $\nen$ is totally determined by the
commutator $[\bel,\del]$ of these basis vectors. Looking at this
commutator in (\ref{comw}), we see that the condition that
$[\bel,\del]$ is in the span of $\bel$ and $\del$ is equivalent to 
$\kappa\equiv\si\equiv 0$. Thus we have
\begin{proposition}\label{petu}
Let $\mathcal N$ be a field of selfdual totally null 2-planes on a
4-dimensional manifold $\mathcal M$ with the metric $g$. Let
$(m,p,n,k)$ be a null frame in ${\mathcal U}\subset\mathcal M$ adapted
to $\mathcal N$. Then the field $\nen={\rm Span}_\bbC(m,k)$ is integrable, 
$[\nen,\nen]\subset\nen$,  
in $\mathcal U$ if and only if the frame connection coefficients
$\Gamma_{144}=\kappa$ and $\Gamma_{141}=\sigma$ vanish identically, 
$\kappa\equiv\si\equiv 0$, in $\mathcal U$. 
\end{proposition}
To interpret the algebraic speciality of the selfdual part of the Weyl tensor, 
we focus on the condition 
\be
C(m,k,m,k)\equiv 0.\label{pndpp}\ee 
Here we consider the Weyl tensor $C_{abcd}$ as a linear 
map $C:\bigotimes^4 {\rm T}^\bbC{\mathcal M}\to\bbC$. Note that, since the so
understood Weyl tensor is \emph{antisymmetric} in the first two
arguments, as well as, independently, in the last two arguments, the
vanishing in equation (\ref{pndpp}), although defined on a
particular basis of $\nen$, is basis independent. Actually, if we 
think of $C$ as a linear map $C:(\bgw^2{\rm T}^\bbC{\mathcal M})\odot(\bgw^2{\rm T}^\bbC{\mathcal M})\to \bbC$, and
identify a 2-dimensional totally null distribution $\nen$ with the
complex line bundle $$\nen_\dz=\{w\in\bgw^2{\rm T}^\bbC{\mathcal M}~|~w=v_1\dz
v_2,~v_1,v_2\in\nen\},$$ then we say that $\nen$ is a 
\emph{principal} totally null distribution iff 
\be
C(\nen_\dz,\nen_\dz)\equiv 0.\label{pndp}
\ee  
\begin{remark}
The quantity
$C(m,k,m,k)$ is a null counterpart of the sectional curvarture
  from Riemannian geometry. In fact, given a 2-dimensional vector
  space $V={\rm Span}_\bbR(X,Y)$, the sectional curvature associated
  with $V$ is $$K=K(X,Y)=\frac{g(R(X,Y)X,Y)}{|X\dz Y|^2}.$$ 
The appearence of the
  denominator $|X\dz Y|^2=g(X,X)g(Y,Y)-g(X,Y)^2$ in this expression
    makes this quantity 
  independent of the choice of $X$, $Y$ in $V$. The notion of
  sectional curvature loses its meanning for vector spaces $V$ which
  are totally null, since for them the metric $g$ when restricted to
  $V$ vanishes, making the denominator $|X\dz Y|^2\equiv 0$ for all
  $X,Y\in V$. To incorporate totally null vector spaces $V$, one needs
  to generalize the notion of sectional curvature, removing the
  denominator from its definition. This leads to the quantity 
$$K_0=K_0(X,Y)=g(R(X,Y)X,Y).$$ This, although basis dependent, 
  transforms in a homogeneous fashion, $$K_0(X,Y)\to
  (ad-bc)^2K_0(X,Y),$$ under the change of basis 
$X\to aX+bY$, $Y\to cX+dY$. Thus
  vanishing or not of $K_0$ is an invariant property of any 2-dimensional 
vector space $V\subset T_x{\mathcal M}$. This property of having
$K_0$ equal or not equal to zero, characterizes $V$ and is well
defined regardless of the fact if the metric is real or complex,
including the cases when $V$ is totally null. 

Now, passing to the specific situation of 4-dimensional manifolds, we
can choose $V$ to be a field of selfdual totally null 2-planes
$\nen$. More specifically, if $\nen={\rm Span}_\bbC(m,k)$, we
  easily check (see (\ref{ca22})) that
  $K_0(m,k)=C(m,k,m,k)=\Psi_0$. Thus $K_0(m,k)$ is the $\Psi_0$ component of the
  \emph{selfdual} part of the Weyl tensor.  For an \emph{antiselfdual} totally null plane
  $\nen'={\rm Span}_\bbC(p,k)$ we have $K_0(p,k)=C(p,k,p,k)=\Psi_0'$,
  which is the corresponding component of the \emph{antiselfdual} part
  of the Weyl tensor.  
This shows that the
  principal selfdual totally null 2-planes are just those for which the 
  quantity $\Psi_0$ vanishes. Thus, in a sense, the principal selfdual
  totally null 2-planes have vanishing sectional curvature. (We have
  also an analogous statement for the principal antiselfdual 2-planes;
  they are related to the antiselfdual part of the Weyl tensor, and are defined by 
the vanishing of the quantity $\Psi_0'$.)   
\end{remark}
Let us now choose a Newman-Penrose coframe $(M,P,N,K)$ which is 
not related to any particular choice of $\nen$. Thus we have 
$g=2(MP+NK)$. Then, at every
point of $\mathcal M$, we have
two families $\nen_z$ and $\nen_{z'}$ of 
2-dimensional totally null planes \cite{twistor}. These two families
are parametrized by a complex parameter $z$ or $z'$,
respectively, and the 2-planes parametrized by $z$ are selfdual, and
those parametrized by $z'$ are antiselfdual. In
terms of the frame $(e_1,e_2,e_3,e_4)=(m,p,n,k)=(\bel,\pel,\vel,\del)$
dual to $(M,P,N,K)$, they are given by 
\be
\nen_z=\Span_\bbC(m+zn,k-zp),\quad\quad\quad z\in\bbC,\label{plan}
\ee
and
\be
\nen_{z'}=\Span_\bbC(p+z'n,k-z'm),\quad\quad\quad z'\in\bbC.
\ee
Adding a totally null plane $\nen_\infty=\Span_\bbC(n,p)$ to the first
family, and  $\nen_{\infty'}=\Span_\bbC(n,m)$ to the second family, we
have two \emph{spheres} of 2-dimensional totally null planes at each
point of $\mathcal M$. The first sphere consists of the selfdual 2-planes,
the second of the antiselfdual 2-planes.
 
Now we find the \emph{principal} 2-planes in each of these
spheres. The principal 2-planes in the first sphere correspond to
those $z$
such that
\be
C(m+zn,k-zp,m+zn,k-zp)=0.\label{ppn}\ee
The left hand side of this equation is a \emph{fourth} order
polynomial in the \emph{complex} variable $z$, thus (\ref{ppn}) 
treated as an equation for $z$, has \emph{four} roots, some of which
may be \emph{multiple} roots. Moreover, equation (\ref{ppn}) written
explicitly in terms of the Newman-Penrose Weyl coeffcients $\Psi$s
and $\Psi'$s, involves only the `unprimed' quantities. Explicitly:
$$C(m+zn,k-zp,m+zn,k-zp)=\Psi_4 z^4-4\Psi_3 z^3+6\Psi_2 z^2+4\Psi_1 z+\Psi_0,$$
where we have used the conventions of the previous section, such as
$C(m,k,m,k)=\Psi_0$, etc. Similar considerations for the second sphere
lead to the following proposition:
\begin{proposition}\label{44}
A selfdual totally null 2-plane $\nen_z=\Span_\bbC(m+zn,k-zp)$ is principal at
$x\in\mathcal M$ iff $z$ is a root of the equation
\be\Psi_4 z^4-4\Psi_3 z^3+6\Psi_2 z^2+4\Psi_1 z+\Psi_0=0.\label{4}\ee
An antiselfdual totally null 2-plane $\nen_{z'}=\Span_\bbC(m+z'k,n-z'p)$ is principal at
$x\in\mathcal M$ iff $z'$ is a root of the equation
\be
\label{4a}\Psi_4' {z'}^4-4\Psi_3' {z'}^3+6\Psi_2' {z'}^2+4\Psi_1' z'+\Psi_0'=0.\ee
\end{proposition}
Thus at every point of $\mathcal M$ we have at most four selfdual
principal null 2-planes and at most four antiselfdual principal null
2-planes. If a principal null 2-plane corresponds to a multiple root of
(\ref{4}) or (\ref{4a}), then such a 2-plane is called a \emph{multiple}
principal null 2-plane. A selfdual
or antiselfdual part of the Weyl tensor with multiple principal 2-planes at a
point is called \emph{algebraically special} at this point.

We also note that the number and the 
multiplicity of the roots in (\ref{4}) or (\ref{4a}) is a
\emph{conformal invariant} of the metric at a point. Thus the
algebraically special cases can be further stratified according to the number
of the roots and their multiplicities. 

The possibilities here for
(\ref{4}) are: a)
three distinct roots, b) two distinct roots, with one of multiplicity
three, c) two distinct roots,
each with multiplicity two, d) one root of multiplicity four, e) selfdual part of the Weyl
tensor is zero. We have also the corresponding possibilities a'), b')
c'), d') and e') for (\ref{4a}). 
\begin{definition}\label{pett}
The selfdual part of the Weyl tensor \emph{is of
Petrov type {\rm II, III, D, N}, or {\rm 0} at a point}, if equation (\ref{4}) has
roots as in the respective cases a), b), c), d) and e) at this
point. If the Petrov type of the selfdual part of the Weyl tensor varies in $\mathcal M$, 
from point to point, but only between the types II and D, we say that it is
of type $\overline{\rm II}$. 
The analogous classification holds also for the antiselfdual
part of the Weyl tensor.   
\end{definition}

\begin{remark}
Suppose that the selfdual part of the Weyl tensor of $({\mathcal M},g)$ does not vanish at
each point of a neighbourhood ${\mathcal U}'\subset\mathcal M$. 
Thus at
every point of ${\mathcal U}'$ we have at least one principal totally null
2-plane. We now take the principal null 2-plane which at
$x\in{\mathcal U}'$ has the
smallest multiplicity $1\leq q\leq 4$. There always exists a
neighbourhood ${\mathcal U}\subset{\mathcal U}'$ of $x$ in which this
principal totally null 2-plane extends to a field $\mathcal N$ of 
principal totally null 2-planes of multiplicity not bigger than $q$. 
In $\mathcal U$ we choose a null frame $(m,p,n,k)$ 
in such a way that ${\rm Span}_\bbC(m,k)=\mathcal N$. In this frame
the definition (\ref{plan}) shows that ${\mathcal  N}={\mathcal N}_0$,
i.e. that the corresponding  $z=0$ in $\mathcal U$. Moreover since $\mathcal N$, as a field of 
principal null 2-planes in $\mathcal U$ satisfies (\ref{4}), then  
 $\Psi_0\equiv 0$ everywhere in this frame.
\end{remark}

This proves the following 
\begin{proposition}
Around every point $x$ of a manifold $({\mathcal M},g)$ 
with nowhere vanishing selfdual part of the Weyl tensor, there exists a
neighbourhood ${\mathcal U}$ and a null frame $(m,p,n,k)$ in
${\mathcal U}$ in which $\Psi_0\equiv 0$ everywhere. 
\end{proposition}  

Now if the selfdual part of the Weyl tensor is algebraically special
of type II in $\mathcal U$, with $\mathcal N$ the corresponding 
principal multiple field of totally null 2-planes, then in $\mathcal
U$ we choose a null frame $(m,p,n,k)$ adapted to $\mathcal N$. In this
frame ${\mathcal
  N}={\mathcal N}_0={\rm Span}(m,k)$, the value $z=0$ is a double root of (\ref{4}), and since this is true at every
point of $\mathcal U$, we have $\Psi_0\equiv\Psi_1\equiv 0$. Performing similar considerations for
types III and N, and forcing $z=0$ to be a root of the
equation (\ref{4}) with the respective locally constant multiplicity
$q=1,2,3$ and 4, we get the following  
\begin{proposition}\label{pett1}
Let ${\mathcal N}$ be a field of principal totally null 2-planes for
the selfdual part of the Weyl tensor of a metric $g$ on a
4-dimensional manifold $\mathcal M$. Assume that $\mathcal N$ has a
constant multiplicity $q$ in a neighbourhood $\mathcal U$ in $\mathcal
M$. Then one can choose a null frame $(m,p,n,k)$ in $\mathcal U$, with
${\mathcal N}={\rm Span}(m,k)$ and $g=2(MP+NK)$, so that  
\begin{itemize}
\item if $q=1$ then in this frame $\Psi_0\equiv 0$ and $\Psi_1\neq 0$,
\item if $q=2$ then in this frame $\Psi_0\equiv \Psi_1\equiv 0$ and
  $\Psi_2\neq 0$, 
\item if $q=3$ then in this frame $\Psi_0\equiv \Psi_1\equiv 
\Psi_2\equiv 0$ and $\Psi_3\neq 0$,
\item if $q=3$ then in this frame $\Psi_0\equiv \Psi_1\equiv
  \Psi_2\equiv\Psi_3\equiv 0$ and $\Psi_4\neq 0$. 
\end{itemize}
Conversely, if we have a null frame in $\mathcal U$ in which 
\begin{itemize}
\item $\Psi_0\equiv \Psi_1\equiv
  \Psi_2\equiv\Psi_3\equiv 0$ and $\Psi_4\neq 0$ then 
${\mathcal N}={\rm Span}(m,k)$ is a field of multiple principal
  2-planes in $\mathcal U$ with multiplicity $q=4$,
\item $\Psi_0\equiv \Psi_1\equiv
  \Psi_2\equiv 0$ and $\Psi_3\neq 0$ then 
${\mathcal N}={\rm Span}(m,k)$ is a field of multiple principal
  2-planes in $\mathcal U$ with multiplicity $q=3$,
\item $\Psi_0\equiv \Psi_1\equiv 0$ and $\Psi_2\neq 0$ then 
${\mathcal N}={\rm Span}(m,k)$ is a field of multiple principal
  2-planes in $\mathcal U$ with multiplicity $q=2$,
\item $\Psi_0\equiv 0$ and $\Psi_1\neq 0$ then 
${\mathcal N}={\rm Span}(m,k)$ is a field of multiple principal
  2-planes in $\mathcal U$ with multiplicity $q=1$.
\end{itemize}
\end{proposition}
This immediately implies
\begin{corollary}
The selfdual part of the Weyl tensor of a metric $g$ on a 4-dimensional manifold
$\mathcal M$ is algebraically special in neighbourhood $\mathcal U$,
with $\mathcal N$ being a field of multiple principal 2-planes in
$\mathcal U$ if and only if there exists a null frame 
$(m,p,n,k)$ in $\mathcal U$ in which $\Psi_0\equiv\Psi_1\equiv 0$ in
$\mathcal U$. In this frame ${\mathcal N}={\rm Span}_\bbC(m,k)$. 
\end{corollary}
\subsection{Generalizing the Przanowski-Pleba\'nski version}
The starting point for our generalizations of 
the Goldberg-Sachs theorem is to replace the
Ricci flat condition from the classical version \cite{gs}, by a condition on only that part of the Ricci 
tensor, which is `visible' to the integrable totally null 2-plane $\nen$.

For this we consider the Ricci tensor of 
$({\mathcal M},g)$ as a 
\emph{symmetric}, possibly \emph{degenerate}, \emph{bilinear form} on
$\mathcal M$. 
We denote it by $Ric$ and extend it to the complexification ${\rm
  T}^\bbC \mathcal M$ by linearity. Now given a complex distribution ${\mathcal
  Z}\subset{\rm T}^\bbC \mathcal M$ we say that the Ricci tensor is degenerate
on ${\mathcal Z}$, 
$$Ric_{|{\mathcal Z}}= 0,\quad\quad{\rm iff}\quad\quad
Ric(Z_1,Z_2)=0,\quad 
\forall Z_1,Z_2\in{\mathcal Z}.$$ 

Then we have the following theorem:
\begin{theorem}\label{gscc}
Let $\tn\subset {\rm T}^\bbC \mathcal M$ be a field of totally null
2-planes on a 4-dimensional manifold
$({\mathcal M},g)$ equipped with a real metric $g$ of any signature. 
Assume that the Ricci tensor $Ric$ of 
$({\mathcal M},g)$, considered as a symmetric bilinear form
on ${\rm T}^\bbC \mathcal M$, is degenerate on $\tn$,
$$Ric_{|\tn}=0.$$
If in addition the field $\tn$ is integrable, $[\tn,\tn]\subset\tn$, everywhere
  on $\mathcal M$, then 
$({\mathcal M},g)$ is algebraically special at every point, with a 
  field of multiple principal totally null 2-planes tangent to $\tn$.
\end{theorem} 

To prove it, we fix a null frame $(m,p,n,k)$ on $\mathcal M$ adapted
to $\nen$. This means that ${\mathcal N}={\rm Span}_\bbC(m,k)$. 

It is then very easy to see that the vanishing of the Ricci tensor on 
$\nen$ is, due to our
conventions, equivalent to the
conditions $$\Rho_{11}\equiv\Rho_{14}\equiv\Rho_{44}\equiv 0.$$

Instead of proving Theorem \ref{gscc}, we prove a theorem that implies
it. This is  
the complex version of the Goldberg-Sachs theorem, which generalizes
the Lorentzian version due to Przanowski and Plebanski \cite{plprz2}. 
When stated in the Newman-Penrose
language, this reads as follows:
\begin{theorem}\label{gst}~\\
(1) Suppose that a 4-dimensional metric $g$ satisfies
$\Rho_{11}\equiv\Rho_{14}\equiv\Rho_{44}\equiv 0$ and    
$\kappa\equiv\si\equiv 0$. Then $\Psi_0\equiv\Psi_1\equiv 0$.\\
(2) If $g$ is Einstein, $Ric(g)=\Lambda g$, and has a nowhere
  vanishing selfdual part of the Weyl
tensor, then $\Psi_0\equiv \Psi_1\equiv 0$ implies
$\kappa\equiv\si\equiv 0$.
\end{theorem}
Before the proof we make the following remarks:
\begin{remark}
It is easy to see that part (1) of the above Theorem is equivalent to
Theorem \ref{gscc}.
\end{remark}
\begin{remark}
Note that $Ric=0$ and more generally $Ric=\Lambda g$ are special cases of our
condition $Ric_{|\tn}=0$.
\end{remark}
   
\begin{proof} (of Theorem \ref{gst}). First we assume that $\kappa$ and $\sigma$ vanish everywhere on
$\mathcal M$. To conclude that $\Psi_0\equiv 0$ is very easy: Actually this conclusion is an immediate
consequence of the Newman-Penrose equation (\ref{np1}). For if
$\kappa$ and $\sigma$ are identically vanishing, then equation
(\ref{np1}) gives $\Psi_0\equiv 0$. Note that this conclusion holds 
even without \emph{any} assumption about the components of the
Schouten tensor
$\Rho$ (or the Ricci tensor).

Now we prove the following 
\begin{lemma}\label{gstl}
Suppose that a 4-dimensional metric $g$ satisfies $\kappa\equiv\sigma\equiv0$
and 
\begin{eqnarray}
&&\bel\Psi_1 \equiv 2 (\bet  + 2 \tau)\Psi_1,\label{pr1}\\
&&\del\Psi_1 \equiv 2 (\eps  - 2\rho) \Psi_1.\label{pr2}
\end{eqnarray} 
Then it also satisfies $$\Psi_1\equiv0.$$
\end{lemma}
\begin{proof}
We use the commutator (\ref{comw}), and the Newman-Penrose
equations (\ref{np2})-(\ref{np4}) to obtain the compatibility
conditions for (\ref{pr1}) and (\ref{pr2}). 
This is a pure calculation. We give its main steps below:

\begin{itemize}
\item applying $[\bel,\del]$ to (\ref{pr1}) and (\ref{pr2}) we get:
$$[\bel,\del]\Psi_1\equiv 2\bel\big((\eps  - 2\rho) \Psi_1\big)-2\del\big((\bet  + 2
  \tau)\Psi_1\big);$$
\item next, using (\ref{comw}), and again  (\ref{pr1}) and
  (\ref{pr2}), we transform this identity into:
\begin{eqnarray}
2(\rho'+\eps'-\eps) (\bet  + 2 \tau)\Psi_1+2(\al'+\bet+\pi')(\eps  - 2\rho)\Psi_1\equiv\nonumber\\2\bel\big((\eps  - 2\rho) \Psi_1\big)-2\del\big((\bet  + 2
  \tau)\Psi_1\big);\label{pr3}
\end{eqnarray}
\item now, the Leibniz rule, and a third use of (\ref{pr1}) and
  (\ref{pr2}), enables us to eliminate of the derivatives of
  $\Psi_1$ in (\ref{pr3});
\item actually, simplifying (\ref{pr3}), and using
  (\ref{pr1}), (\ref{pr2}) we get:
\be
\Big(2\bel(\eps-2\rho)-2\del(\bet+2\tau)+2(\eps-\eps'-\rho')(\bet+2\tau)-2(\al'+\bet+\pi')(\eps-2\rho)\Big)\Psi_1\equiv
0;\label{num}
\ee
\item the last step in the proof of the lemma is to use the Newman-Penrose
  equations 
  (\ref{np2})-(\ref{np4}); 
\item these equations 
eliminate $\bel\eps-\del\bet$, (look at (\ref{np2})), $\bel\rho$, (look
at (\ref{np3})), and $\del\tau$, (look at (\ref{np4})), from the identity (\ref{num}); 
\item this makes the identity (\ref{num}) derivative-free;
\item actually it transforms (\ref{num}) to a remarkable identity:
\be
(10\Psi_1)\Psi_1\equiv 0;\label{eqqq}
\ee
\item the identity (\ref{eqqq}) obviously implies $\Psi_1\equiv 0$;
\end{itemize}
This proves Lemma \ref{gstl}.
\end{proof}
To conclude the proof of the part one of Theorem \ref{gst} we use 
our assumptions $\Rho_{11}\equiv\Rho_{14}\equiv\Rho_{44}\equiv0$, $\kappa\equiv\sigma\equiv0$, 
and their consequence $\Psi_0\equiv 0$, and insert them 
in the Bianchi identities (\ref{b1}) and
(\ref{b2}). This trivially gives the relations (\ref{pr1}) and
(\ref{pr2}), respectively. Then an obvious use of Lemma \ref{gstl} 
finishes the proof of part one of Theorem \ref{gst}.

We now pass to the proof of part two of Theorem \ref{gst}. 

\noindent
When going from $\Big(\Psi_0\equiv\Psi_1\equiv 0\Big)$ to
$\Big(\kappa\equiv\sigma\equiv0\Big)$ we do as follows:
\begin{itemize}
\item Initially we only assume that
  $\Rho_{11}\equiv\Rho_{14}\equiv\Rho_{44}\equiv 0$. 
\item Then the Bianchi identities (\ref{b1}) and (\ref{b2}) give:
\be
 2 \Rho_{13}\kappa + (3 \Psi_2 + \Rho_{12} -\Rho_{34})\si\equiv 0\label{wd1}\ee
and
\be(3 \Psi_2  - \Rho_{12}+ \Rho_{34})\kappa  + 2\Rho_{24}\si\equiv 0,\label{wd2}\ee
respectively.
\end{itemize}
At this stage the following remark is in order:
\begin{remark}
If we were able to conclude that the rank of the matrix 
\be 
m=\bma 2 \Rho_{13}& 3 \Psi_2 + \Rho_{12} -\Rho_{34}\\
 3 \Psi_2  - \Rho_{12}+ \Rho_{34}& 2\Rho_{24}\ema\label{maci}\ee
was \emph{identically} equal to \emph{two}, this would immediately yield 
$\kappa\equiv\sigma\equiv 0$, which would conclude the proof.
On the other extreme, if we were sure that the matrix $m$ was \emph{identically} equal
  to zero (i.e if it had rank identically equal to zero), we would
  argue as follows: The identically zero rank of $m$ means that in addition to
  $\Rho_{11}\equiv\Rho_{14}\equiv\Rho_{44}\equiv 0$ we have:
  $\Rho_{13}\equiv\Rho_{24}\equiv\Rho_{12}-\Rho_{34}\equiv\Psi_2\equiv
  0$. Then, combining the Bianchi identities (\ref{b3}) and
  (\ref{b9}), we get 
$$2\Rho_{33}\kappa+2(\Rho_{23}-3\Psi_3)\si\equiv 0.$$
Similarly, using the Bianchi identities (\ref{b4}) and
  (\ref{b10}) we get:
$$2(\Rho_{23}+3\Psi_3)\kappa+2\Rho_{22}\si\equiv 0.$$
Thus, in such case, the situation is similar to the previously
  considered case with the matrix $m$: Now we have 
$$m_1=\bma  \Rho_{33}&- 3 \Psi_3 + \Rho_{23}\\
 3 \Psi_3  + \Rho_{23}& \Rho_{22}\ema,$$
and if $m_1$ has rank \emph{identically} equal to \emph{two}, we conclude
that $\kappa\equiv\sigma\equiv 0$. If it has rank \emph{identically}
equal to \emph{zero}, we in addition have
$\Rho_{33}\equiv\Rho_{22}\equiv\Rho_{23}\equiv\Psi_3\equiv 0$. This,
due to the Bianchi identities, implies also 
that $\Rho_{12}\equiv\Rho_{34}\equiv{\rm const}$. Comparing this with (\ref{b5})
and (\ref{b6}) leads to  
$$\Psi_4\si\equiv\Psi_4\kappa\equiv 0,$$
which if we assume $\Psi_4\neq 0$, yields $\kappa\equiv\sigma\equiv 0$.
\end{remark}

This remark emphasizes that the local properties of the matrices $m$ and
$m_1$ are crucial for the behaviour of $\kappa$ and $\si$.  
Since we have \emph{no guarantee} that rank of e.g. $m$ is locally constant, returning to our proof, we must strenghten our
assumptions on $g$ by requiring that it satisfies more curvature
conditions than $\Rho_{11}\equiv\Rho_{14}\equiv\Rho_{44}\equiv 0$.
\begin{itemize}
\item The additional conditions which enable us to get
  $\kappa\equiv\sigma\equiv 0$ are:
$$\Rho_{13}\equiv\Rho_{22}\equiv\Rho_{23}\equiv
\Rho_{24}\equiv\Rho_{33}\equiv\Rho_{12}-\Rho_{34}\equiv
    0.$$
These, with the already
assumed $\Rho_{11}\equiv\Rho_{14}\equiv\Rho_{44}\equiv 0$, constitute
the \emph{full set of Einstein conditions} $Ric(g)=\Lambda g$, for the metric
$g$.
\item Under the Einstein assumption $Ric(g)=\Lambda g$
  and the requirement that the selfdual part of the Weyl tensor is
  nonvanishing, we get $\kappa\equiv\si\equiv 0$ in a very easy way, by 
a successive inspection of
  the Bianchi identities (\ref{b1}), (\ref{b2}),
  (\ref{b3}),(\ref{b4}),(\ref{b5}),(\ref{b6}). 
\item Indeed, the assumed Einstein equations
  $\Rho_{11}\equiv\Rho_{14}\equiv\Rho_{44}\equiv\Rho_{22}\equiv\Rho_{24}\equiv\Rho_{13}\equiv\Rho_{23}\equiv\Rho_{33}\equiv\Rho_{12}-\Rho_{34}\equiv
  0$, the algebraical speciality conditions $\Psi_0\equiv\Psi_1\equiv 0$, and the Bianchi identities
  (\ref{b1}), (\ref{b2}), give $\si \Psi_2\equiv 0$ and
  $\kappa\Psi_2\equiv 0$.  This means that whenever $\Psi_2\neq 0$ we
  have $\kappa\equiv\si\equiv 0$. By continuity 
  the points in which $\kappa$ or $\sigma$ are
  nonzero form \emph{open} sets in $\mathcal M$. On these sets
  $\Psi_2\equiv 0$ everywhere. Thus the discussed situation has only
  two possible outcomes: either $\kappa\equiv \sigma\equiv 0$ (which
  finishes the proof), or we
  have $\Psi_2\equiv 0$ in an open set, in addition to the assumed
  $\Psi_0\equiv\Psi_1\equiv 0$.
\item In this latter case we look at the Bianchi identities (\ref{b3}) and
  (\ref{b4}), obtaining: $\si\Psi_3\equiv 0$ and $\kappa\Psi_3\equiv
  0$. This again leads to either $\si\equiv\kappa\equiv 0$ or to
  $\Psi_3\equiv 0$ in addition to
  $\Psi_0\equiv\Psi_1\equiv\Psi_2\equiv 0$.
\item If $\Psi_3\equiv 0$ the Bianchi identities (\ref{b5}) and
  (\ref{b6}) give: $\si\Psi_4\equiv 0$ and $\kappa\Psi_4\equiv
  0$. Thus 
  if we want to have nonvanishing selfdual part of the Weyl tensor, we are forced
  to have $\kappa\equiv\si\equiv 0$.
\item This finishes the proof in this direction.
\end{itemize}  
Thus in going from $\Big(\Psi_0\equiv\Psi_1\equiv 0\Big)$ to
$\Big(\kappa\equiv\sigma\equiv0\Big)$, we are only able to prove the
theorem in the classical (although with a possibly nonzero cosmological
constant) Goldberg-Sachs version, namely Theorem \ref{gst}, (2).
\end{proof}
Whether it is possible to weaken the Einstein assumption above to
$Ric_{|\mathcal N}\equiv 0$ is an
open question.
\subsection{Generalizing the Kundt-Thompson and the Robinson-Schild version}\label{lpo}
As noted by Kundt and Thompson \cite{KT} and Robinson and Schild \cite{robsch}, to achieve 
the algebraic speciality of the metric, when $\kappa\equiv\sigma\equiv0$ has
been assumed, it is sufficient to use weaker
conditions than
$\Rho_{11}\equiv\Rho_{14}\equiv\Rho_{44}\equiv0$. There are various
approaches to obtain these conditions in the General Relativity 
literature (see e.g. \cite{pr}). In this section we present our
approach, which is signature independent.

We first assume that $\Rho_{11}\equiv\Rho_{14}\equiv\Rho_{44}\equiv0$ holds \emph{only
conformally}. Thus we merely assume that there exists a scale
$\Upsilon:{\mathcal M}\to\bbR$ such that the rescaled metric
$\hat{g}={\rm e}^{2\Upsilon}g$ satisfies $$\hat{Ric}_{|\mathcal N}\equiv0,$$
where ${\mathcal N}={\rm Span}_\bbC(m,k)$. This means that choosing a null
coframe $(M,P,N,K)$ for $g$, and the corresponding  rescaled
null coframe $\hat{M}={\rm e}^{\Upsilon}M$,  $\hat{P}={\rm
  e}^{\Upsilon}P$,  $\hat{N}={\rm e}^{\Upsilon}N$ and $\hat{K}={\rm
  e}^{\Upsilon}K$ for $\hat{g}$ we have 
\be
\hat{\Rho}_{11}\equiv\hat{\Rho}_{14}\equiv\hat{\Rho}_{44}\equiv0\label{cr}.\ee 
Note that for this to be satisfied we do not need to assume
$\Rho_{11}\equiv\Rho_{14}\equiv\Rho_{44}\equiv0$. Our aim now is to deduce what
restrictions on $g$ are imposed by equations (\ref{cr}).  

As it is well known (see e.g. \cite{govnur}) 
the rescaled Schouten tensor $\hat{\Rho}$ is related to $\Rho$ via:
$$\hat{\Rho}_{ab}=\Rho_{ab}-\nabla_a\Upsilon_b+\Upsilon_a\Upsilon_b-\tfrac12
\Upsilon_c\Upsilon^c g_{ab},$$
with $\Upsilon_a=\nabla_a\Upsilon$. 
Now, applying the covariant derivative 
$\nabla_c$ on both sides of this equation, 
antisymmetrizing over the indices $\{ca\}$ and using again this
equation to eliminate the covariant derivatives of $\Upsilon_a$ we
get
\be
\nabla_{[c}\hat{\Rho}_{a]b}+\Upsilon_{[a}\hat{\Rho}_{c]b}+\Upsilon^d\hat{\Rho}_{d[a}g_{c]b}\equiv\tfrac12
(A_{bca}+C_{acb}^{\quad d}\Upsilon_d).\label{bi}\ee
 Here $A_{abc}$ is the Cotton tensor 
$$A_{abc}= 2\nabla_{[b}\Rho_{c]a},$$
and $C_{abcd}$ is the Weyl tensor. Note that in addition to 
$$
A_{abc}=-A_{acb},
$$
as a consequence of the
first and the second Bianchi identities, we also have:
\be
A_{abc}+A_{cab}+A_{bca}=0,\label{420}\ee
and
\be
A_{abc}=\nabla_dC^d_{~abc},\quad\quad\quad A^a_{~ab}=0,\label{421}
\ee
respectively.

The obtained identity (\ref{bi}) is a generalization of the
identity known in the theory of conformally Einstein spaces (see
e.g. \cite{govnur}). It is interesting on its own, but it is particularly
useful in our situation of equations (\ref{cr}).   

Let as assume that in addition to (\ref{cr}) the distribution of
totally null planes $\mathcal N$ is integrable. This means that in the
frame $(m,p,n,k)$ we have
$$\kappa\equiv\sigma\equiv0,$$ 
which is the same as assuming that the respective connection
coefficients satisfy 
\be
\Gamma_{414}\equiv\Gamma_{411}\equiv 0.\label{gami}
\ee
As we proved in the previous section this 
implies that the Weyl tensor
coefficient $$\Psi_0\equiv0.$$ Now, using the frame $(m,p,n,k)$ and our
assumptions (\ref{cr}) and (\ref{gami}) on the l.h.s of the identity (\ref{bi}),  
we directly check that the following proposition is true:
\begin{proposition}\label{dodat}
Suppose that a distribution of totally null 2-planes ${\mathcal N}$ on
$({\mathcal M},g)$ be
integrable, $[{\mathcal N},{\mathcal N}]\subset {\mathcal N}$, and
that the Schouten tensor $\hat{P}$ of the rescaled metric
$\hat{g}={\rm e}^{2\Upsilon}g$ is degenerate on $\mathcal N$,
$\hat{P}_{|\mathcal N}\equiv 0$. Then 
for every three vector fields $X,Y,Z\in{\mathcal N}$ we have:
$$
X^aY^bZ^c(\nabla_{[c}\hat{\Rho}_{a]b}+\Upsilon_{[a}\hat{\Rho}_{c]b}+\Upsilon^d\hat{\Rho}_{d[a}g_{c]b})\equiv
0.
$$
\end{proposition} 
Since, in addition, in the coframe $(m,p,n,k)$ the Weyl tensor coefficient 
$\Psi_0\equiv 0$, the r.h.s. of (\ref{bi}), after being contracted with
vectors $X,Y,Z$ from $\mathcal N$, includes only the Weyl tensor
coefficient $\Psi_1$. Thus the considered identity, when restricted to
$\mathcal N$, reduces to two complex equations:
\be
A_{141}-\Psi_1\delta\Upsilon\equiv 0,\label{es1}
\ee 
and
\be
A_{441}-\Psi_1D\Upsilon\equiv 0.\label{es2}
\ee 
This relates the components $\{141\}$ and $\{441\}$ of the Cotton
tensor \emph{algebraically} to the Weyl tensor coefficient $\Psi_1$,
and proves the following 
\begin{proposition}\label{pro}
A metric $g$ with an integrable field of selfdual totally null 2-planes
$\mathcal N$ on a 4-dimensional manifold $\mathcal M$ 
admits a conformal scale $\Upsilon:{\mathcal M}\to\bbR$ such that the
rescaled metric $\hat{g}$ has Ricci tensor $\hat{Ric}$ degenerate on
$\mathcal N$,
$$\hat{Ric}_{|\mathcal N}\equiv0,$$
only if the Cotton tensor $A$ of the original metric satisfies
equations (\ref{es1}), (\ref{es2}) in a null coframe in which $\kappa\equiv\sigma\equiv0$.  
\end{proposition}
It is interesting that the expressions (\ref{es1}) and
(\ref{es2}) appear also in the following
\begin{proposition}\label{pra}
Suppose that a metric $g$ admits an integrable maximal totally null field of
2-planes. Then the Cotton tensor components $A_{141}$ and $A_{441}$ 
in the null coframe $(M,P,N,K)$ in which $\kappa\equiv\sigma\equiv0$ are related
to the Cotton tensor  components $\hat{A}_{\hat{1}\hat{4}\hat{1}}$ and $\hat{A}_{\hat{4}\hat{4}\hat{1}}$ 
in the null coframe $({\rm e}^{\Upsilon}M,{\rm e}^{\Upsilon}P,{\rm
  e}^{\Upsilon}N,{\rm e}^{\Upsilon}K)$ of the rescaled metric
$\hat{g}={\rm e}^{2\Upsilon}g$ via 
\begin{eqnarray}
\hat{A}_{\hat{1}\hat{4}\hat{1}}={\rm e}^{-3\Upsilon}(A_{141}-\Psi_1\delta\Upsilon),\label{pra1}\\
\hat{A}_{\hat{4}\hat{4}\hat{1}}={\rm e}^{-3\Upsilon}(A_{441}-\Psi_1D\Upsilon).\label{pra2}
%%/home/pawel/notebooks/stonybrook/goldbegsachs/nplorentzgspaperroddennygsweaknowe1conformalcottontrans.nb
\end{eqnarray}   
\end{proposition}
The proof of this fact is straightforward. For example it can be checked in the Newman-Penrose formalism with
$\kappa=\sigma=0$, in which the relevant components of the Cotton
tensor read:
\begin{eqnarray}
&&A_{141}=D\Rho_{11}-\bel\Rho_{14}+(2\epsilon'-2\epsilon+\rho')\Rho_{11}+(2\bet+2\pi')\Rho_{14}-\la'\Rho_{44},\label{kot1}\\
&&A_{441}=D\Rho_{14}-\bel\Rho_{44}-\kappa'\Rho_{11}+(2\rho'-2\epsilon)\Rho_{14}+(2\alpha'+2\bet+\pi')\Rho_{44}.\label{kot2}
\end{eqnarray}

Now, treating the Cotton tensor $A$ as a linear map  
${\rm T}{\mathcal M}\times{\rm T}{\mathcal M}\times{\rm T}{\mathcal
  M}\to\bbR$, we recall that $A$ is \emph{degenerate} on a vector
distribution $\mathcal Z$, $A_{|\mathcal Z}=0$,  iff $A(Z_1,Z_2,Z_3)=0$ for all
$Z_1,Z_2,Z_3\in \mathcal Z$. 
Then, if we take ${\mathcal N}={\rm
  Span}_\bbC(m,k)$, where $(m,p,n,k)$ is a null frame, we see that
$A_{441}=A_{141}=0$ if and only if $A_{|\mathcal N}=0$.   
This together with Propositions \ref{pro} and \ref{pra} imply the following 
\begin{corollary}
Suppose that a metric $g$ admits an \emph{integrable} 
maximal totally null field $\mathcal N$ of 2-planes. If the
metric can be conformally 
rescaled to $\hat{g}$ so that the rescaled Ricci tensor $\hat{Ric}$ 
is degenerate
on $\mathcal N$, $\hat{Ric}_{|\mathcal N}\equiv0$, then in
this scale the rescaled Cotton tensor $\hat{A}$ is degenerate on $\mathcal N$,
$\hat{A}_{|\mathcal N}\equiv0$.   
\end{corollary}
\begin{remark}
We note that given an integrable totally null field of 
2-planes $\mathcal N$ the condition $\hat{A}_{|\mathcal N}\equiv0$ is
\emph{weaker} than $\hat{Ric}_{|\mathcal N}\equiv0$. We saw that 
$\hat{Ric}_{|\mathcal N}\equiv0$ implies $\hat{A}_{|\mathcal N}\equiv0$, but 
the converse is not guaranteed.
%%/home/pawel/notebooks/stonybrook/goldbegsachs/nplorentzgspaperroddennygsweaknowe1conformal.nb
\end{remark}
Now we use the Bianchi identities (\ref{b11}) and (\ref{b21}), which
we display here as the following  
\begin{lemma}\label{bil}
On any 4-dimensional manifold with a metric $g$ as in (\ref{npg}) we have
\begin{eqnarray}
A_{141}&\equiv&\vel\Psi_0 +(\mu- 4 \gam) \Psi_0 -\bel\Psi_1 + 2 (2\tau+\bet) \Psi_1 - 3\si \Psi_2\label{ib1} \\
A_{414}&\equiv& \pel\Psi_0 -(\pi+ 4 \al) \Psi_0  +\del\Psi_1+ 2 (2\rho-\eps) \Psi_1 + 
   3 \kappa \Psi_2 \label{ib2}.\end{eqnarray}
\end{lemma}
\begin{proof}
This is proved in the Appendix, but we can also see this by observing
that subtracting (\ref{ib1}) from (\ref{kot1}) and, respectively
(\ref{ib2}) from (\ref{kot2}) we obtain the respective Bianchi
identities (\ref{b1}) and (\ref{b2}).   
\end{proof}
This Lemma is crucial for the rest of our arguments in this
section. It has various consequences, the first being the following 
much sharper version of part
one of Theorem \ref{gst}:
 \begin{theorem}\label{gsccc}
Let $\tn\subset {\rm T}^\bbC \mathcal M$ be a field of totally null
2-planes on a 4-dimensional manifold
$({\mathcal M},g)$ equipped with metric $g$. 
Assume that the Cotton tensor $A$ of the metric $g$, considered as a 
threelinear form
on ${\rm T}^\bbC \mathcal M$, is degenerate on $\tn$,
$$A_{|\tn}\equiv0.$$
If in addition the field $\tn$ is integrable, $[\tn,\tn]\subset\tn$, everywhere
  on $\mathcal M$, then 
$({\mathcal M},g)$ is algebraically special at every point, with a
 
  field of multiple principal totally null 2-planes tangent to $\tn$.
\end{theorem}
\begin{proof}
In an adapted null coframe $(M,P,N,K)$ our integrability assupmtion is 
$\kappa\equiv\sigma\equiv 0$, which  as we know, implies $\Psi_0\equiv
0$. The assumption about the degeneracy of the Cotton tensor means
$A_{141}\equiv A_{441}\equiv
0$, which together with $\Psi_0\equiv 0$ and Lemma \ref{bil} gives the
identities: $\bel\Psi_1 \equiv 2 (\bet  + 2 \tau)\Psi_1$ and 
$\del\Psi_1 \equiv 2 (\eps  - 2\rho) \Psi_1$. This implies
$\Psi_1\equiv 0$ by Lemma \ref{gstl}. Thus the field of (principal) 
totally null 2-planes $\mathcal N$ is multiple.
\end{proof}
\begin{remark}
Note that as a result of this theorem, the assumption 
$A_{|\mathcal N}\equiv 0$ is \emph{conformal}. Without knowing that
$\kappa\equiv\sigma\equiv 0$ and $A_{141}\equiv A_{441}\equiv 0$ imply
$\Psi_1\equiv 0$, the assumption $A_{141}\equiv A_{441}\equiv 0$ seemed to
be \emph{not} conformal, because of the inhomogeneous terms in the
transformations (\ref{pra1})-(\ref{pra2}). But since under the
assumptions 
$\kappa\equiv\sigma\equiv 0$ and $A_{141}\equiv A_{441}\equiv 0$ we
were able to discover that actually $\Psi_1\equiv 0$, then $A_{141}$
and $A_{441}$ transform homogeneously under the conformal rescaling. 
Thus in such case the
condition $A_{|\mathcal N}\equiv 0$ \emph{is} conformal. 
\end{remark}
The second application of Lemma \ref{bil} is included in the following 
\begin{remark}\label{krem}
Suppose that we would like to have a still sharper (than in Theorem
\ref{gsccc}) version of part one of Theorem \ref{gst}. Thus
instead of assuming $Ric_{|\mathcal N}\equiv 0$, or the weaker
condition 
$A_{|\mathcal N}\equiv 0$, we would like to have an assumption about
vanishing of still higher order derivatives of the curvature, that
together with $\kappa\equiv\sigma\equiv 0$ would imply $\Psi_1\equiv
0$. Then Lemma (\ref{bil}) assures that it is impossible, and the
condition $A_{|\mathcal N}\equiv 0$ can not be weakened. Indeed,
denoting such hypothetical condition by $S\equiv 0$, we would have
$(\kappa\equiv\sigma\equiv 0~\&~S\equiv 0)\Rightarrow (\Psi_1\equiv
0)$. But since $\kappa\equiv\sigma\equiv 0$, in addition, implies that
$\Psi_0\equiv 0$, then Lemma \ref{bil} implies $A_{141}\equiv
A_{441}\equiv 0$. Thus the hypothethically weaker than $A_{|\mathcal
  N}\equiv 0$ condition $S\equiv
0$, in turn, implies $A_{|\mathcal N}\equiv 0$. Since this alone, according to Theorem
\ref{gsccc}, is already sufficient to imply $\Psi_1\equiv 0$, we do not
need condition $S\equiv 0$ to obtain the desired result. This
proves the following      
\end{remark}
\begin{theorem}\label{west}
The weakest curvature condition which together with the integrability
condition $[{\mathcal N},{\mathcal N}]\subset\mathcal N$, implies that
the field of totally null 2-planes $\mathcal N$ is principal and
multiple is the degeneracy of the Cotton tensor on $\mathcal N$,
$A_{|\mathcal N}\equiv 0$. 
\end{theorem}
\begin{example}\label{pk}
An example of a condition $S\equiv 0$ which is \emph{a priori} 
weaker than $A_{|\mathcal N}\equiv
0$  may be obtained as follows. The 
procedure used in the proof of Lemma \ref{gstl} may be
equally applied to
the situation in which the conditions (\ref{pr1})-(\ref{pr2}) are replaced
by the Bianchi identities (\ref{ib1}) and (\ref{ib2}). Then, under the
assumption that $\kappa\equiv\sigma\equiv 0$, and hence
$\Psi_0\equiv 0$, we literally repeat all the steps from 
the proof of Lemma \ref{gstl}. Indeed, starting with 
the application of $\bel$ on both sides of $A_{141}\equiv-\bel\Psi_1
+ 2 (2\tau+\bet) \Psi_1$ and $\del$ on
both sides of $A_{414}\equiv\del\Psi_1+ 2 (2\rho-\eps) \Psi_1$, 
after subtraction and use of the commutator 
(\ref{comw}), we obtain the following \emph{identity}:
\begin{eqnarray}
&&\quad\quad\quad\quad\quad\quad\quad\quad\quad\quad\quad\quad-10\Psi_1^2\equiv\label{ideos}\\
&&\del A_{141}-\bel
A_{441}-(3\epsilon-\rho'-\epsilon'-4\rho)A_{141}+(3\bet+\alpha'+\pi'+4\tau)A_{441}.\nonumber
\end{eqnarray}
This, is satisfied \emph{always} when $\kappa\equiv\sigma\equiv
0$. Thus \emph{the vanishing of the r.h.s of}
(\ref{ideos}) implies $\Psi_1\equiv 0$. 
Moreover, since when $\kappa\equiv\sigma\equiv 0$ the
vanishing of $\Psi_1$ is a conformal property, then the vanishing of
the r.h.s. of (\ref{ideos}) is a \emph{conformal} property. In fact a direct calculation shows that if 
in a null coframe $(M,P,N,K)$ we have $\kappa\equiv\sigma\equiv0$ and  
\be
S=\del A_{141}-\bel
A_{441}-(3\epsilon-\rho'-\epsilon'-4\rho)A_{141}+(3\bet+\alpha'+\pi'+4\tau)A_{441}\label{s}\ee 
then in the conformally rescaled metric $\hat{g}={\rm e}^{2\Upsilon}g$
and in the corresponding rescaled null 
coframe  $({\rm e}^{\Upsilon}M,{\rm e}^{\Upsilon}P,{\rm
  e}^{\Upsilon}N,{\rm e}^{\Upsilon}K)$ we have
$\hat{\kappa}\equiv\hat{\sigma}\equiv 0$ and 
$$\hat{S}={\rm e}^{-4\Upsilon}S.$$
Now using the
explicit formulae for the covariant derivatives of the Cotton tensor
components $A_{141}$ and $A_{441}$: 
\begin{eqnarray*}
&&\nabla_4 A_{141}=\del A_{141}-(3\epsilon-\epsilon')A_{141}+\pi'A_{441}\\
&&\nabla_1 A_{441}=\bel A_{141}-(3\bet+\alpha')A_{441}-\rho'A_{141}, 
\end{eqnarray*}
solving this for $\del A_{141}$ and $\bel A_{141}$ and inserting in
(\ref{s}), we get 
\be
S=\nabla_4 A_{141}-\nabla_1
A_{441}+4\rho A_{141}+4\tau A_{441}.\label{s'}\ee
We thus have a condition $S\equiv 0$, which together with
$\kappa\equiv\sigma\equiv 0$ is \emph{conformal} and \emph{implies}
that $\Psi_1\equiv 0$. It is always satisfied when $A_{441}\equiv
A_{141}\equiv 0$, i.e. we have $(A_{441}\equiv
A_{141}\equiv 0)\Rightarrow S\equiv 0$, and at the first glance there
is no reason for the implication  $(S\equiv 0)\Rightarrow (A_{441}\equiv
A_{141}\equiv 0)$ However, this implication is true, on the ground of 
the discussion in Remark \ref{krem}. As a consequence we have
\begin{proposition}
Under the assumption that the distribution of selfdual totally null
2-planes $\mathcal N$ is integrable, $[\nen,\nen]\subset \nen$, the
following two, conformally invariant, conditions are equivalent
\begin{itemize}
\item the Cotton tensor of the metric $g$ is degenerate on $\nen$,
  $A_{|\nen}\equiv 0$
\item the scalar $S$ of the metric $g$, as defined in (\ref{s'}),
  identically vanishes, $S\equiv 0$.
\end{itemize} 
\end{proposition}
\end{example}

To discuss the next application of Lemma \ref{bil} we introduce
\begin{definition}\label{lop}
A metric $g$ on a 4-dimensional manifold $\mathcal M$ is called
\emph{II-generic} if and only if the points in which its
selfdual part of the Weyl 
tensor degenerates to Petrov types
\emph{III}, \emph{N} or 0 are rare, in the sense that they belong to closed sets without interior in
$\mathcal M$.
\end{definition}

In particular every metric with 
selfdual part of the Weyl tensor being at each point of $\mathcal M$ algebraically
general, or of mixed type: algebraically general on some subsets and type
II or type D on their complements, is II-generic; a metric which is
e.g. of type III in an open set of $\mathcal M$ is not II-generic.

Now we are ready to discuss a slight generalization of the known
\emph{conformal} versions of the Goldberg-Sachs theorem. In the Lorentzian
case such versions were given by Kundt and Thompson \cite{KT} and Robinson and Schild
\cite{robsch}. Penrose and Rindler \cite{pr} gave a complex
(spinorial) version of
the Kundt-Thompson/Robinson-Schild theorem. Here we quote \emph{our} complex version, which is a slight
generalization: 

\begin{theorem}\label{pol}
%Let $\mathcal M$ be a 4-dimensional manifold with a metric $g$ whose
%selfdual Weyl tensor is at each point of $\mathcal M$ 
%in the complement of the algebraic (Petrov) types \emph{III}, \emph{N} or
%\emph{0}. Let $\mathcal N$ be a field of totally null 2-planes on $\mathcal
%M$.
Let $\mathcal M$ be a 4-dimensional manifold with a 
II-generic metric $g$. Let $\mathcal N$ be a field of selfdual totally null 2-planes on $\mathcal
M$. Then any two of the following imply the third:
\begin{itemize}
\item[(0)] The Cotton tensor of $g$ is degenerate on ${\mathcal N}$,
  $A_{|\mathcal N}\equiv 0$.    
\item[(i)] $\mathcal N$ is integrable, $[{\mathcal N},{\mathcal
    N}]\subset\mathcal N$.
\item[(ii)] The selfdual part of the Weyl tensor is algebraically special on
  $\mathcal M$ with $\mathcal N$ being a multiple principal field of
  selfdual totally null 2-planes.
\end{itemize}
\end{theorem}

\begin{proof}
First we observe that the implication $\Big({\rm (0)}~\&~{\rm
  (i)}\Big)\Rightarrow {\rm (ii)}$ is true, as a simple application
of Theorem \ref{gsccc}. Note that for this we do not need the 
genericity assumption about the Weyl tensor.

To prove the other two implications we choose a null coframe on 
$({\mathcal M},g)$ so that 
${\mathcal N}={\rm Span}_\bbC(m,k)$ and $g=2(MP+NK)$ as in
(\ref{npg}). 
Then \begin{itemize}
\item the condition (0) is: $A_{141}\equiv A_{441}\equiv 0$, 
\item the
condition (i) is:
$\kappa\equiv\sigma\equiv 0$, 
\item the condition   
(ii) is: $\Psi_0\equiv\Psi_1\equiv 0$.\end{itemize} 
Now, the proof of $\Big({\rm (i)}~\&~{\rm
  (ii)}\Big)\Rightarrow {\rm (0)}$ is an immediate consequence of
Lemma \ref{bil}, since the assumptions ${\rm (i)}~\&~{\rm (ii)}$
imply the identical vanishing of the r.h.s. of identities 
(\ref{ib1})-(\ref{ib2}), which means that also their
l.h.s. identically vanish, $A_{141}\equiv A_{441}\equiv 0$. Note that 
also in the proof of this statement the genericity assumption
about the Weyl tensor was not needed. 

This assumption is however needed to get the
last implication  $\Big({\rm (0)}~\&~{\rm
  (ii)}\Big)\Rightarrow {\rm (i)}$. Indeed assuming 
${\rm (i)}~\&~{\rm (ii)}$, the identities (\ref{ib1})-(\ref{ib2}) from
Lemma \ref{bil} reduce to the identities $-3\sigma \Psi_2\equiv 0$
and $3\kappa\Psi_2\equiv0$. Now, similarly as in the proof of part two
of the Theorem \ref{gst}, to conclude that $\kappa\equiv \sigma\equiv
0$ in a neighbourhood ${\mathcal U}\subset{\mathcal M}$, 
it is enough to assume that $\Psi_2\neq 0$ on the complement of the
closed sets without interior in $\mathcal U$. Since in our
coframe in $\mathcal U$, according to (ii), we have $\Psi_0\equiv\Psi_1\equiv
0$, Proposition \ref{pett1} assures that the coefficient $\Psi_2$ of the Weyl tensor is nonvanishing on the complement of the
closed sets without interior in $\mathcal U$ if and
only if the metric is II-generic in $\mathcal U$. Since this
is the main aasumption of Theorem \ref{pol} we 
we see that $3\sigma \Psi_2\equiv 0$
and $3\kappa\Psi_2\equiv0$ imply $\kappa\equiv \sigma \equiv 0$ in $\mathcal U$. This  
proves the part $\Big({\rm (0)}~\&~{\rm
  (ii)}\Big)\Rightarrow {\rm (i)}$ of the theorem. 
\end{proof} 
As a consequence of this proof we also have the
following 
\begin{corollary}\label{coco}
Let $\mathcal M$ be a 4-dimensional manifold with a metric $g$ and 
let $\mathcal N$ be a field of selfdual totally null 2-planes on $\mathcal
M$.  Assume that $\mathcal N$ is integrable, $[{\mathcal N},{\mathcal
    N}]\subset\mathcal N$, and that 
the selfdual part of the Weyl tensor is algebraically special on
  $\mathcal M$, with $\mathcal N$ being a multiple principal field of
  selfdual totally 2-planes.
Then the Cotton tensor of $g$ is degenerate on ${\mathcal N}$,
  $A_{|\mathcal N}\equiv 0$.    
\end{corollary}

To discuss the sharpening of the Theorem \ref{pol} with respect to the
implication $\Big({\rm (0)}~\&~{\rm
  (ii)}\Big)\Rightarrow {\rm (i)}$ we introduce two more notions
analogous to the II-generiticity. 

\begin{definition}
A metric $g$ on a 4-dimensional manifold $\mathcal M$ is called
\emph{III-generic} if and only if the points in which its 
selfdual part of the Weyl tensor degenerates to Petrov types \emph{N} or 0 belong
to closed sets without interior in
$\mathcal M$. Similarly, a metric $g$ on a 4-dimensional manifold $\mathcal M$ is called
\emph{N-generic} if and only if the points in which its selfdual part
of the Weyl tensor 
vanishes belong to closed sets without interior in
$\mathcal M$. 
\end{definition}  

For the III-generic metrics we have the following 
\begin{theorem}\label{pol1}
Let $\mathcal M$ be a 4-dimensional manifold with a III-generic metric
$g$, whose selfdual part of the Weyl tensor is in addition algebraically special
at all points of $\mathcal M$. 
Let $\mathcal N$ be the corresponding field of 
multiple principal totally null 2-planes on $\mathcal M$. 
If the Cotton tensor $A$ of the metric $g$ satisfies 
$$A(~\cdot~,Z_1,Z_2)\equiv 0,\quad\quad\quad\forall Z_1,Z_2\in\mathcal N$$
then the field $\mathcal N$ is integrable, $[{\mathcal N},{\mathcal
    N}]\subset\mathcal N$, on $\mathcal M$.
\end{theorem}
Similarly for the N-generic metrics we have
\begin{theorem}\label{pol2}
Let $\mathcal M$ be a 4-dimensional manifold with an N-generic metric
$g$, whose selfdual part of the Weyl tensor is in addition algebraically special
at all points of $\mathcal M$. 
Let $\mathcal N$ be the corresponding field of 
multiple principal totally null 2-planes on $\mathcal M$. 
Consider the 2-forms $A_Z=A(Z,~\cdot~,~\cdot~)$, where $A$ is the
Cotton tensor 
of the metric $g$ and $Z$ is a complex-valued vector field $Z$ on $\mathcal M$.
If for every vector field $Z\in\mathcal N$ the two form $A_Z$ 
is antiselfdual at each point of $\mathcal M$, then the field $\mathcal N$ is integrable, $[{\mathcal N},{\mathcal
    N}]\subset\mathcal N$, on $\mathcal M$.
\end{theorem}
We first prove Theorem \ref{pol1}.
\begin{proof}
Again we choose a null coframe on 
$({\mathcal M},g)$ so that 
${\mathcal N}={\rm Span}_\bbC(m,k)$ and $g=2(MP+NK)$ as in
(\ref{npg}). Since ${\mathcal N}$ consists of \emph{multiple}
principal null 2-planes, according to Proposition \ref{pett1}, we have
$\Psi_0\equiv\Psi_1\equiv 0$ in this coframe. Moreover, 
in this coframe the condition $A(~\cdot~,Z_1,Z_2)\equiv
0$ $\forall Z_1,Z_2\in\mathcal N$ means that the coframe components
$A_{i41}$, $i=1,2,3,4$ satisfy
\be
A_{141}\equiv A_{214}\equiv A_{341}\equiv A_{414}\equiv 0.\label{sas}\ee 
Now we again use the Bianchi identities (\ref{b11})-(\ref{b21}) which
reduce to 
$$-3\sigma\Psi_2\equiv 0,\quad\quad{\rm and}\quad\quad
3\kappa\Psi_2\equiv 0.$$

Similarly as in the proof of the second part of the 
Theoerm \ref{gst} this yields $\kappa\equiv\sigma\equiv 0$, with the
exception when $\Psi_2\equiv 0$. In such a case we have 
\be\Psi_0\equiv\Psi_1\equiv\Psi_2\equiv 0,\label{ssa}\ee
and these two Bianchi
identities are tautologies. Thus to conclude something about $\kappa$
and $\sigma$ we need to use another pair of Bianchi identities. These
are given by (\ref{b31})-(\ref{b41}) and refer to the respective 
components $A_{341}$ and $_{214}$ of the Cotton tensor. Now, with the
assumed (\ref{sas}) and (\ref{ssa}) these identities reduce to
$$2\sigma \Psi_3\equiv 0,\quad\quad{\rm
  and}\quad\quad\-2\kappa\Psi_3\equiv 0.$$ 
This does not yield $\kappa\equiv\sigma\equiv 0$ only if $\Psi_3\equiv
0$ in the neighbourhood. But this is forbidden by our assumption that
the metric is III-generic in the considered neighbourhood. 

Thus if the metric is III-generic in the neighbourhood we proved that
$\kappa\equiv\sigma\equiv 0$ in a frame adapted to $\mathcal N$, which
according to Proposition \ref{petu}, means that $\mathcal N$ is
integrable.  
\end{proof}

\begin{proof} \emph{of Theorem \ref{pol2}}. Choosing the null frame as
  in the above proof we first interpret the condition about the Cotton
  tensor 2-forms $A_Z$ being all antiselfdual. Since $\mathcal N$ is
  spanned by $m$ and $k$ we only need to consider
  the 2-forms $A_m=A(m,~\cdot~,~\cdot~)$ and 
$A_k=A(k,~\cdot~,~\cdot~)$. We have:
\begin{eqnarray*}
A_m=A_{123}\theta^2\dz\theta^3+\tfrac12(A_{112}-A_{134})(\theta^1\dz\theta^2-\theta^3\dz\theta^4)+A_{114}\theta^1\dz\theta^4+\\
A_{113}\theta^1\dz\theta^3+\tfrac12(A_{112}+A_{134})(\theta^1\dz\theta^2+\theta^3\dz\theta^4)+A_{124}\theta^2\dz\theta^4
\end{eqnarray*}  
and 
\begin{eqnarray*}
A_k=A_{423}\theta^2\dz\theta^3+\tfrac12(A_{412}-A_{434})(\theta^1\dz\theta^2-\theta^3\dz\theta^4)+A_{414}\theta^1\dz\theta^4+\\
A_{413}\theta^1\dz\theta^3+\tfrac12(A_{412}+A_{434})(\theta^1\dz\theta^2+\theta^3\dz\theta^4)+A_{424}\theta^2\dz\theta^4. 
\end{eqnarray*} 
So  looking at the bases
  (\ref{sela}) and (\ref{asel}) of the
  selfdual and antiselfdual 2-forms $\Sigma$ and $\Sigma'$, we
  conclude that these 2-forms are antiselfdual iff the following six
  conditions for the coframe components of the Cotton tensor are 
satisfied:
\begin{eqnarray}
&&A_{114}\equiv
  A_{414}\equiv 0\quad\&\label{uh1}\\
&&A_{112}- A_{134}\equiv 0\quad\&\label{uhh}\\
&&A_{412}-
  A_{434}\equiv 0\quad\&\label{uh}\\
&&A_{123}\equiv A_{423}\equiv 0.\label{uh2}
\end{eqnarray}
Now we use the symmetries of the Cotton tensor to give equivalent
forms of the conditions (\ref{uhh})-(\ref{uh}). Using (\ref{421}) we get
$A_{112}\equiv A_{341}-A_{413}$ and using (\ref{420}) we get
$A_{134}\equiv -A_{413}-A_{341}$. Subtracting the latter from the
former we get the identity $$A_{112}-A_{134}\equiv 2A_{341}.$$ In the
similar way we prove the identity $$A_{412}-A_{434}\equiv 2A_{214}.$$ 
Comparing these two identities with (\ref{uh1})-(\ref{uh2}) we
conclude that the condition that $A_Z$ is antiselfdual for all
$Z\in\mathcal N$, in our coframe, is equivalent to the six conditions
\begin{eqnarray*}
&&A_{114}\equiv A_{414}\equiv 0\quad\&\\
&&A_{341}\equiv A_{214}\equiv0\quad\&\\
&&A_{123}\equiv A_{423}\equiv 0.
\end{eqnarray*} 
Since
the first four conditions are precisely $A(~\cdot~,Z_1,Z_2)\equiv 0$
for $Z_1,Z_2\in\mathcal N$, we now use Theorem \ref{pol1} to conclude
$\kappa\equiv\sigma\equiv 0$, provided that we are not in the situation
when 
\be\label{dfd}\Psi_0\equiv\Psi_1\equiv\Psi_2\equiv\Psi_3\equiv 0\ee
in the neighbourhod. If this is the case, to show that we still have
$\kappa\equiv\sigma\equiv 0$ we need the additional assumption
(\ref{uh2}). With this and (\ref{dfd}) being assumed, using the Bianchi identities
(\ref{b51})-(\ref{b61}), we easilly obtain 
$$-\sigma\Psi_4\equiv 0\quad\quad{\rm and}\quad\quad\kappa\Psi_4\equiv
0.$$
This implies that $\kappa\equiv\sigma\equiv 0$ in the neighbourhood,
on the ground of the N-genericity of the metric. This
finishes the proof.
\end{proof}  
As a counterpart to Corollary \ref{coco} we have
\begin{corollary}
Let $\mathcal M$ be a 4-dimensional manifold with a metric $g$ and 
let $\mathcal N$ be a field of selfdual totally null 2-planes on $\mathcal
M$.  Assume that $\mathcal N$ is integrable, $[{\mathcal N},{\mathcal
    N}]\subset\mathcal N$, and that 
the selfdual part of the Weyl tensor is algebraically special on
  $\mathcal M$ with $\mathcal N$ being a multiple principal field of
  selfdual totally null 2-planes.
Then if $\mathcal N$ has multiplicity equal to \emph{three} 
the Cotton tensor of $g$ satsifies $A(~\cdot~,Z_1,Z_2) \equiv 0$ for
all $Z_1,Z_2\in\mathcal N$. If $\mathcal N$ has multiplicity 
equal to \emph{four} the 2-form $A_Z$ of the Cotton tensor $A$ 
of $g$ is antiselfdual.   
\end{corollary} 
\begin{proof}
The proof is an immediate application of the Bianci identities (\ref{b11})-(\ref{b61}).
\end{proof}
\section{Interpretation in terms of a characteristic connection}\label{cchh}
The terms $4\rho A_{141}+4\tau A_{441}$ that appear in formula
(\ref{s'})
defining $S$ in Example \ref{pk} suggests that to describe 
the geometry of manifolds with
$\kappa\equiv\sigma\equiv 0$ it would be useful to have a 
\emph{vectorial} object, say $B_a$, with components $B_a$ being
roughly 
\be B_{a}=(B_1,B_2,B_3,B_4)=(4s^{-1}\tau,B_2,B_3,-4s^{-1}\rho),\label{asde}\ee where $s$ is a complex
constant. If we were able to find a geometric way of distinguishing such
$B_a$, then the formula for $S$ would be $S=(\nabla_4-sB_4)
A_{141}-(\nabla_1-sB_1)A_{441}$ and would have an explicit geometric
meaning. Note that the values of components $B_2$ and $B_3$
are totally irrelevant here! In this section we show how to
geometrically distinguish such (partially determined) $B_a$. 
\subsection{Characteristic connection of a totally null 2-plane}
Let us chose an arbitrary 1-form $B=B_a\theta^a$ on $({\mathcal
  M},g=g_{ab}\theta^a\theta^b)$. Given a choice of $B$ one defines a
new connection $\nw$ on $\mathcal M$, which is related to the 
Levi-Civita connection as follows. 

Let $\Gamma_{ab}=\Gamma_{abc}\theta^c$,
be the Levi-Civita connection 1-forms as given in (\ref{ca1}). Define
\be
\gw_{abc}=\Gamma_{abc}+\tfrac12
(g_{ca}B_b-g_{cb}B_a+g_{ab}B_c).\label{we4}\ee
Then the new connection $\nw$ is defined on $\mathcal M$ by
\be
\nw_Xe_b=X^c\nw_ce_b=X^c\gw^a_{~bc}e_a,\quad\quad
\gw^a_{~bc}=g^{ad}\gw_{dbc},\label{we5}\ee
where $(e_a)$ is a frame dual to the coframe $(\theta^a)$,
$e_a\hook\theta^b=\delta^b_{~a}$.

The connection $\nw$ is called the \emph{Weyl connection}. It is the unique
\emph{torsionless} connection satisfying 
\be
\nw g=-B g.\label{wg}
\ee
It has the nice property of being \emph{conformal} in the sense that if
the metric $g$ undergoes a transformation $g\to\hat{g}={\rm
  e}^{2\phi}g$, then equation (\ref{wg}) is preserved, 
$$\nw\hat{g}=-\hat{B}\hat{g},$$ with a mere
change $B\to\hat{B}=B-2\der\phi$. 
 
The conformal properties of Weyl connections would be very interesting
for our purpose of describing conformal conditions for the
Goldberg-Sachs theorem, provided that, we
were able to associate a unique Weyl form $B$ with the main object of this
theorem namely a field of totally null 2-planes $\mathcal N$. The
following theorem shows that although such a natural way of chosing $B$ is
possible only \emph{partially}, it nevertheless enables us to define a
canonical connection on $\mathcal N$, which encodes its conformal
properties.

\begin{theorem}\label{cch1}
Let $\mathcal N$ be a field of totally null 2-planes on $({\mathcal
  M},g)$, where $g$ is a 4-dimensional metric of \emph{any} (including
complex) signature. Let us assume that $\mathcal N$ is integrable
$[{\mathcal N},{\mathcal N}]\subset\mathcal N$. Then there exists a \emph{unique} connection $\check{\nabla}$ 
on $\mathcal N$,
which encodes the conformal properties of this field of totally null 2-planes. 
\end{theorem}

\begin{proof}
We define the conection $\check{\nabla}$ in two steps.

{\bf Step One}.  We first look for a Weyl connection $\nw$ on
$\mathcal M$, as in (\ref{we4})-(\ref{we5}), which has the property
that it preserves $\mathcal N$. This means that we ask if there exists
a Weyl connection $\nw$ on $\mathcal N$, such that 
\be
\nw_Y X \in {\mathcal N} \quad\quad \forall ~X\in {\mathcal N}\quad
\&\quad\forall ~Y\in {\rm T}{\mathcal M}~?\label{qey}\ee
To answer this question, we work in the adapted null frame
$(e_1,e_2,e_3,e_4)=(m,p,n,k)$, with the usual dual coframe
$(\theta^1,\theta^2,\theta^3,\theta^4)=(M,P,N,K)$, so that the field
of totally null 2-planes $\mathcal N$ is ${\mathcal
  N}=\Span(e_1,e_4)=\Span(m,k)$. Then the question (\ref{qey}) is
equivalent to the question of existence of $\nw$ such that
$$(\nw_{c}e_1)\dz e_1\dz e_4=0,\quad\&\quad (\nw_{c}e_4)\dz
e_1\dz e_4=0,\quad \forall~c=1,2,3,4,$$
where we abreviated $\nw_{e_c}$ to  $\nw_{e_c}=\nw_c$.         
It is very easy to see that, since in the chosen frame the coefficients of the metric $g_{ab}$ are
all zero, except $g_{12}=g_{21}=g_{34}=g_{43}=1$, then these
conditions are equivalent to:
$$
\begin{aligned}
&(\nw_ce_1)\dz e_1\dz e_4=\gw_{11c}e_2\dz e_1\dz
  e_4+\gw_{41c}e_3\dz e_1\dz e_4=0\\
&(\nw_ce_4)\dz e_1\dz e_4=\gw_{14c}e_2\dz e_1\dz
  e_4+\gw_{44c}e_3\dz e_1\dz e_4=0
\end{aligned}\quad\forall ~c=1,2,3,4
$$
or, what is the same,
$$\gw_{11c}=\gw_{14c}=\gw_{44c}=0,\quad \forall c=1,2,3,4.$$
Comparing these last equations with (\ref{we4}), we easily see that
$\gw_{11c}=\gw_{44c}=0$ is \emph{automatically} satisfied for all $c=1,2,3,4$, and then, by considering the remaining conditions 
$\gw_{14c}=\gw_{41c}=0$, we see that (\ref{qey}) is equivalent to:
\be
\Gamma_{14c}+\tfrac12(g_{c1}B_4-g_{c4}B_1)=0\quad \forall~ c=1,2,3,4.\label{ret}\ee
Now examining these equations for $c=1$ and $c=4$ we get the
conditions that the Levi-Civita connection coefficients $\Gamma_{141}$
and $\Gamma_{144}$ must satisfy
\be
\Gamma_{141}=\Gamma_{144}=0.\label{ret1}\ee 
Examining the equations (\ref{ret}) for $c=2$ and $c=3$, we get the
relations between the components $B_1$ and $B_4$ of the 1-form $B$ and
the Levi-Civita connection coefficients $\Gamma_{143}$ and
$\Gamma_{142}$. These are:
\be
B_1=2\Gamma_{143},\quad\quad B_4=-2\Gamma_{142}.\label{ret2}\ee
Thus, the requirement that there is a Weyl connection preserving
$\mathcal N$ is equivalent to the fact that in a coframe adapted to
$\mathcal N$, we have (\ref{ret1}) and (\ref{ret2}). Since $\Gamma_{141}$ and
$\Gamma_{144}$, in the coframe adapted to $\mathcal N$, are 
$\Gamma_{141}=\sigma$ and $\Gamma_{144}=\kappa$, then we see that the
connection $\nw$ exists only if the field of totally null 2-planes
$\mathcal N$ is \emph{integrable}. When $\mathcal N$ is integrable then, in the adapted coframe $(\theta^i)$, the
two of the components of the Weyl 1-form $B$, namely $B_1$ and $B_4$,
are totally determined. They are equal to 
$$B_1=2\tau,\quad\quad B_4=-2\rho,$$
as desired in (\ref{asde}), with $s=2$. 

Concluding this part of the proof, we say that the condition
(\ref{qey}) that the Weyl connection preserves $\mathcal N$ determines
this connection only up to the terms $B_2$ and $B_3$ in the Weyl
1-form. In step two of the proof we restrict this connection to
$\mathcal N$.

{\bf Step two}. Since $\nw$ preserves $\mathcal N$ in \emph{any}
direction then, in particular, it preserves it along $\mathcal N$.
Thus $\nw$, with \emph{any} choice of $B_2$ and $B_3$, restricts
naturally to $\mathcal N$. But \emph{apriori} this restriction may
depend on the choice of $B_2$ and $B_3$. That this is \emph{not} the
case follows from the following.

First observe that because of (\ref{ret1}), we 
have 
$$
\begin{aligned}
&\gw_{211}=\Gamma_{211}+B_1,\quad \gw_{111}=0,\quad
  \gw_{411}=0,\quad\gw_{311}=\Gamma_{311}\\
&\gw_{214}=\Gamma_{214}+\tfrac12 B_4,\quad \gw_{114}=0,\quad
  \gw_{414}=0,\quad\gw_{314}=\Gamma_{314}+\tfrac12 B_1\\
&\gw_{241}=\Gamma_{241}+\tfrac12 B_4 ,\quad \gw_{141}=0,\quad
  \gw_{441}=0,\quad\gw_{341}=\Gamma_{341}+\tfrac12 B_1\\
&\gw_{244}=\Gamma_{244},\quad \gw_{144}=0,\quad \gw_{444}=0,\quad\gw_{344}=\Gamma_{344}+B_4.
\end{aligned}
$$
Thus the covariant derivatives
$$
\begin{aligned}
&\nw_1
  e_1=\gw^c_{11}e_c=\gw_{211}e_1+\gw_{111}e_2+\gw_{411}e_3+\gw_{311}e_4\\
&\nw_4
  e_1=\gw^c_{14}e_c=\gw_{214}e_1+\gw_{114}e_2+\gw_{414}e_3+\gw_{314}e_4\\
&\nw_1
  e_4=\gw^c_{41}e_c=\gw_{241}e_1+\gw_{141}e_2+\gw_{441}e_3+\gw_{341}e_4\\
&\nw_4
  e_4=\gw^c_{44}e_c=\gw_{244}e_1+\gw_{144}e_2+\gw_{444}e_3+\gw_{344}e_4
\end{aligned}$$
of vectors $(e_1,e_4)$ in the directions $e_1$ and $e_4$ spanning
$\mathcal N$, are expressible purely in terms of the Levi-Civita 
connection coeffiecients $\Gamma_{abc}$ and the totally determined
part of $B$. In these relations the unknown coefficients of $B$,
namely $B_2$ and $B_3$, do not appear!

Thus $\nw$ restricts to a \emph{unique} and totally determined
connection on $\mathcal N$. We define
$$\check{\nabla}=\nw_{|{\mathcal N}}\quad\quad{\rm
  on}\quad\quad{\mathcal N}.$$
Since this connection is constructed with only conformal objects, it
is manifestly conformal.

The formulae for this cconnection in the Newman-Penrose formalism are:
\be\begin{aligned}
&\check{\nabla}_mm=(\bet-\al'+2\tau)m-\la'k\\
&\check{\nabla}_km=(\varepsilon-\varepsilon'-\rho)m+(\tau-\pi')k\\
&\check{\nabla}_mk=(\rho'-\rho)m+(\al'+\bet+\tau) k\\
&\check{\nabla}_kk=\kappa'm+(\varepsilon+\varepsilon'-2\rho)k.
\end{aligned}\label{chc}\ee
\end{proof}  
The connection $\check{\nabla}$ defined in Theorem \ref{cch1} is called the 
\emph{characteristic connection} of an integrable totally null 2-plane $\mathcal N$ field.

Now, having any three (complex-valued) vector fields $X,Y,Z\in\mathcal
N$, we define the torsion $\check{T}$ and the curvature $\check{R}$ of $\check{\nabla}$ via
the usual:
\begin{eqnarray}
\check{T}(X,Y)=\check{\nabla}_XY-\check{\nabla}_YX-[X,Y],\label{tut}\\
\check{R}(X,Y)Z=
[\check{\nabla}_X,\check{\nabla}_Y]Z-\check{\nabla}_{[X,Y]}Z.\label{cuc}
\end{eqnarray}
By construction these are \emph{conformal} tensors defined on $\mathcal N$. Since both $\check{T}$ and
$\check{R}$ are antisymmetric in $X,Y$ they may have at most two,
respectively four, independent components. Actually we have the
following 
%/home/pawel/notebooks/stonybrook/goldbergsachs/malakoneksja.nb
\begin{theorem}
The characteristic connection
$\check{\nabla}$ of an integrable $\mathcal N$ is \emph{torsionless}, 
$$\check{T}\equiv 0.$$ Its curvature, $\check{R}$, is given by
\begin{eqnarray}
\check{R}(m,k)m&=&4\Psi_1 m,\label{cry1}\\
\check{R}(m,k)k&=&4\Psi_1 k,\label{cry2}
\end{eqnarray}  
where $\Psi_1$ is the Weyl tensor coefficient
of the Levi-Civita connection as defined in (\ref{ca22}).
\end{theorem}
\begin{proof}
The torsionless property of the connection and the formulae
(\ref{cry1})-(\ref{cry2}) for the curvature can be
checked by a 
direct calculation. Indeed, for the torsionless we only have to show that
$\check{T}(m,k)=0$. One checks that this is a direct consequence of
the definitions (\ref{tut}), (\ref{chc}) and the commutation relation
$[\bel,\del]$ from (\ref{comw}). To check (\ref{cry1}) one
uses the definition (\ref{cuc}), the commutator $[\bel,\del]$ and the Newman-Penrose equations
(\ref{np2}), (\ref{np4}), (\ref{np5}), (\ref{np6}) and
(\ref{np7}). Similarly, to check (\ref{cry2}) one
uses (\ref{cuc}), (\ref{comw}) and the Newman-Penrose equations
(\ref{np2}), (\ref{np3}), (\ref{np6}), (\ref{np9}) and (\ref{np8}). In
all of these expressions one has to put the integrability conditions
$\kappa\equiv\sigma\equiv 0$. The
rest of the proof is easy pure algebra. 
\end{proof}   
Thus we see that the curvature of $\check{\nabla}$ has only \emph{one}
independent component, which is a constant multiple of
$\Psi_1$. Moreover, the entire curvature, which may be identified with
the \emph{curvature operator} $\check{R}(m,k):{\mathcal N}\to{\mathcal N}$, satisfies 
$$\check{R}(m,k)=(4\Psi_1){\rm Id}_{{\mathcal N}}.$$ 
Recalling 
that $\Psi_1$ is that part of the selfdual 
part of the Weyl tensor, which if vanishes, makes it algebraically special, we have 
the following 
\begin{corollary}
A 4-dimensional manifold $\mathcal M$ with a metric $g$ and an
integrable field of totally null 2-planes $\mathcal N$ is
algebraically special if and only if the characteristic connection
$\check{\nabla}$ of $\mathcal N$ is flat, i.e. iff its curvature
$\check{R}\equiv 0$.  
\end{corollary}

This proves the following Proposition.
\begin{proposition}\label{jhj}
A 4-dimensional manifold $({\mathcal M},g)$ is algebraically special
iff it posesses an integrable field of totally null 2-planes whose
characteristic connection is flat.
\end{proposition}

\subsection{Characteristic connection and the sharpest Goldberg-Sachs theorem}

Given an integrable field of totally null 2-planes $\mathcal N$ we
have the corresponding characteristic connection $\check{\nabla}$. Let
$(f_A)=(f_1,f_2)$ be a frame in $\mathcal N$. In the previous section
we found that the curvarture of $\check{\nabla}$ in the basis
$(f_A)=(m,k)$ is
$$\check{R}^A_{~BCD}=4\Psi_1\delta^A_{~B}\epsilon_{CD},$$
where $A,B,C,D=1,2$, $(\delta^A_{~B})=\bma 1&0\\0&1\ema$, and
$(\epsilon_{CD})=\bma 0&1\\-1&0\ema$. Thus, in particular, the `Ricci
tensor' $\check{R}_{AB}=\check{R}^C_{~ACB}$ of this connection is \emph{antisymmetric} and equal to 
$$\check{R}_{AB}=4\Psi_1\epsilon_{AB}.$$
Since the curvature has only one component, it is obvious that the other possible contraction, namely
$\check{R}^C_{~CAB}$, is proportional to $\check{R}_{AB}$: $\check{R}^C_{~CAB}=2\check{R}_{AB}$. 
Using this Ricci tensor we are able to formulate the following, quite 
elegant, strengthening of the generalization of the Goldberg-Sachs theorem
given in Theorem \ref{gsccc}.
\begin{theorem}\label{ghj}
Let $\tn\subset {\rm T}^\bbC \mathcal M$ be an integrable field of totally null
2-planes on a 4-dimensional manifold
$({\mathcal M},g)$ equipped with metric $g$. 
Assume that the tensor
$\check{\nabla}_{[C}\check{\nabla}_{D]}\check{R}_{AB}$ vanishes
everywhere on $\mathcal M$,
\be
\check{\nabla}_{[C}\check{\nabla}_{D]}\check{R}_{AB}\equiv 0.\label{ghjj}\ee
Then $({\mathcal M},g)$ is algebraically special at every point of
  $\mathcal M$, with a
  multiple 
  field of principal totally null 2-planes tangent to $\tn$.
\end{theorem} 
\begin{proof}
For every connection $\nabla_A$, the action of the 
operator $\nabla_{[C}\nabla_{D]}$ on any tensor is a suitable 
linear action of the curvature of $\nabla_A$ on this tensor. 
Since for $\check{\nabla}_A$ the curvature has only one component $\Psi_1$,
the quantity $\check{\nabla}_{[C}\check{\nabla}_{D]}\check{R}_{AB}$
only involves a constant coefficient sum of terms of the form $\Psi_1\check{R}_{AB}$. Since $\check{R}_{AB}$ itself is proportional to
$\Psi_1$, because of the symmetry, we conclude that     
$$\check{\nabla}_{[C}\check{\nabla}_{D]}\check{R}_{AB}=c~\Psi_1^2\epsilon_{AB}
\epsilon_{CD},\quad\quad\quad c={\rm const}.$$

The constant $c$ may be calculated in a particular basis, e.g. in the
basis $(f_A)=(m,k)$. Using this basis, the definitions (\ref{chc}) and the
Newman-Penrose equations from the Appendix, it is a matter of algebra to
check that $c=-16$.

Now, if $\check{\nabla}_{[C}\check{\nabla}_{D]}\check{R}_{AB}\equiv
0$, then also $\Psi_1^2\equiv 0$, and hence
$\Psi_1\equiv 0$. Since $\tn$ is integrable, then we also have
$\Psi_0\equiv 0$, which means that $\tn$ is a multiple totally null
2-plane. This finishes the proof.
\end{proof}

%In the intrinsinc language of the characteristic connection 
%this theorem can be restated as:
%\begin{theorem}\label{hgn}
%Suppose that a characteristic connection $\check{\nabla}$ of an
%integrable totally null 2-dimensional distribution satisfies 
%$$\check{\nabla}_{[C}\check{\nabla}_{D]}\check{R}_{AB}= 0$$
%everywhere on a 4-dimensional manifold $({\mathcal M},g)$. Then this
%connection is flat.
%\end{theorem}
%
%
%\begin{proof}
%The proof of this statement is an immediate consequence of Theorem
%\ref{ghj} and Proposition \ref{jhj}.
%\end{proof}
\begin{remark}
Since $S$ as in (\ref{s'}) is equal to $-10\Psi_1^2$, and this is turn
is 8/5 of the only component of the conformal tensor
$\check{\nabla}_{[C}\check{\nabla}_{D]}\check{R}_{AB}$, it is now
clear why an `adhoc' defined object $S$ in (\ref{s}) is a \emph{weighted scalar}.
\end{remark}

\begin{remark}
According to the discussion in Example \ref{pk}, the assumption about
the conformal tensor 
$\check{\nabla}_{[C}\check{\nabla}_{D]}\check{R}_{AB}\equiv 0$,
replacing the Ricci flatness condition from the original
Goldberg-Sachs theorem, can not
be weakened if one wants to get the implication
$(\kappa\equiv\sigma\equiv 0)\Rightarrow (\Psi_0\equiv\Psi_1\equiv
0)$. Thus, although the connection $\check{\nabla}$ provides plenty of  
a priori ``weaker'' conditions, such as for example
$\check{\nabla}_{[E}\check{\nabla}_{F]}\check{\nabla}_{[C}\check{\nabla}_{D]
}\check{R}_{AB}\equiv 0$, 
or conditions with more iterations of the curvature operator 
$\check{\nabla}_{[C}\check{\nabla}_{D]}$,  
they all are equivalent to the simplest condition 
$\check{\nabla}_{[C}\check{\nabla}_{D]}\check{R}_{AB}\equiv 0$.  
\end{remark}
\section{Generalizations of the 
Goldberg-Sachs theorem for real metrics}\label{rc}
Theorems \ref{gst}, \ref{gsccc}, \ref{west}, \ref{pol}, \ref{pol1} and
\ref{pol2} were proved 
assuming that the metric $g$ is \emph{complex}. The proofs \emph{also work} 
when $g$ is \emph{real}. To see this it is enough to look at the proofs 
assuming one of the reality conditions ($L$), ($E$), ($S_c$) or
($S_r$) of Remarks \ref{r1}, \ref{r2},
\ref{r3} and \ref{r4}. They impose relations between the components of the
Weyl tensor $\Psi_\mu$ and $\Psi'_\nu$, between the Schouten tensor
components $\Rho_{ab}$ and between the Cotton tensor components
$A_{abc}$. These relations are harmless for the arguments in the
proofs. They however 
may be used to \emph{shorten} the proofs and may cause that some
assumptions appearing in the complex versions can be dropped off. 

We first discuss the Euclidean case. 
\subsection{Euclidean case} In this case, in every null coframe
$(M,P,N,K)$, as in (\ref{npg}), the
reality conditions ($E$) imply that in particular: 
\be
\Psi_4=\bar{\Psi}_0,\quad  \Psi_3=\bar{\Psi}_1,\quad \Psi_2=\bar{\Psi}_2,\quad\Psi_4'=\bar{\Psi}_0',\quad  \Psi_3'=\bar{\Psi}_1',\quad
    \Psi_2'=\bar{\Psi}_2'.\label{ssw}\ee
In the rest of this section we consider the selfdual 
part of the Weyl tensor and principal null 2-planes associated with it. The
analysis of the antiselfdual case is analogous.

Relations (\ref{ssw}), when compared with the equations (\ref{4}) 
defining the principal 2-planes, imply the following:
\begin{proposition}
If $z=z_1$ is a solution of 
$\Psi_4 z^4-4\Psi_3 z^3+6\Psi_2 z^2+4\Psi_1 z+\Psi_0=0$
then is so $z_2=-\frac{1}{\bar{z}_1}$.
\end{proposition}
\begin{proof}
Inserting (\ref{ssw}) and $z=z_1$ 
in the equation defining the principal null 2-planes (\ref{4}) we get
$$\bar{\Psi}_0 z_1^4-4\bar{\Psi}_1 z_1^3+6\Psi_2 z_1^2+4\Psi_1
z_1+\Psi_0=0.$$
Now dividing this by $z_1^{-4}$ and taking the complex conjugation of
the result, we get 
$$\bar{\Psi}_0 z_2^4-4\bar{\Psi}_1 z_2^3+6\Psi_2 z_2^2+4\Psi_1
z_2+\Psi_0=0,$$ 
which finishes the proof.
\end{proof}
Comparing this with Proposition \ref{eus} we have 
\begin{corollary}
Principal null 2-planes always appear in pairs corresponding to 
pairs of solutions $(z_1,z_2)=(z_1,\frac{-1}{\bar{z}_1})$ of equation
(\ref{4}).\\
A pair of solutions $(z_1,z_2)=(z_1,\frac{-1}{\bar{z}_1})$ of equation
(\ref{4}) at a point $x$ distinguishes a pair
$(J(z_1),J(\frac{-1}{\bar{z}_1}))$ of principal 
hermitian structures $J(z_1)$ and $J(\frac{-1}{\bar{z}_1})$ at $x$,
which are conjugate to each other, 
$\overline{J(\frac{-1}{\bar{z}_1})}=-J(z_1)$.\end{corollary}

\begin{proof}
The only thing to be proven is
$\overline{J(\frac{-1}{\bar{z}_1})}=-J(z_1)$. By definition of these
two structures we have
$J(z)(m+zn)=i(m+zn)$, $J(z)(k-zp)=i(k-zp)$
and
$J(\frac{-1}{\bar{z}})(m-\frac{1}{\bar{z}}n)=i(m-\frac{1}{\bar{z}}n)$,
$J(\frac{-1}{\bar{z}})(k+\frac{1}{\bar{z}}p)=i(k+\frac{1}{\bar{z}}p).$
The second set of equations is equivalent to 
$J(\frac{-1}{\bar{z}})(\bar{z}m-n)=i(\bar{z}m-n)$ and
$J(\frac{-1}{\bar{z}})(\bar{z}k+p)=i(\bar{z}k+p)$, which after
conjugation and the use of the reality conditions ($E$) gives:
\begin{eqnarray*}
&&\overline{J(\tfrac{-1}{\bar{z}})}(k-zp)=-i(k-zp)=-J(z)(k-zp),\\
&&\overline{J(\tfrac{-1}{\bar{z}})}(m+zn)=-i(m+zn)=-J(z)(m+zn).
\end{eqnarray*}
\end{proof}
This corollary implies that at each point $x$ of $\mathcal M$ the
selfdual part of the Weyl tensor may be in one of the following Petrov types:
\begin{itemize}
\item[type G:] the generic type, in which the selfdual part of the Weyl tensor
  does not vanish at $x$, and in which we have two distinct pairs 
$(z_1,z_2)=(z_1,\frac{-1}{\bar{z}_1})$ and
  $(z_3,z_4)=(z_3,\frac{-1}{\bar{z}_3})$, $z_1\neq z_3$, of solutions
  of equation (\ref{4}). In such case the pairs $(z_1,z_2)$ and
  $(z_3,z_4)$ correspond to two pairs
  of \emph{different} mutually conjugate principal hermitian
  structures $(J(z_1),J(z_2))$ and
  $(J(z_3),J(z_4))$ at $x$.
\item[type D:] this is the degeneracy of type G. It occurs when $z_1$ is a
  double root of (\ref{4}), i.e. when $z_3=z_1$. In such case we have 
only \emph{one} pair of double
  principal hermitian structures $(J(z_1),J(z_2))$ at $x$. 
\item[type 0:] this is the antiselfdual type in which the selfdual
  part of the Weyl tensor \emph{vanishes} at $x$. In this case the sphere of
  selfdual 2-planes has no distinguished points. 
\end{itemize}
Note that always we may choose a Newman-Penrose frame in which
$\Psi_0=0$ at $x$. In types G or D it is achieved by choosing
the Newman-Penrose vectors $m$ and $k$ such that they span the
principal null 2-plane corresponding to $z_1$. Then, in such a frame, 
the algebraically special type D is characterized by $\Psi_1=0$ and
$\Psi_2\neq 0$ at
$x$. If in such a frame $\Psi_1\neq 0$, then the selfdual part of the 
Weyl tensor
is algebraically general (of type G) at $x$.
    
This proves the following 
\begin{theorem}
At every point of a 4-dimensional manifold $\mathcal M$ equipped with 
a real euclidean-signature metric $g$ the selfdual part of the 
Weyl tensor may be of
one of the types G, D, and
0, with the analogous types for the antiselfdual part of the Weyl tensor. Thus, at
eavery point of a 4-manifold $\mathcal M$ equipped with a \emph{euclidean signature}
metric $g$ we have $3\times 3=9$ `Petrov' types. 
\end{theorem}

Thus the Euclidean reality conditions ($E$) imply that the number of
possible 
Petrov types in the Euclidean case is much smaller than in the complex
case. This implies that the complex theorems of the previous Section
have much stronger Euclidean versions. In particular, the proof of 
Theorem \ref{pol}, when the reality conditions ($E$) are assumed, 
goes through as in the complex version, with the only exception, that 
the II-generiticity property of $g$ may now be weakend to the
assumption that the selfdual part of the Weyl tensor is nowhere vanishing (or even
to a still weaker assumption that the points at which the selfdual
part of the Weyl tensor
vanishes form closed sets without interior). Indeed, in the Euclidean
case, the assumption $\Psi_0\equiv\Psi_1\equiv 0$ and $\Psi_2\neq 0$, which is 
needed for the conclusion that $\kappa\equiv\sigma\equiv 0$, means only that
the selfdual part of the Weyl tensor is nonvanishing, since now
$\Psi_0\equiv\Psi_1\equiv 0$ implies that $\Psi_4\equiv\Psi_3\equiv
0$. This proves the Riemannian version of the Goldberg-Sachs 
Theorem \ref{polo}.

One of the corollaries from the complex Theorem
\ref{pol} is also the following 
\begin{corollary}
If the selfdual part of the Weyl tensor of a real metric $g$ of Riemannian
signature does 
not vanish on a 4-dimensional manifold $\mathcal M$, then modulo
complex conjugation, such a metric admits at most two
hermitean structures that agree with the orientation. If such
hermitean structures exist
their spaces of (1,0) vectors coincide with the selfdual 
principal totally null 2-planes. In particular, in type D we may have
only one hermitean structure, which exists if and only if the
Cotton tensor for $g$ vanishes on its space of (1,0) vectors.  
\end{corollary}

The Euclidean version of Theorem \ref{gst} is also worth quoting. We
have
\begin{corollary}\label{gstteh}
Assume that a 4-dimensional manifold $\mathcal M$ equipped with a real
metric of Riemannian signature $g$ has a nonvanishing selfdual part of
the Weyl
tensor $C^+$. Suppose that it admits a hermitean structure $J$ which agrees with
the orientation, and that its Ricci tensor vanishes on the space
$\mathcal N$ of (1,0) vectors of $J$. Then $C^+$ is of type D, 
with $\mathcal N$ being the only principal selfdual null 2-plane.  
\end{corollary}  
\subsection{Split signature case} To spell out all the possible
Petrov types and their interpretations in this case we first consider 
the Newman-Penrose coframe $(M,P,N,K)$ with the reality conditions
($S_c$) from Remarks \ref{r1} and \ref{r2}. In this coframe the sphere
of selfdual totally null 
2-planes ${\mathcal N}_z$ is spanned by $m+zn$ and $k-zp$
as in (\ref{plan}). Now, having the reality conditions $S_c$, we ask
which values of $z\in\bbC$ correspond to the \emph{nongeneric} 
selfdual totally null 
2-planes which have real index equal to \emph{two}. We have
the following 
\begin{proposition}\label{lkp}
A selfdual 2-plane ${\mathcal N}_z$ has real index equal to two if and
only if the complex parameter $z\in\bbC$ lies on the unit circle $z\bar{z}=1$. 
\end{proposition}
\begin{proof}
Due to the reality conditions ($S_c$) 
a \emph{real} nonvanishing 
vector $v=a(m+zn)+b(k-zp)$ from ${\mathcal N}_z$ must satisfy 
$$a(m+zn)+b(k-zp)=\bar{a}(p-\bar{z}k)+\bar{b}(-n-\bar{z}m).$$ 
Equating to zero the respective coefficients at $m,p,n,k$ we easily
get that this is possible if and only if $z\bar{z}=1$. 
Thus ${\mathcal N}_z$ includes real nozero vectors if and only if
$z\bar{z}=1$. We further observe that if $z\bar{z}=1$ then $v$ is
real if and only if $b=-\bar{a}\bar{z}$. Thus, when $z$ is fixed, we
have a 1-complex-parameter-family $v=v(a)$ of real vectors in ${\mathcal
  N}_z$. Choosing two diffrent values of $a$ we get
$$v(a)\dz v(a')=(a\bar{a}'-a'\bar{a})(m\dz p-\bar{z}m\dz k-zp\dz n-n\dz k).$$
This shows that ${\mathcal N}_z$ with $z\bar{z}=1$ includes
\emph{independent} real vectors (take e.g. $a=1$ and $a'=i$), thus it has real inedex two. This
finishes the proof.  
\end{proof}
Let us now choose a Newman-Penrose coframe as in (\ref{npg}). Then the
reality conditions ($S_c$) imply that we have:
\be
\Psi_4=\bar{\Psi}_0,\quad  \Psi_3=-\bar{\Psi}_1,\quad \Psi_2=\bar{\Psi}_2,\quad\Psi_4'=\bar{\Psi}_0',\quad  \Psi_3'=-\bar{\Psi}_1',\quad
    \Psi_2'=\bar{\Psi}_2',\label{swc}\ee
and the reality conditions ($S_r$) mean that all Weyl tensor
coeffcients $\Psi$ and $\Psi'$ are real:
\be\Psi_0=\bar{\Psi}_0,\quad  \Psi_1=\bar{\Psi}_1,\quad
\Psi_2=\bar{\Psi}_2,\quad\Psi_3=\bar{\Psi}_3,\quad
\Psi_4=\bar{\Psi}_4,\label{swr}\ee
(we also have anologous relations for $\Psi'$). 
 
We pass to the split signature version of the Petrov
classification. We perform the analysis for the selfdual part of the Weyl tensor;
the classification for the antiselfdual case is anlogous. 

Let us fix a point $x\in\mathcal M$. Let $(M,P,N,K)$ be a 
Newman-Penrose coframe around $x$ satsifying the reality conditions
($S_c$), and as a consequence (\ref{swc}). We have the following 
\begin{proposition}
If $z=z_1$ is a solution of 
$\Psi_4 z^4-4\Psi_3 z^3+6\Psi_2 z^2+4\Psi_1 z+\Psi_0=0$
then is so $z_2=\frac{1}{\bar{z}_1}$.
\end{proposition}
\begin{proof}
Inserting (\ref{swc}) and $z=z_1$ 
in the equation defining the principal null 2-planes (\ref{4}) we get
$$\bar{\Psi}_0 z_1^4+4\bar{\Psi}_1 z_1^3+6\Psi_2 z_1^2+4\Psi_1
z_1+\Psi_0=0.$$
Now dividing this by $z_1^{-4}$ and taking the complex conjugation of
the result, we get 
$$\bar{\Psi}_0 z_2^4+4\bar{\Psi}_1 z_2^3+6\Psi_2 z_2^2+4\Psi_1
z_2+\Psi_0=0,$$ 
which finishes the proof.
\end{proof}
Comparing this Proposition with Proposition \ref{lkp} we get
\begin{corollary}
Selfdual principal null 2-planes always appear in pairs corresponding to pairs 
of solutions $(z_1,z_2)=(z_1,\frac{1}{\bar{z}_1})$ of equation
(\ref{4}). The situation in which $z_1=z_2$ happens only if the
principal selfdual null 2-plane has real index two.
\end{corollary}
Using Proposition \ref{eus} we may also reinterpret this corollary as follows
\begin{corollary}
If equation (\ref{4}) at a point $x$ admits a principal selfdual null 2-plane 
of real index zero, then at this point we have two distinguished
hermitian structures $J(z_1)$ and $J(\frac{1}{\bar{z}_1})$ associated
with the solution $z_1$ of (\ref{4}). Moreover these two structures
are conjugate to each other.
\end{corollary} 
\begin{proof}
The only thing to be proven is
$\overline{J(\frac{1}{\bar{z}_1})}=-J(z_1)$. By definition of these
two structures we have
$J(z)(m+zn)=i(m+zn)$, $J(z)(k-zp)=i(k-zp)$
and
$J(\frac{1}{\bar{z}})(m+\frac{1}{\bar{z}}n)=i(m+\frac{1}{\bar{z}}n)$,
$J(\frac{1}{\bar{z}})(k-\frac{1}{\bar{z}}p)=i(k-\frac{1}{\bar{z}}p).$
The second set of equations is equivalent to 
$J(\frac{1}{\bar{z}})(\bar{z}m+n)=i(\bar{z}m+n)$ and
$J(\frac{1}{\bar{z}})(\bar{z}k-p)=i(\bar{z}k-p)$, which after
conjugation and the use of the reality conditions ($S_c$) gives:
\begin{eqnarray*}
&&\overline{J(\tfrac{1}{\bar{z}})}(k-zp)=-i(k-zp)=-J(z)(k-zp),\\
&&\overline{J(\tfrac{1}{\bar{z}})}(m+zn)=-i(m+zn)=-J(z)(m+zn).
\end{eqnarray*}
\end{proof}

Because of quite different
reality conditions (\ref{swc}) and
(\ref{swr}) at each point $x\in\mathcal M$ 
we need to consider separately two different cases: the
generic one a) in which the selfdual part of the Weyl tensor admits \emph{at least one} 
principal totally null 2-plane of real index \emph{zero} at $x$, and the less
generic one b) in which \emph{all} principal null planes have real
index \emph{two} at $x$. 

In the case a) we chose a Newman-Penrose coframe $(M,P,N,K)$ around
$x$ such that it satisfies the reality conditions ($S_c$) and that the
principal totally null 2-plane of real index zero corresponds to the
solution $z=0$ of (\ref{4}). Then in such a coframe $\Psi_0= 0$,
and the equation defining the principal null 2-planes becomes 
$4\bar{\Psi}_1 z^3+6\Psi_2 z^2+4\Psi_1z=0$, or 
\be2\bar{\Psi}_1 z^2+3\Psi_2 z+2\Psi_1=0.\label{qua}\ee
Thus in this coframe we have two solutions $(z_1,z_2)=(0,\infty)$
corresponding to the mutually conjugate principal (almost) hermitian
structures associated with two fields of principal 2-planes of index
zero, and the rest of the principal 2-planes has to be determined as 
solutions to the quadratic equation (\ref{qua}). The roots of this
equations are obviously 
$$z_{3,4}=\frac{-3\Psi_2\pm\sqrt{9\Psi_2^2-16\Psi_1\bar{\Psi}_1}}{4\bar{\Psi}_1}.$$
The interpretation depends on the sign of
$9\Psi_2^2-16\Psi_1\bar{\Psi}_1$ and on whether $\Psi_1$ vanishes or
not.
It follows that at each point $x\in\mathcal M$ we have now \emph{four}
cases:
\begin{itemize}
\item[type G:] the generic case in which $z_3\neq z_4=\frac{1}{\bar{z}_3}$, $z_3\bar{z}_3\neq
  1$, $z_3\neq 0$ and $z_3\neq\infty$. In such case we have two pairs
  of \emph{different} mutually conjugate principal hermitian
  structures at $x$ corresponding to $(J(0),J(\infty))$ and
  $(J(z_3),J(\frac{1}{\bar{z}_3}))$. This case happens when 
  $9\Psi_2^2>16\Psi_1\bar{\Psi}_1$ and $\Psi_1\neq 0$ at $x$.
\item[type SG:] in this case $z_3\neq z_4=\frac{1}{\bar{z}_3}$, $z_3\bar{z}_3=
  1$. Here, in addition to the pair
  of mutually conjugate principal hermitian
  structures $(J(0),J(\infty))$ at $x$, we have \emph{two different}
  principal totally null 2-planes of real index \emph{two} at
  $x$. These real 2-planes are associated with the solultions $z_3$ and
  $z_4$, which lie on the circle $z\bar{z}=1$. This case happens when 
  $9\Psi_2^2<16\Psi_1\bar{\Psi}_1$ at $x$.
\item[type II:] this is the degenerate case of the type SG. It happens when 
$9\Psi_2^2=16\Psi_1\bar{\Psi}_1$ and $\Psi_1\neq 0$ at $x$, and the equation (\ref{qua})
  has double root $z_3=z_4$ at $x$. We necessarily have
  $z_3\bar{z}_3=1$ in this case, and thus, in addition to the pair of
  mutually conjugate principal hermitian structures
  $(J(0),(J(\infty))$ we have also \emph{one} double principal null 2-plane of
  real index \emph{two} at $x$. 
\item[type D:] this is another degeneration of the type G. Now $\Psi_1=0$ at
  $x$ and we have $z_3=0$ and $z_4=\infty$ as solutions of
  (\ref{qua}). Thus in this case the points $z=0$ and $z=\infty$ have 
  multiplicity \emph{two}, and we have only \emph{one} pair of double
  principal hermitian structures $(J(0),J(\infty))$ at $x$. 
\end{itemize}
We now pass to the cases in which we do \emph{not} have a single
principal null 2-plane which has a real index zero at $x$. The
analysis here could still be performed in the Newman-Penrose coframe
satisfying the reality conditions $S_c$, but since now all the solution of
equation (\ref{4}) would have to satsify $z\bar{z}=1$, we would not be able to
choose the frame in such a way that $\Psi_0$ would be zero at $x$. This would
lead to the analysis of the roots of the \emph{quartic} equation
(\ref{4}), and it is why it is
now much easier to reason in the coframe that satsifies the reality
conditions ($S_r$). So now, we choose a
Newman-Penrose coframe $(M,P,N,K)$ around $x$, which satisfies the 
reality conditions ($S_r$) and, since now we have at least one
principal null 2-plane of real index two at $x$, we may assume that we have $\Psi_0=0$ at $x$. In
this coframe our principal totally null 2-plane of real index two
corresponds to $z_1=0$ and the other principal 2-planes are determined
by 
$$\Psi_4 z^3-4\Psi_3 z^2+6\Psi_2 z^2+4\Psi_1=0.$$
Here all the $\Psi_1$, $\Psi_2$, $\Psi_3$ and $\Psi_4$ are \emph{real}
and we admit only \emph{real solutions} for $z$. (If the solution is
complex, it corresponds to a 2-plane with real index zero, and
corresponds to one of the cases G, SG, II, or D, considered earlier.)

Now, a XVth century substitution $z\to z-\frac{4\Psi_3}{3\Psi_4}$, brings this
equation into the form $z^3+pz+q=0$, which has
\emph{three real} roots for $z$ iff $27p^4+4q^3\geq 0$. This inequality gives
the restriction on the Weyl tensor, which determines the situation we
are talking about here. If the selfdual part of the Weyl tensor satisfies this
restriction, the equation (\ref{4})
has \emph{four real roots}. This, in addition to G, SG, II and D, 
defines the five new Petrov types:
\begin{itemize}
\item[type G$_r$:] equation (\ref{4}), written in the coframe with reality
conditions ($S_r$), has \emph{four different real} roots, meaning that
we have \emph{four} different principal null 2-planes of real index
two at x,
\item[type II$_r$:] equation (\ref{4}), written in the coframe with reality
conditions ($S_r$), has \emph{one double and two different real} roots, meaning that
we have \emph{three} different principal null 2-planes of real index
two at x, one of them with multiplicity two,
\item[type III$_r$:] equation (\ref{4}), written in the coframe with reality
conditions ($S_r$), has \emph{one triple and one distinct real} roots,
meaning that 
we have \emph{two} different principal null 2-planes of real index
two at x, one of them with multiplicity three, 
\item[type N$_r$:] equation (\ref{4}), written in the coframe with reality
conditions ($S_r$), has \emph{one quadruple} root, meaning that
we have a \emph{single quadruple} principal null 2-plane of real index
two at x,
\item[type D$_r$:] equation (\ref{4}), written in the coframe with reality
conditions ($S_r$), has \emph{two distinct double real} roots, meaning that
we have \emph{two} different principal null 2-planes of real index
two at x, each of them having multiplicity two. 
\end{itemize} 
Finally we have the Petrov type corresponding to the situation when
the \emph{selfdual part of the Weyl tensor vanishes} at $x$ (the metric is
\emph{antiselfdual at $x$}).

This proves the following 
\begin{theorem}
At every point of a 4-dimensional manifold $\mathcal M$ equipped withh
a real split-signature metric $g$ the selfdual part of the Weyl tensor may be of
one of the types G, SG, II, D, G$_r$, II$_r$, III$_r$, N$_r$, D$_r$,
0, with the analogous types for the antiselfdual part of the 
Weyl tensor. Thus, at
every point of a 4-manifold $\mathcal M$ equipped with a \emph{split signature}
metric $g$ we have $10\times 10=100$ `Petrov' types. 
\end{theorem}
The above analysis also suggest the following terminology: 
the name \emph{algebraically
  special} for the selfdual part of the Weyl tensor in the split signature case is 
reserved to the types II, D,  II$_r$, III$_r$, N$_r$, D$_r$ and
0 only. Although the types SG
and G$_r$ are algebraically (and geometrically!) distinguished from
the most general case $G$, we also call them \emph{algebraically
  general}. With
this terminology, Theorems \ref{gstts} and \ref{gstts1} 
follow from our Theorem \ref{gst}.

Because of the huge number of the algebraically special cases to be 
considered, we skip the discussion of the split signature versions 
of further theorems from Section \ref{prof} here. Such a discussion 
deserves a separate paper. This should also answer several interesting
questions, such as for example, the following: `are there split-signature 
Einstein metrics of type II?', `is it possible to
have a split signature Einstein 4-manifold on which an integrable totally null
2-plane can change its real index from 0 to 2?', etc. 

We close this section by mentioning the recent paper \cite{law}. It 
is entirely devoted to the Newman-Penrose formalism
adapted to the split signature situation, and it provides a 
version of the split-signature Goldberg-Sachs theorem.
\subsection{Lorentzian case} Here the Petrov types are precisely the
same as in the complex case described by the Definition
\ref{pett}, i.e. we have types $G$, $II$, $D$, $III$, $N$ and $0$
here. The Lorentzian 
reality conditions ($L$) \emph{do not make any
restriction} on the Weyl tensor coefficients $\Psi_\mu$. What they do is,
they give a simple ralation between the self-dual part of the Weyl tensor and 
the antiselfdual one. We have 
$\Psi_\mu'=\bar{\Psi}_\mu$, so here the antiselfdual part of the Weyl tensor is
totally determined by the selfdual one. Since in the proofs in Section
\ref{prof} the coefficients $\Psi_\mu'$ never appear, and only $\Psi_\mu$s
matter, all the proofs, and the theorems presented in
Section \ref{prof} restrict naturally to the Lorentzian case without
any alteration. 

However, since in the Lorentzian signature the fields of totally null
2-planes have always real index one, it is customary to formulate
the Lorentzian theorems in terms of the real vector field $k$ such that
$\Span_{\bbC}(k)=\nen\cap\bar{\nen}$. In particular, such a null real
vector field is said to be \emph{geodesic and shear-free}
\cite{rt} 
if it satsifies 
\be
{\mathcal L}_kg=ag+g(k)\om,
\label{rt}
\ee
with a function $a$ and a 1-form $\om$ on $\mathcal M$. Here $g(k)$ is a
1-form on $\mathcal M$ such that $X\hook g(k)=g(k,X)$ for any vector field
  $X\in{\rm T}\mathcal M$. When written in terms of the field $\tn$ of the
  associated totally null 2-planes, condition (\ref{rt}) is equivalent
  to  
$$[\tn,\tn]\subset\tn,$$
i.e. to the formal integrability condition for $\tn$.

Suppose now the Weyl tensor $C_{abcd}$ of $({\mathcal M},g)$ is 
\emph{nonvanishing}. It is well known \cite{car} that the algebraic equation   
\be
k_{[e}C_{a]bc[d}k_{f]}k^b k^c=0,\label{pnd}
\ee
for a null vector $k$ has at most \emph{four} solutions at every point
$x\in \mathcal M$. The solutions $k$ of equation (\ref{pnd}) at $x\in
\mathcal M$ are
called the \emph{principal null directions} (PNDs) at $x$. If equation
(\ref{pnd}) admits \emph{exactly} four PNDs at $x\in \mathcal M$ then
$({\mathcal M},g)$
is said to be \emph{algebraically general} at $x$. If the number $q$ of
solutions to (\ref{pnd}) at $x\in \mathcal M$ is $1\leq q\leq 3$ then
$({\mathcal M},g)$
is called \emph{algebraically special} at $x$. In such case the
\emph{quartic} equation (\ref{pnd}) has at least one 
\emph{multiple root}, and
the solution $k$ corresponding to it is called a \emph{multiple}
PND. This notion of the algebraical speciality coincides with the one
in terms of the principal null 2-planes, since on a Lorentzian
oriented and time oriented 4-manifold $\mathcal M$, there is one to
one correspondence between fields of totally null 2-planes in the
complexification and real null vector fields, defined by the
intersection of the 2-planes with their complex conjugations.

Having said this, we present the Lorentzian version of our complex
Theorem \ref{gscc}.

\begin{theorem}\label{gsc}
Let $\tn\subset {\rm T}^\bbC \mathcal M$ be a field of totally null 2-planes on a Lorentzian 4-dimensional manifold
$({\mathcal M},g)$. Assume that the Ricci tensor $Ric$ of 
$({\mathcal M},g)$, considered as a symmetric bilinear form
on ${\rm T}^\bbC \mathcal M$, is degenerate on $\tn$,
$$Ric_{|\tn}=0.$$
If in addition the field $\tn$ is integrable, $[\tn,\tn]\subset\tn$, everywhere
  on $\mathcal M$, or what is the same, if $k$ such that ${\rm Span}_\bbC(k)=\tn\cap\bar{\tn}$, is
  geodesic and shear-free, then 
$({\mathcal M},g)$ is algebraically special at every point, with a
  multiple 
  PND tangent to $k$.
\end{theorem}

\begin{remark}
In \cite{lhn} we used Theorem \ref{gsc} without proof, since it would have
made an already long paper even longer. Actually some statements 
equivalent to Theorem \ref{gsc} are known to a
few general relativists, see e.g. 
Lemma 2.2 on p. 577 of \cite{plprz2}. Since this equivalence is not
easy to decipher, we decided to present this theorem here, as a
corollary from the complex Theorem \ref{gscc}.
\end{remark}
\subsection{Counterexample to Trautman's conjecture }\label{cout}
Trautman in \cite{at} asked if there exists an example of a
4-dimensional Bach flat metric with nonvanishing selfdual part of the 
Weyl tensor $C^+$, 
for which an integrable field of selfdual totally null 2-planes would not be
principal for $C^+$. He \emph{conjectured} that the answer to this question
is `no'. Although the question was formulated in
the Lorentzian setting, it makes sense in any signature. It is also very
closely related to the Goldberg-Sachs theorem. 

Our analysis of this
theorem from Section \ref{lpo}, especially the discussion in Example
\ref{pk}, suggests that the examples Trautman asks about, should be
possible. This is because, the conditions needed for `if and only if'
between conditions (i) and (ii) in Theorem \ref{pol} are related to
those derivatives of the Cotton tensor that are not present in the
Bach tensor. This is clear from Example \ref{pk}: the 
integrability condtions for $A_{|\mathcal N}\equiv 0$, give $S\equiv
0$, where $S$ is given by (\ref{s'}). And although the Bach tensor
components may be obtained by differentiating some components of the
Cotton tensor, the derivatives of the Cotton tensor appearing in $S$ are
not (at least algebraically) 
expressible in terms of the components of the Bach tensor.     

Below in this Section we present a simple example of a metric with
\emph{Euclidean} signature which is Bach-flat, 
admits an \emph{integrable hermitean structure} which \emph{agrees with the
orientation}, and whose selfdual part of the Weyl tensor is of general type G.

On $\bbR^4$, with local coordinates $(x^1,x^2,x^3,x^4)$, consider
$z=x^1+i x^2$ and $w=x^3+i x^4$, and a complex-valued function
$f=f(w,z)$ \emph{holomorphic} in \emph{both} arguments $w$ and
$z$. Given $f$ define a Riemannian metric 
%%/home/pawel/notebooks/stonybrook/goldbergsachs/countertrautman.nb
$$g=2\Big(\der w\der\bar{w}+{\rm
  exp}\big(f(w,z)+\bar{f}(\bar{w},\bar{z})\big)\der
z\der\bar{z}\Big).$$
Now introduce the Newman-Penrose coframe by setting 
$$M=\der\bar{w},\quad P=\der w,\quad N={\rm e}^f\der z,\quad K={\rm
  e}^{\bar{f}}\der\bar{z}.$$
They obviously satisfy the Euclidean reality conditions ($E$). A short
calculation shows, that modulo the complex conjugation, the only nonvanishing Newman-Penrose
coefficients are:
$$\al=-\tfrac12\pi=\beta'=-\tfrac12\tau'=\tfrac14f_w.$$
In particular $\kappa=\sigma=0$, which is obvious since the field of
selfdual totally null 2-planes
$\mathcal N$ spanned by
$m=\partial_{\bar{w}}$ and $k={\rm e}^{-\bar{f}}\partial_{\bar{z}}$ is integrable.
Now our main point is that the only nonvanishing components of the Weyl
tensor are:
$$\Psi_3=\bar{\Psi}_1=\tfrac14{\rm e}^{-f}f_{wz}.$$
This in particular means that the field $\mathcal N$ is principal
(since $\Psi_0\equiv 0$), but when $f_{wz}\neq 0$ it is \emph{not}
multiple ($\Psi_3\neq 0\neq\Psi_1$). Moreover, since
$\Psi_0'\equiv\Psi_1'\equiv\Psi_2'\equiv\Psi_3'\equiv\Psi_4'\equiv 0$,
i.e. the full antiselfdual part of the Weyl tensor identically vanishes, the
metric is \emph{Bach flat}. This answers in positive the question of Trautman
we mentioned at the begining of this Section. Moreover, if $f_{wz}\neq
0$, due to the
Corollary \ref{gstteh}, this selfdual metric can not have Ricci tensor
vanishing on $\mathcal N$, and as such is never
conformal to an Einstein metric.   

\subsection{Characteristic connection in real signatures}\label{chch}
We now reexamine the arguments from Section \ref{cchh} from the point
of view of the reality conditions. 

From Step one of the proof of Theorem \ref{cch1} we know that the Weyl
form $B$ of the Weyl connection which preserves an integrable $\tn$,
in an adapted to $\mathcal N$ coframe is given by 
$B=2\tau M+ B_2 P+ B_3 N-2\rho K$.
Thus in the complex case (or in the real cases in which we do not
insist on $B$ to be real) the Weyl 1-form is not totally determined by
$\tn$. 

The situation is quite different in the Riemannian ($E$) and the split
signature ($S_c$). In these two cases, the requirements that $B$
is real determines it completely! Indeed, it is easy to see that 
the reality conditions ($E$) or ($S_c$) together with the requirement
that $B$ be real implies that $B$ is equal to
\be
B=2\tau M+ 2\pi P-2\mu N-2\rho K\label{ebe}\ee
or, what is the same,
$$\tfrac12 B=\Gamma_{143}\theta^1+\Gamma_{234}\theta^2+\Gamma_{321}\theta^3+\Gamma_{412}\theta^4.$$

This proves the following theorem
\begin{theorem}\label{pshh}
Let $\mathcal N$ be a field of totally null 2-planes on $({\mathcal
  M},g)$, where $g$ is a 4-dimensional metric of Riemannian or split 
signature. Let us assume that $\mathcal N$ is integrable
$[{\mathcal N},{\mathcal N}]\subset\mathcal N$ and that it has a real
index 0 everywhere on $\mathcal M$. Then there exists a \emph{canonical}
Weyl connection $\nw$ 
on $\mathcal M$,
which encodes the conformal properties of the structure $({\mathcal
  M},g,\tn)$.

The connection $\nw$ is uniquely determined by the requirements that
\begin{itemize}
\item it is real,
\item it is torsionless,
\item it satisfies: $\nw g=-B g$,
\item it satisfies: $\nw_X\tn\subset\tn$ for all $X\in T\mathcal M$.
\end{itemize} 
In terms of a coframe $(\theta^a)$ adapted to $\tn$ and the 
connection 1-forms $\gw^a_{~b}=g^{ad}\gw_{dbc}\theta^c$ the connection
$\nw$ is given by 
$$\gw_{abc}=\Gamma_{abc}+\tfrac12
(g_{ca}B_b-g_{cb}B_a+g_{ab}B_c)$$ with 
$$\tfrac12
B=\Gamma_{143}\theta^1+\Gamma_{234}\theta^2+\Gamma_{321}\theta^3+\Gamma_{412}\theta^4.$$
Here $\Gamma_{abc}$ are the Levi-Civita connection coefficients in the
adapted coframe.
\end{theorem}

\begin{definition}
Let $J$ be a hermitean (or pseudohermitean) structure on an $2n$-dimensional manifold
$({\mathcal M},g)$ with a metric of Riemannian (or
split) signature. A torsionless conection $\nwh$ on
$({\mathcal M},g,J)$ is called (pseudo)hermitean-Weyl iff
\begin{itemize}
\item[-~~] $\nwh J=0$,
\item[-~~] and $\nwh g=-B g$ for some real 1-form $B$ on $\mathcal M$.
\end{itemize}
\end{definition}

According to our discussion in Section \ref{se2}, integrable 
totally null 2-planes
of real index 0 on a 4-dimensional manifold $({\mathcal M},g)$ are in one-to-one
correspondence with (pseudo)hermitean structures $J$ on $({\mathcal
  M},g)$, thus Theorem \ref{pshh} can be reformulated as:
\begin{theorem}\label{npsh}
Every 4-dimensional (pseudo)hermitean manifold $({\mathcal M},g,J)$
defines a canonical (pseudo)hermitean-Weyl connection $\nwh$. This connection encodes
the conformal properties of the structure $({\mathcal M},g,J)$. It is
given by $\nwh=\nw$, where $\nw$ is as in Theorem \ref{pshh}.
\end{theorem} 
Thus in the (pseudo)hermitean case there is a better connection,
namely $\nw$,  
than the characteristic connection $\check{\nabla}$. It is better, since it enables to
differentiate \emph{any} vector from the tangent space of $\mathcal M$ along \emph{any} other vector
from T$\mathcal M$. The connection $\check{\nabla}$ enables for the differentiation along
${\mathcal N}={\rm T}^{(1,0)}{\mathcal M}$ only. And, $\nw$ is better, 
because it contains \emph{much more information} than
$\check{\nabla}$. In particluar,  
$\check{\nabla}$ is simply the restriction of $\nw$ to
$\mathcal N$.

We now pass to the (pseudo)hermitean part of our elegant Goldberg-Sachs Theorem
\ref{ghj}. 
%/home/pawel/notebooks/stonybrook/goldbergsachs/weylkoneksjahermit.nb

We need some preparations:

Given the (pseudo)hermitean-Weyl connection $\nw$, as in Theorem
\ref{npsh}, we use the formula (\ref{we4})
to pass to the connection 1-forms $\gw_{ab}=\gw_{abc}\theta^c$. Here
$(\theta^c)$ is a coframe \emph{adapted} to $J$. The word `adapted' 
(in accordance with the discussion in Section \ref{se2}) means that
the considered coframe is adapted to
$\tn={\rm T}^{(1,0)}\mathcal M$ as in the
definition of this notion at the begining of Section \ref{prof}. Now, there is a sequence of
definitions, which closely mimics the situation in Riemannian
geometry:

Having the connection 1-forms $\gw_{ab}$, the metric $g$ and its inverse, represented by
$g^{ab}$, we also have the 1-forms $\gw^a_{~b}=g^{ac}\gw_{cb}$. Using
them, we define the curvature of the connection $\nw$. We do it, in terms of the curvature
2-forms $\ow^a_{~b}$, analogous to those given in the formula
(\ref{ca2}), by:
$$\tfrac12\rw^a_{~bcd}\theta^c\dz\theta^d=\der\gw^a_{~b}+\gw^a_{~c}\dz\gw^c_{~b}.$$
Here $\rw^a_{~bcd}$ are the curvature coefficents in the coframe
$(\theta^a)$. Then we define the Ricci tensor $$\rw_{ab}=\rw^c_{~acb},$$
and its scalar $$\rw=g^{ab}\rw_{ab}.$$ The next step is to define the
Schouten tensor 
$$\row_{ab}=\tfrac12 \rw_{ab}-\tfrac{1}{12}\rw g_{ab}$$ and the Cotton
tensor 
$$\aw_{abc}=2\nw_{[b}\row_{c]a}.$$ This defines a linear map 
$$\aw:{\rm T}{\mathcal M}\times{\rm T}{\mathcal M}\times{\rm
  T}{\mathcal M}\to\bbR$$ given by $$\aw=\tfrac12 \aw_{abc}\theta^a\otimes(\theta^b\dz\theta^c).$$  Then the (pseudo)hermitean part of Theorem
\ref{ghj} is:
\begin{theorem}\label{scsi}
Let  $({\mathcal M},g,J)$ be a 4-dimensional (pseudo)hermitean manifold
and let $\nw$ be its canonical (pseudo)hermitean-Weyl connection $\nwh$.
Assume that 
\be\nw_{X}\aw(Y,X,Y)\equiv \nw_{Y}\aw(X,X,Y)\quad{\rm for~ all~
  vectors}\quad X,Y\in \tn={\rm
  T}^{(1,0)}{\mathcal M}. \label{scs}\ee
Then the selfdual part of the Weyl tensor for $({\mathcal M},g)$ is algebraically special at every point of
  $\mathcal M$, with $J$ being the 
  multiple principal hermitean structure on $\tn$.
\end{theorem} 
\begin{proof}
The proof of this Theorem consists of straightforward calculations
using the above definitions. The key point in these calculations is
that $\nw_{X}\aw(Y,X,Y)-\nw_{Y}\aw(X,X,Y)$, when $X,Y$ run through all
the vectors from $\tn$, is always proportional
to $\nw_4\aw_{141}-\nw_1\aw_{441}$. Here the indices $1$ an $4$ are
the components from the coframe adapted to $J$, in which $e_1=m$ and
$e_4=k$. By a direct calculation one can check that 
$\nw_4\aw_{141}-\nw_1\aw_{441}=16\Psi_1^2$. Thus, when 
$\nw_{X}\aw(Y,X,Y)\equiv\nw_{Y}\aw(X,X,Y)$, as assumed, $\Psi_1\equiv
0$, which proves the theorem.  
\end{proof}
\begin{remark}
When calculating $\nw_4\aw_{141}-\nw_1\aw_{441}$, during the proof of
the above theorem, we observed that the relation
$\nw_4\aw_{141}-\nw_1\aw_{441}=16\Psi_1^2$ is true even without the
(pseudo)hermitean reality conditions ($E$) or ($S_c$). For this
crucial relation to be true, we need to take $B$ as in (\ref{ebe}) and
to assume the integrability of $\tn$, i.e. to assume
$\kappa\equiv\sigma\equiv 0$. If these two assumptions are satisfied
then $\nw_4\aw_{141}-\nw_1\aw_{441}=16\Psi_1^2$ irrespective of the
signature of the metric. It is even true when the metric is complex!
Thus the Weyl connection $\nw$ with $B$ as in (\ref{ebe}) seems to
be meaningful in case of $g$ being complex, or having any
signature. The only trouble with such a connection is that in the
Lorentzian case it is complex. If one can live with this, one can
replace the condition (\ref{ghjj}) in Theorem \ref{ghj} by (\ref{scs})
and Theorem \ref{ghj} will be true for complex metrics, as well for
metrics of all the other real signatures.
\end{remark}
\section{Appendix}
The 36 signature independent Newman-Penrose equations, which include 16 first Bianchi identities, are:
\begin{eqnarray}
&& \bel\kappa = \del \si+
 \al' \kappa + 3 \bet \kappa  + \kappa \pi'  - 3 \eps \si + \eps' \si + \rho \si + 
  \rho' \si + \kappa \tau+ \Psi_0\label{np1}\end{eqnarray}\begin{eqnarray*}
&& \pel\kappa' = \del \si'+
 \al \kappa' + 3 \bet' \kappa'  + \kappa' \pi   - 3 \eps' \si'+ \eps \si' + \rho' \si'+ 
  \rho \si'  + \kappa' \tau'+ \Psi_0'\end{eqnarray*}
\begin{eqnarray}
&& \del \bet = \bel\eps -\al' \eps - \bet \eps' - \gam \kappa  - \kappa \mu - \eps \pi' - \bet \rho' - \al \si + \pi \si - 
  \Psi_1\label{np2}
\end{eqnarray}
\begin{eqnarray*}
&& \del \bet' =  \pel\eps'- \al \eps'-\bet' \eps  - \gam' \kappa' - \kappa' \mu' - \eps' \pi - \bet' \rho - \al' \si' + 
  \pi' \si' - \Psi_1'
\end{eqnarray*}
\begin{eqnarray}
&& \bel\rho = 
 \pel\si + \kappa \mu'- \kappa \mu  + \al' \rho + \bet \rho  - 3 \al \si + 
  \bet' \si - \rho' \tau+ \rho \tau  - \Psi_1- \Rho_{14}\label{np3}
\end{eqnarray}
\begin{eqnarray*}
&& \pel\rho' = \bel\si'+
 \kappa' \mu - \kappa' \mu' + \al \rho' + \bet' \rho' - 
  3 \al' \si' + \bet \si' - \rho \tau' + \rho' \tau' - \Psi_1' -
  \Rho_{24} 
\end{eqnarray*}
\begin{eqnarray}
&& \del \tau = \vel \kappa - \gam' \kappa-3 \gam \kappa   + \pi' \rho
  + \pi \si - \si \tau' - \eps' \tau + \eps \tau- \rho \tau - \Psi_1+
  \Rho_{14} \label{np4}
\end{eqnarray}
\begin{eqnarray*}
&& \del \tau' =  \vel \kappa'-\gam \kappa' - 
  3 \gam' \kappa'+ \pi \rho' + \pi' \si' - \si' \tau - 
  \eps \tau' + \eps' \tau' - \rho' \tau' - \Psi_1' + \Rho_{24} 
\end{eqnarray*}
\begin{eqnarray*}
&& \vel \rho = 
 \pel\tau - \kappa \nu + \gam \rho + \gam' \rho - \mu' \rho  - 
  \la \si - \al \tau + \bet' \tau - \tau \tau'- \Psi_2 - \Rho_{12} -
  \Rho_{34}
\end{eqnarray*}
\begin{eqnarray*}
&& \vel \rho' = 
 \bel\tau' - \kappa' \nu'+ \gam' \rho' + \gam \rho'  - \mu \rho' -
   \la' \si' - \al' \tau' + \bet \tau' - \tau \tau' - \Psi_2' - \Rho_{12} - \Rho_{34}\end{eqnarray*}
\begin{eqnarray*}
&& \vel \al =  \pel\gam+
 \bet' \gam + \al \gam' - \bet \la - \al \mu' - \eps \nu  + \nu \rho - 
  \la \tau - \gam \tau'+ \Psi_3
\end{eqnarray*}
\begin{eqnarray*}
&& \vel \al' = \bel\gam'+ \bet \gam'+
 \al' \gam  - \bet' \la'  - \al' \mu - \eps' \nu' + 
  \nu' \rho' - \la' \tau' - \gam' \tau + \Psi_3'
\end{eqnarray*}
\begin{eqnarray*}
&& \vel \la = \pel\nu-3 \gam \la + \gam' \la  - \la \mu -
   \la \mu' + 3 \al \nu + \bet' \nu - \nu \pi - \nu \tau'- \Psi_4 \end{eqnarray*}
\begin{eqnarray*}
&& \vel \la' = \bel\nu' - 3 \gam' \la'+
 \gam \la' - \la' \mu'- \la' \mu   + 3 \al' \nu' + \bet \nu' -
   \nu' \pi'  - \nu' \tau- \Psi_4'\end{eqnarray*}
\begin{eqnarray*}
&&\del \la =\pel\pi -3 \eps \la + \eps' \la  - \kappa' \nu + \al \pi - \bet' \pi - \pi^2 - 
  \la \rho  - \mu \si'- \Rho_{22}
\end{eqnarray*}
\begin{eqnarray}
&& \del \la' = \bel\pi' - 3 \eps' \la'+
 \eps \la'  - \kappa \nu' + \al' \pi' - \bet \pi' - \pi'^2 - \la' \rho' - \mu' \si - 
  \Rho_{11}\label{np5}
\end{eqnarray}
\begin{eqnarray*}
&& \del \mu = \bel\pi -\eps \mu - \eps' \mu - \kappa \nu - 
  \al' \pi + \bet \pi - \pi \pi'- \mu \rho' - 
  \la \si - \Psi_2 - \Rho_{12} - \Rho_{34} 
\end{eqnarray*}
\begin{eqnarray*}
&& \del \mu' = 
 \pel\pi'- \eps' \mu' - \eps \mu'  - \kappa' \nu' - \al \pi' + \bet' \pi' - \pi \pi' - \mu' \rho - \la' \si' - 
  \Psi_2' - \Rho_{12} - \Rho_{34}
\end{eqnarray*}
\begin{eqnarray*}
&& \del \al = \pel\eps + 
  \al \eps'-2 \al \eps - \bet' \eps - \gam \kappa' - \kappa \la  - \eps \pi - \al \rho + \pi \rho  - 
  \bet \si'+ \Rho_{24}
\end{eqnarray*}
\begin{eqnarray}
&& \del \al' = \bel\eps'+
 \al' \eps - 2 \al' \eps' - \bet \eps' - \gam' \kappa - \kappa' \la'  - \eps' \pi'- \al' \rho' + \pi' \rho' - \bet' \si + 
  \Rho_{14} \label{np6}
\end{eqnarray}
\begin{eqnarray*}
&& \vel \bet = \bel\gam+
 \al' \gam + 2 \bet \gam - \bet \gam' - \al \la'  - \bet \mu - \eps \nu'+ \nu \si - \gam \tau - \mu \tau - 
  \Rho_{13} 
\end{eqnarray*}
\begin{eqnarray*}
&& \vel \bet' = \pel\gam' + \al \gam' + 
  2 \bet' \gam'-\bet' \gam - \al' \la  - \bet' \mu' - \eps' \nu  + \nu' \si' - 
  \gam' \tau' - \mu' \tau'- \Rho_{23}
\end{eqnarray*}
\begin{eqnarray}
&& \del \rho = \pel\kappa-3 \al \kappa - \bet' \kappa  - \kappa \pi + 
  \eps \rho + \eps' \rho - \rho^2 - \si \si' - 
  \kappa' \tau- \Rho_{44} \label{np9}
\end{eqnarray}
\begin{eqnarray}
&& \del \rho' =  \bel\kappa'-3 \al' \kappa' - \bet \kappa'  - \kappa' \pi' + \eps' \rho' + 
  \eps \rho' - \rho'^2 - \si \si' - \kappa \tau' - \Rho_{44}\label{np8}
\end{eqnarray}
\begin{eqnarray*}
&& \vel \mu =  \bel\nu -\la \la' - \gam \mu - 
  \gam' \mu - \mu^2+ \al' \nu + 3 \bet \nu - \nu' \pi - 
  \nu \tau- \Rho_{33} \end{eqnarray*}
\begin{eqnarray*}
&& \vel \mu' = \pel\nu'-\la \la'  - \gam' \mu' - \gam \mu' - {\mu'}^2 + 
  \al \nu' + 3 \bet' \nu' - \nu \pi'- \nu' \tau' - \Rho_{33} \end{eqnarray*}
\begin{eqnarray*}
&& \del \nu = 
 \vel \pi - \eps' \nu - 3 \eps \nu + \la \pi' - \gam' \pi+ \gam \pi + \mu \pi- \mu \tau' - \la \tau  + \Psi_3 - 
  \Rho_{23} \end{eqnarray*}
\begin{eqnarray*}
&& \del \nu' = 
 \vel \pi' - \eps \nu' - 3 \eps' \nu' + \la' \pi - \gam \pi' + \gam' \pi' + \mu' \pi'- \mu' \tau - \la' \tau' + 
  \Psi_3' - \Rho_{13} \end{eqnarray*}
\begin{eqnarray*}&& 
\del \gam = 
  \vel \eps-2 \eps \gam - \eps' \gam - \eps \gam' - \kappa \nu + \bet \pi + \al \pi' - \al \tau + \pi \tau - 
  \bet \tau' - \Psi_2 + \Rho_{34}\end{eqnarray*}\begin{eqnarray*}
&& \del \gam' = \vel \eps'- 2 \eps' \gam'- \eps \gam' -\eps' \gam   - \kappa' \nu'  + \bet' \pi'+
   \al' \pi - \al' \tau'+ 
  \pi' \tau' - \bet' \tau  - \Psi_2' + \Rho_{34}\end{eqnarray*}\begin{eqnarray*}
&& \pel\mu = \bel\la-\al' \la + 3 \bet \la - \al \mu - \bet' \mu  + \mu \pi -
   \mu' \pi - \nu \rho  + \nu \rho'- \Psi_3 - \Rho_{23}\end{eqnarray*}\begin{eqnarray*}
&& \bel\mu' = \pel\la'-\al \la' + 
  3 \bet' \la' - \al' \mu' - \bet \mu' + \mu' \pi'- \mu \pi'- \nu' \rho'   + 
  \nu' \rho - \Psi_3'- \Rho_{13} \end{eqnarray*}\begin{eqnarray}
&& \bel\tau = 
 \vel \si + \kappa \nu' + \la' \rho  - 3 \gam \si + \gam' \si + \mu \si - \al' \tau +
   \bet \tau + \tau^2+ \Rho_{11}\label{np7}\end{eqnarray}\begin{eqnarray*}
&& \pel\tau' = 
 \vel \si' + \kappa' \nu  + \la \rho'  - 3 \gam' \si'+ \gam \si' + \mu' \si' - 
  \al \tau' + \bet' \tau' + \tau'^2+ \Rho_{22}\end{eqnarray*}\begin{eqnarray*}
&& \bel\al = \pel\bet+
 \al \al' - 2 \al \bet + \bet \bet'  - \eps \mu + \eps \mu' + \gam \rho + 
  \mu \rho  - \gam \rho' - \la \si - \Psi_2+ \Rho_{12}\end{eqnarray*}\begin{eqnarray*}
&& \pel\al' = \bel\bet'+
 \al \al' - 2 \al' \bet' + \bet \bet' - \eps' \mu' + \eps' \mu  + \gam' \rho' + \mu' \rho' - 
  \gam' \rho  - \la' \si' - \Psi_2'+ \Rho_{12}
\end{eqnarray*}

The 20 second Bianchi identities are:
\begin{eqnarray}
\bel\Psi_1 &=&\vel\Psi_0  -\del\Rho_{11} + \bel\Rho_{14} - 4 \gam \Psi_0 + \mu \Psi_0 + 2 \bet \Psi_1 - 3\si \Psi_2  + 4 \tau\Psi_1  -\nonumber\\&& 2 \kappa \Rho_{13}  + 
   2 \eps \Rho_{11} - 2 \eps' \Rho_{11} - 2 \bet \Rho_{14} - 2 \pi' \Rho_{14} + 
  \la' \Rho_{44} - \rho'\Rho_{11} - \si\Rho_{12} + \si\Rho_{34}\label{b1} \end{eqnarray}\begin{eqnarray*}
\pel\Psi_1' &=& \vel\Psi_0' -\del\Rho_{22}  + \pel\Rho_{24}- 4 \gam' \Psi_0'+ \mu' \Psi_0' + 
   2 \bet' \Psi_1' - 3 \si'\Psi_2' + 4 \tau'\Psi_1'  -\\&& 2 \kappa' \Rho_{23} + 2 \eps' \Rho_{22}   - 2 \eps \Rho_{22}- 
   2 \bet' \Rho_{24} - 2 \pi \Rho_{24} + \la \Rho_{44}- \rho \Rho_{22}  -\si' \Rho_{12}  + \si'
   \Rho_{34}  \end{eqnarray*}\begin{eqnarray}
\del\Psi_1& =& - \pel\Psi_0 -\del\Rho_{14} + \bel\Rho_{44} + 4 \al \Psi_0 + \pi \Psi_0 + 2 \eps \Psi_1 - 
   3 \kappa \Psi_2 - 4 \Psi_1 \rho +\kappa' \Rho_{11}  + \nonumber\\&&\kappa \Rho_{12} + 2 \eps \Rho_{14} - \kappa \Rho_{34} - 
   2 \al' \Rho_{44} - 2 \bet \Rho_{44} - \pi' \Rho_{44} - 2 \rho'\Rho_{14}  - 2 \si\Rho_{24} \label{b2}\end{eqnarray}\begin{eqnarray*}
\del\Psi_1' &=&  - \bel\Psi_0'-\del\Rho_{24} + \pel\Rho_{44} + 4 \al' \Psi_0' + \pi' \Psi_0' + 
   2 \eps' \Psi_1' - 3 \kappa' \Psi_2' - 
   4 \Psi_1' \rho' + \kappa \Rho_{22} +\\&& \kappa' \Rho_{12} + 2 \eps' \Rho_{24} - \kappa' \Rho_{34} - 2 \al \Rho_{44} - 2 \bet' \Rho_{44} - \pi \Rho_{44}  - 
   2 \rho \Rho_{24}- 2 \si'\Rho_{14} \end{eqnarray*}\begin{eqnarray}
\vel\Psi_1 &=&   \bel\Psi_2+\del\Rho_{13} - \bel\Rho_{34} + \nu \Psi_0 + 2 \gam \Psi_1 - 2 \mu \Psi_1- 2\si \Psi_3 - 3 \tau\Psi_2  - 
  \nonumber \\&& \pi \Rho_{11} - \pi' \Rho_{12} +2 \eps' \Rho_{13} + \mu \Rho_{14}
   + \la' \Rho_{24} + \kappa \Rho_{33} + 
   \pi' \Rho_{34} +\rho' \Rho_{13}   +\si \Rho_{23} \label{b3}\end{eqnarray}\begin{eqnarray*}
\vel\Psi_1' &=& 
  \pel\Psi_2'+ \del\Rho_{23} - \pel\Rho_{34} + \nu' \Psi_0' + 2 \gam' \Psi_1' - 2 \mu' \Psi_1'- 2\si' \Psi_3'  - 3 \tau'\Psi_2'  - 
   \\&& \pi' \Rho_{22}- \pi \Rho_{12}+ 2 \eps \Rho_{23}+ \mu' \Rho_{24}  + \la \Rho_{14}  + 
   \kappa' \Rho_{33} + \pi \Rho_{34} + \rho \Rho_{23} + \si'\Rho_{13} \end{eqnarray*}\begin{eqnarray}
\pel\Psi_1& = &-\del\Psi_2 + \del\Rho_{12} - \bel\Rho_{24} + \la \Psi_0 + 2 \al \Psi_1 + 2 \pi \Psi_1 + 
   2 \kappa \Psi_3 - 3 \rho\Psi_2  + \nonumber\\&&\kappa' \Rho_{13} + \pi \Rho_{14} + \kappa \Rho_{23} + 2 \al' \Rho_{24} +
    \pi' \Rho_{24} - \mu \Rho_{44} + \rho'\Rho_{12}  - \rho'\Rho_{34}
    + 
\si\Rho_{22}\label{b4}\end{eqnarray}\begin{eqnarray*}
\bel\Psi_1' &=& -\del\Psi_2' + \del\Rho_{12} - \pel\Rho_{14} + \la' \Psi_0' + 2 \al' \Psi_1' + 2 \pi' \Psi_1' + 2 \kappa' \Psi_3' - 3 \rho'\Psi_2'  +  \\&&\kappa \Rho_{23}+\pi' \Rho_{24} + \kappa' \Rho_{13}  +  2 \al \Rho_{14} + 
   \pi \Rho_{14} - \mu' \Rho_{44}+ \rho\Rho_{12}  - \rho \Rho_{34} + 
   \si'\Rho_{11} \end{eqnarray*}\begin{eqnarray}
\vel\Psi_2 &=&  - \bel\Psi_3+\vel\Rho_{12} - \pel\Rho_{13} + 2 \nu \Psi_1 - 3 \mu \Psi_2 - 2 \bet \Psi_3 + 2 \tau\Psi_3  + \si\Psi_4  + 
   \nonumber\\&& \la \Rho_{11} + \mu' \Rho_{12} - 2 \bet' \Rho_{13} +\nu \Rho_{14}
   + \nu' \Rho_{24} - \rho 
\Rho_{33} - 
   \mu' \Rho_{34} + \tau\Rho_{23} + \tau'\Rho_{13} \label{b5}\end{eqnarray}\begin{eqnarray*}
\vel\Psi_2'& = &
  - \pel\Psi_3'+ \vel\Rho_{12} - \bel\Rho_{23} + 2 \nu' \Psi_1' - 3 \mu' \Psi_2' - 2 \bet' \Psi_3' + 2\tau' \Psi_3' +\si' \Psi_4'  + 
    \\&&  \la' \Rho_{22}+ \mu \Rho_{12}-2 \bet \Rho_{23} +
   \nu' \Rho_{24} + \nu \Rho_{14}- 
 \rho'  \Rho_{33} - 
\mu \Rho_{34}  + \tau'\Rho_{13} + \tau\Rho_{23} \end{eqnarray*}
\begin{eqnarray}
\del\Psi_3 &= & \pel\Psi_2+\vel\Rho_{24} - \pel\Rho_{34} - 2 \la \Psi_1 - 3 \pi \Psi_2 - 2 \eps \Psi_3 + 
   \kappa \Psi_4 - 2  \rho\Psi_3 + \la \Rho_{14} +\nonumber\\&& \rho \Rho_{23} - 2 \gam' \Rho_{24} + \mu' \Rho_{24} + 
   \nu \Rho_{44} + \si' \Rho_{13} - \tau\Rho_{22}  -\tau' \Rho_{12}
   +\tau' \Rho_{34} \label{b6}\end{eqnarray}
\begin{eqnarray*}
\del\Psi_3' &=& 
  \bel\Psi_2'+ \vel\Rho_{14} - \bel\Rho_{34} - 2 \la' \Psi_1' - 3 \pi' \Psi_2' - 2 \eps' \Psi_3' + 
   \kappa' \Psi_4' - 
   2\rho' \Psi_3' +\la' \Rho_{24} +\\&& \rho' \Rho_{13} - 2 \gam \Rho_{14} + \mu \Rho_{14} + \nu' \Rho_{44}  +  \si \Rho_{23} -\tau' \Rho_{11} -  \tau\Rho_{12} +  \tau\Rho_{34}\end{eqnarray*}
\begin{eqnarray}
\vel\Psi_3 &=& - \bel\Psi_4  -\vel\Rho_{23} + \pel\Rho_{33}- 3 \nu \Psi_2 - 2 \gam \Psi_3 - 
   4 \mu \Psi_3 - 4 \bet \Psi_4+ 
   \tau\Psi_4  + \nu \Rho_{12} - \nonumber\\&& 2 \la \Rho_{13} +\nu' \Rho_{22} - 
   2 \gam \Rho_{23} - 2 \mu' \Rho_{23} + 2 \al \Rho_{33} + 2 \bet'
   \Rho_{33} - \nu \Rho_{34}  - \tau' \Rho_{33}\label{b7}\end{eqnarray}
\begin{eqnarray*}
\vel\Psi_3' &=&- \pel\Psi_4' -\vel\Rho_{13}  + \bel\Rho_{33} - 3 \nu' \Psi_2' - 2 \gam' \Psi_3' - 
   4 \mu' \Psi_3' - 4 \bet' \Psi_4'+ \tau'\Psi_4'  + \nu' \Rho_{12} - \\&&  2 \la' \Rho_{23}+\nu \Rho_{11}-2 \gam' \Rho_{13} - 
   2 \mu \Rho_{13} + 2 \al' \Rho_{33} + 2 \bet \Rho_{33} - \nu' \Rho_{34} - 
  \tau \Rho_{33}  \end{eqnarray*}
\begin{eqnarray}
\pel\Psi_3& =& 
  \del\Psi_4 - \vel\Rho_{22} + \pel\Rho_{23} - 3 \la \Psi_2 - 2 \al \Psi_3 + 4 \pi \Psi_3 + 
   4 \eps \Psi_4 + \rho\Psi_4  + \nonumber\\&& 2 \gam' \Rho_{22} -2 \gam \Rho_{22}  - 
   \mu' \Rho_{22}- \la \Rho_{12} + 2 \al \Rho_{23} - 2 \nu \Rho_{24} + \la \Rho_{34} + \si'\Rho_{33}  - 
   2 \tau'\Rho_{23}\label{b8}\end{eqnarray}
\begin{eqnarray*}
\bel\Psi_3'& = &
  \del\Psi_4' - \vel\Rho_{11} + \bel\Rho_{13} - 3 \la' \Psi_2' - 2 \al' \Psi_3' + 4 \pi' \Psi_3' + 
   4 \eps' \Psi_4' +\rho' \Psi_4' + \\&& 2 \gam \Rho_{11} - 2 \gam' \Rho_{11} -\mu \Rho_{11} - \la' \Rho_{12} + 
   2 \al' \Rho_{13} - 2 \nu' \Rho_{14} + \la' \Rho_{34}  +\si \Rho_{33}  - 
   2\tau \Rho_{13} \end{eqnarray*}
\begin{eqnarray}
\bel\Rho_{12}& =& \del\Rho_{13} + \vel\Rho_{14} + \pel\Rho_{11} - 2 \bel\Rho_{34} - 2 \al \Rho_{11} + 2 \bet' \Rho_{11} - 
   \pi \Rho_{11} - \pi' \Rho_{12} + 2 \eps' \Rho_{13} + \nonumber\\&&2 \rho \Rho_{13} - 2 \gam \Rho_{14} + 
   \mu \Rho_{14} + 2 \mu' \Rho_{14} + \la' \Rho_{24} + \kappa \Rho_{33} + \pi' \Rho_{34} + \nu' \Rho_{44} + 
  \rho' \Rho_{13}  +\label{b9}\\&& \si \Rho_{23} -  \tau\Rho_{12} +
  \tau\Rho_{34}  -  
\tau'\Rho_{11}\nonumber\end{eqnarray}\begin{eqnarray}
\pel\Rho_{12}& =& 
  \del\Rho_{23} + \vel\Rho_{24} + \bel\Rho_{22} - 2 \pel\Rho_{34} - 2 \al' \Rho_{22} + 
   2 \bet \Rho_{22} - \pi' \Rho_{22}- \pi \Rho_{12} + 2 \eps \Rho_{23} +\nonumber\\&&  2 \rho'\Rho_{23}   - 2 \gam' \Rho_{24} + \mu' \Rho_{24} + 
   2 \mu \Rho_{24} +\la \Rho_{14}  + \kappa' \Rho_{33} + \pi \Rho_{34} + \nu \Rho_{44} +\rho \Rho_{23}+\nonumber\\&& 
   \si' \Rho_{13}  - \tau'\Rho_{12}  +  \tau'\Rho_{34}-  \tau\Rho_{22}\nonumber
\end{eqnarray}
\begin{eqnarray}
\del\Rho_{34} &= &-2 \del\Rho_{12} + \vel\Rho_{44} + \pel\Rho_{14} + \bel\Rho_{24} - \rho \Rho_{12} - \kappa' \Rho_{13} - 
   2 \al \Rho_{14} - \pi \Rho_{14} - \kappa \Rho_{23} - \nonumber\\&&2 \al' \Rho_{24} - \pi' \Rho_{24} + \rho \Rho_{34} - 
   2 \gam \Rho_{44} - 2 \gam' \Rho_{44} + \mu \Rho_{44} + \mu' \Rho_{44} -  \rho'\Rho_{12} + 
   \rho' \Rho_{34} - \label{b10}\\&&\si \Rho_{22} - \si' \Rho_{11} - 2 \tau \Rho_{24}
   - 2 \tau' \Rho_{14}\nonumber
\end{eqnarray}
\begin{eqnarray}
\vel\Rho_{34}& =& \del\Rho_{33} - 2 \vel\Rho_{12} + \pel\Rho_{13} + \bel\Rho_{23} - \la \Rho_{11} - \mu \Rho_{12} - 
   \mu' \Rho_{12} + 2 \bet' \Rho_{13} - 2 \pi \Rho_{13} -\nonumber\\&& \nu \Rho_{14} - \la' \Rho_{22} + 
   2 \bet \Rho_{23} - 2 \pi' \Rho_{23} - \nu' \Rho_{24} + 2 \eps \Rho_{33} + 2 \eps' \Rho_{33} + 
   \rho \Rho_{33} + \mu \Rho_{34} + \nonumber\\
&&\mu' \Rho_{34} +  \rho'\Rho_{33} -  \tau\Rho_{23} -  
\tau'\Rho_{13}\nonumber\end{eqnarray}
Using relations (\ref{420})-(\ref{421}) we can reexpress identities
(\ref{b1})-(\ref{b8}) in terms of the components of the Cotton
tensor. After this the Cotton tensor components `hide' the terms with
the Schouten tensor components $\Rho_{ij}$, and the respective
identities assume a more compact form as follows:
\begin{eqnarray}
A_{141}&=&\vel\Psi_0 +(\mu- 4 \gam) \Psi_0 -\bel\Psi_1 + 2 (2\tau+\bet) \Psi_1 - 3\si \Psi_2 \label{b11} \end{eqnarray}
\begin{eqnarray}
A_{414}&=& \pel\Psi_0 -(\pi+ 4 \al) \Psi_0  +\del\Psi_1+ 2 (2\rho-\eps) \Psi_1 + 
   3 \kappa \Psi_2 \label{b21}\end{eqnarray}
\begin{eqnarray}
A_{341}&=&\vel\Psi_1+ 2 (\mu-\gam) \Psi_1-  \bel\Psi_2 + 3 \tau\Psi_2- \nu \Psi_0+ 2\si \Psi_3  \label{b31}\end{eqnarray}
\begin{eqnarray}
A_{214}& = &\pel\Psi_1 - 2 (\al+\pi) \Psi_1 +\del\Psi_2 + 3 \rho\Psi_2  - \la \Psi_0 - 
   2 \kappa \Psi_3 \label{b41}\end{eqnarray}
\begin{eqnarray}
A_{132} &=&\vel\Psi_2  +3 \mu \Psi_2+ \bel\Psi_3  + 2 (\bet-\tau) \Psi_3 - 2 \nu \Psi_1- \si\Psi_4  \label{b51}\end{eqnarray}
\begin{eqnarray}
A_{423}&=&  \pel\Psi_2  - 3 \pi \Psi_2-\del\Psi_3 - 2 (\eps+\rho) \Psi_3 - 2 \la \Psi_1+ 
   \kappa \Psi_4  \label{b61}\end{eqnarray}
\begin{eqnarray}
A_{323} &=&\vel\Psi_3 + 2 (\gam+2\mu) \Psi_3 + \bel\Psi_4   +(4 \bet-\tau) \Psi_4+ 3 \nu \Psi_2 \label{b71}\end{eqnarray}
\begin{eqnarray}
A_{223}& =&\pel\Psi_3 + 2 (\al-2\pi) \Psi_3   
  -\del\Psi_4- 
   (\rho+4 \eps )\Psi_4   + 3 \la \Psi_2, \label{b81}\end{eqnarray}
with the analogous identities for the primed quantities.

\end{document}